\documentclass[a4paper,10pt]{article}

\usepackage{latexsym}
\usepackage{mathrsfs}
\usepackage{amsfonts,amsmath,amssymb}
\usepackage{indentfirst}
\usepackage{subeqnarray}
\usepackage{verbatim}
\usepackage[pdftex]{graphicx}
\usepackage{epstopdf}
\usepackage{stmaryrd}
\usepackage{multirow}
\usepackage{diagbox}
\usepackage{subcaption}
\usepackage{pdfpages}
\usepackage{bm}

\newcommand{\bx}{\bm{x}}
\newcommand{\bX}{\bm{X}}
\newcommand{\bz}{\bm{z}}
\newcommand{\bpsi}{\bm{\psi}}
\newcommand{\bphi}{\bm{\phi}}

\newcommand{\bb}{\bm{b}}
\newcommand{\bB}{\bm{B}}
\newcommand{\bh}{\bm{h}}
\newcommand{\bH}{\bm{H}}
\newcommand{\bmm}{\bm{m}}
\newcommand{\br}{\bm{r}}

\newcommand{\bv}{\bm{v}}
\newcommand{\bw}{\bm{w}}
\newcommand{\bW}{\bm{W}}
\newcommand{\bp}{\bm{p}}
\newcommand{\bq}{\bm{q}}
\newcommand{\by}{\bm{y}}
\newcommand{\bcW}{\mathbfcal{W}}
\newcommand{\bcB}{\mathbfcal{B}}
\newcommand{\cB}{\mathcal{B}}
\newcommand{\cU}{\mathcal{U}}

\newcommand{\cA}{\mathcal{A}}
\newcommand{\cG}{\mathcal{G}}
\newcommand{\cL}{\mathcal{L}}
\newcommand{\cN}{\mathcal{N}}
\newcommand{\dE}{\mathbb{E}}
\newcommand{\dP}{\mathbb{P}}
\newcommand{\balpha}{\bm{\alpha}}
\newcommand{\bbeta}{\bm{\beta}}
\newcommand{\bgamma}{\bm{\gamma}}
\newcommand{\bxi}{\bm{\xi}}
\newcommand{\dbx}{\mathrm{d}\bx}
\newcommand{\ds}{\mathrm{d}s}
\newcommand{\dy}{\mathrm{d}y}
\newcommand{\Real}{\mathbb{R}}
\newcommand{\set}[1]{\left\{#1\right\}}

\newtheorem{theorem}{Theorem}[section]
\newtheorem{lemma}[theorem]{Lemma}
\newtheorem{assumption}[theorem]{Assumption}
\newtheorem{example}{Example}
\newtheorem{remark}[theorem]{Remark}
\numberwithin{equation}{section}
\DeclareMathAlphabet\mathbfcal{OMS}{cmsy}{b}{n}

\newenvironment{proof}[1][Proof]{\textbf{#1.} }
{\ \rule{0.75em}{0.75em}\smallskip}

\usepackage[pdftex,unicode,colorlinks,linkcolor=blue,citecolor=red]{hyperref}

\begin{document}

\begin{center}
\Large\bf \textbf{Adaptive-Growth Randomized Neural Networks for PDEs: Algorithms and Numerical Analysis}
\end{center}

\begin{center}
    {\large\sc Haoning Dang}\footnote{School of Mathematics and Statistics, Xi'an Jiaotong University, Xi'an, Shaanxi 710049, China. E-mail: {\tt haoningdang.xjtu@stu.xjtu.edu.cn}},\quad
    {\large\sc Fei Wang}\footnote{School of Mathematics and Statistics \& State Key Laboratory of Multiphase Flow in Power Engineering, Xi'an Jiaotong University, Xi'an, Shaanxi 710049, China. The work of this author was partially supported by the National Natural Science Foundation of China (Grant No.\ 92470115) and Shaanxi Fundamental Science Research Project for Mathematics and Physics (Grant No.\ 22JSY027). Email: {\tt feiwang.xjtu@xjtu.edu.cn}}, \quad
    {\large\sc Song Jiang}\footnote{Institute of Applied Physics and Computational Mathematics, Beijing 100088, China. E-mail: {\tt jiang@iapcm.ac.cn}}  
\end{center}

\medskip
\begin{quote}

{\bf Abstract. } 
Randomized neural network (RaNN) methods have been proposed for solving various partial differential equations (PDEs), demonstrating high accuracy and efficiency. However, initializing the fixed parameters remains challenging. Additionally, RaNNs often struggle to approximate PDE solutions with sharp gradients or discontinuities when using smooth activations and shallow architectures. In this paper, we propose an Adaptive-Growth Randomized Neural Network (AG-RaNN) to address these challenges. We first design a frequency-based initialization for a shallow RaNN. Using the residual as an error indicator, we then adaptively grow the network in width (neuron growth) and depth (layer growth) to improve the accuracy of the numerical solution. The weights and biases of new neurons are constructed rather than trained, which enhances the approximation power without additional nonlinear optimization. To handle discontinuities, we further introduce a domain splitting strategy. We also establish a unified error analysis covering approximation, statistical, and optimization errors. Extensive numerical experiments demonstrate the efficiency and accuracy of AG-RaNN.

\end{quote}

\medskip


\section{Introduction}

Partial differential equations (PDEs) are widely used in fields such as physics, engineering, economics, and biology. Traditional numerical methods, such as finite element, finite difference, finite volume, and spectral methods, have been well-established in solving PDEs. For oscillatory or discontinuous solutions, discontinuous Galerkin (DG) methods are often employed. While these methods are effective, they rely on mesh generation, which can be time-consuming, especially in complex geometries. Meshes may also distort under large deformations or moving interfaces, requiring remeshing. Moreover, traditional methods suffer from the ``curse of dimensionality," making them less efficient in high-dimensional problems. Even in relatively low dimensions (e.g., $d \leq 3$), solving equations such as the Helmholtz equation with high wave numbers can be computationally expensive due to the fine discretization required to capture oscillations. These challenges motivate the search for alternative, mesh-free approximation frameworks for PDEs.

Recently, deep neural networks have been proposed to solve PDEs, such as the deep Ritz method (\cite{E2018DRM}), physics-informed neural networks (\cite{Raissi2019Physics}), and the deep Galerkin method (\cite{Sirignano2018DGM}), among others (\cite{Xu2020FiniteNeuron,Liao2021DNM}). These approaches are mesh-free and therefore particularly suited for PDEs on complex geometries, where generating high-quality meshes is difficult or expensive. Their strong approximation capability also offers the potential to alleviate the curse of dimensionality. However, training neural networks involves solving nonlinear and nonconvex optimization problems, which can result in large optimization errors and lower accuracy, especially for stiff or high-frequency PDEs.

To reduce optimization errors, Randomized Neural Networks (RaNNs, \cite{Pao1994Learning, Igelnik1995Stochastic}) offer a simpler approach by fixing some parameters and training only the linear coefficients in the output layer. This turns the problem into a least-squares solution, which is both efficient and accurate. Extreme Learning Machines (ELM, \cite{Huang2006Extreme}), a type of RaNN, have been used to solve PDEs (\cite{Yang2018ELM, Fabiani2021ELM, Wang2024ELM}), along with methods such as local ELM (\cite{Dong2021Local, Dong2023Method}), the random feature method (\cite{Chen2022RFM}) and local RaNN (LRaNN) with finite difference method (\cite{Li2023LRNN}), which couple different neural networks in subdomains. RaNNs have also been effective when combined with weak form methods such as Petrov-Galerkin and DG to solve various equations, such as RaNN Petrov-Galerkin method (\cite{Shang2023Randomized, Shang2023RNNPG}), local RaNN-DG (LRaNN-DG) methods (\cite{Sun2022LRNNDG,Sun2023DVWE,Sun2025KdV}), LRaNN hybridized discontinuous Petrov-Galerkin method (\cite{Dang2024HDPG}). These developments show that RaNN-based models can significantly reduce optimization difficulty while retaining strong approximation capability.

Despite these advantages, generating informative, task-relevant fixed parameters for RaNNs remains a major challenge. Some studies employ optimization techniques, such as differential evolution or Gaussian random fields, to address this (\cite{Dong2022Parameters,Zhang2024Transferable}). Growing neural networks (\cite{Zemouri2020}), in which layers and neurons are added based on error estimators, have also been proposed, and adaptive-growth ELMs (\cite{Lan2010CELM,Zhang2012GELM}) provide constructive strategies for generating new neurons. Nevertheless, these approaches are typically limited to single hidden layers and do not systematically exploit depth, which may be insufficient for solving more complex PDEs with sharp fronts, boundary/internal layers, or discontinuities.

In this study, we propose a novel Adaptive-Growth Randomized Neural Network (AG-RaNN) framework that addresses the above challenges at both the parameter and architecture levels. The main structure of AG-RaNN is summarized as follows. First, we propose a frequency-based parameter initialization technique to generate a one-hidden-layer RaNN, using discrete Fourier transforms to construct the initial parameters. We then introduce neuron growth, where additional neurons are appended to the existing layer using similar frequency-based techniques. Since the RaNN methodology reduces optimization error at the cost of approximation flexibility, it can struggle to approximate sharp or highly localized solutions. To address this, we introduce a deeper network structure through layer growth: new hidden layers are added to the existing network, with weights and biases determined by residual error estimators. Finally, to tackle discontinuous solutions, we propose a domain-splitting strategy that partitions the domain according to the range of the approximate solution, allowing several neural networks to approximate different parts of the solution and thereby enhancing the overall approximation power. This approach is particularly suited for solving PDEs with sharp or discontinuous solutions.

From the viewpoint of approximation theory, several works have studied randomized neural networks. Classical results for RaNN-type models (\cite{Igelnik1995Stochastic}) establish existential universal approximation theorems in the $L^2$ norm with rates of order $1/2$. These results define families of sampling distributions (e.g., uniform distributions on bounded domains for weights and biases), yet the key hyperparameters—such as scaling factors and distribution parameters—are selected only through existence arguments depending on unknown properties of the target function. A similar situation arises in random Fourier features (\cite{Rahimi2007RFF}), as well as in ELM and related random-basis constructions (\cite{Huang2006Extreme,Lan2010CELM,Zhang2012GELM}). In all these cases, the theory is ``constructive'' only at the level of specifying the distributional form, but does not provide an explicit algorithm for choosing the hyperparameters or adapting the sampling distribution. To the best of our knowledge, no existing RaNN framework offers a fully algorithmic and residual-driven mechanism for updating both the network architecture and the sampling distribution of random parameters.

In this work, we pursue a more constructive approach. We categorize the theoretical frameworks into three conceptual tiers: the first refers to purely existential universal approximation theorems, which merely guarantee the existence of an appropriate distribution; the second describes formulations that define an explicit family of sampling distributions, yet still rely on non-constructive hyperparameter selection; and the highest tier comprises fully algorithmic frameworks in which both the sampling distribution and network architecture are updated through explicit, residual-driven rules with convergence guarantees. Existing theories of randomized neural networks largely reside in the second category. In contrast, the proposed AG-RaNN is conceived as a fully algorithmic approach: it leverages frequency information and residual-based error indicators to adaptively expand the hidden layers and implicitly adjust the sampling distribution of random parameters.

\textbf{Main contributions.}
The main contributions of this paper are summarized as follows:
\begin{itemize}
  \item We propose the first AG-RaNN framework, which combines frequency-based initialization, neuron growth, layer growth, and domain splitting, providing a systematic way to construct deep, adaptive randomized neural networks for PDEs.
        
  \item We develop a unified error analysis for RaNN-type PDE solvers in an abstract Hilbert-space framework. By establishing a graph-norm equivalence between the PDE residual and the natural energy norm, we obtain a clean decomposition of the total error into approximation, statistical, and optimization components. To the best of our knowledge, such a graph-norm-based analysis has not previously been available for randomized neural network PDE methods.
  
  \item We demonstrate that AG-RaNN is particularly effective for PDEs with sharp or discontinuous solutions, achieving higher accuracy than baseline RaNN methods in a series of numerical experiments.
\end{itemize}

The remainder of the paper is organized as follows. Section~\ref{Sec:NN} revisits the basics of neural networks and presents shortcut connections between different layers. Section~\ref{Sec:strategy} introduces strategies to construct adaptive-growth RaNNs to solve PDEs. Section~\ref{sec:analysis} provides a comprehensive theoretical analysis of RaNN methods. Section~\ref{Sec:NumerExper} offers numerical experiments demonstrating the efficiency of the proposed method. Finally, Section~\ref{sec:sum} summarizes our findings and outlines directions for future research.

\section{Neural Networks} \label{Sec:NN}

In this section, we present basic knowledge of neural networks and introduce necessary notation.

\subsection{Fully Connected Feedforward Neural Networks}

A fully connected feedforward neural network (NN) consists of $L+2$ layers, where $L$ represents the number of hidden layers, and the $l$-th layer contains $m_l\,(l=0,\cdots,L+1)$ neurons. The activation function in the $l$-th layer denoted by $\rho_l(\bx)$ acts elementwise on vector $\bx\in\Real^{m_l}$, where $l=1,\cdots,L$. The weights between the $i$-th and $j$-th layers are $\bW^{(i,j)}\in\Real^{m_j\times m_i}$, where $i<j$ and $i,j=0,\cdots,L+1$. The biases in the $l$-th layer are denoted by $\bb^{(l)}\in\Real^{m_l\times 1}$, where $l=1,\cdots,L+1$. With input $\bx\in\Real^d$, the fully connected feedforward NN can be mathematically described as follows:
\begin{align} 
    \by^{(0)}&=\bx, \\
    \by^{(l)}&=\rho_l\left(\bW^{(l-1,l)}\by^{(l-1)}+\bb^{(l)}\right), \quad l=1,\cdots,L, \\
    \bz&=\bW^{(L,L+1)}\by^{(L)}+\bb^{(L+1)}, 
\end{align}
where $\by^{(l)}\,(l=1,\cdots,L)$ are the hidden layers, and $\bz$ is the network output. Throughout this paper we set $b^{(L+1)} = 0$ for simplicity.

\subsection{Shortcut Connections}

We extend the feedforward neural network by adding shortcut connections between layers, inspired by residual neural networks (\cite{He2016ResNet}). The network with shortcuts is described as:
\begin{align} 
    \by^{(1)}&=\rho_1\left(\bW^{(0,1)}\bx+\bb^{(1)}\right), \\
    \by^{(l)}&=\rho_l\left(\sum_{k=1}^{l-1}\,\bW^{(k,l)}\by^{(k)}+\bb^{(l)}\right), \quad l=2,\cdots,L, \\
    \bz&=\sum_{l=1}^L\,\bW^{(l,L+1)}\by^{(l)}.
\end{align}

\begin{figure}[!htbp]
    \centering
    \vspace{-0.7cm}
    \includegraphics[width=0.8\textwidth]{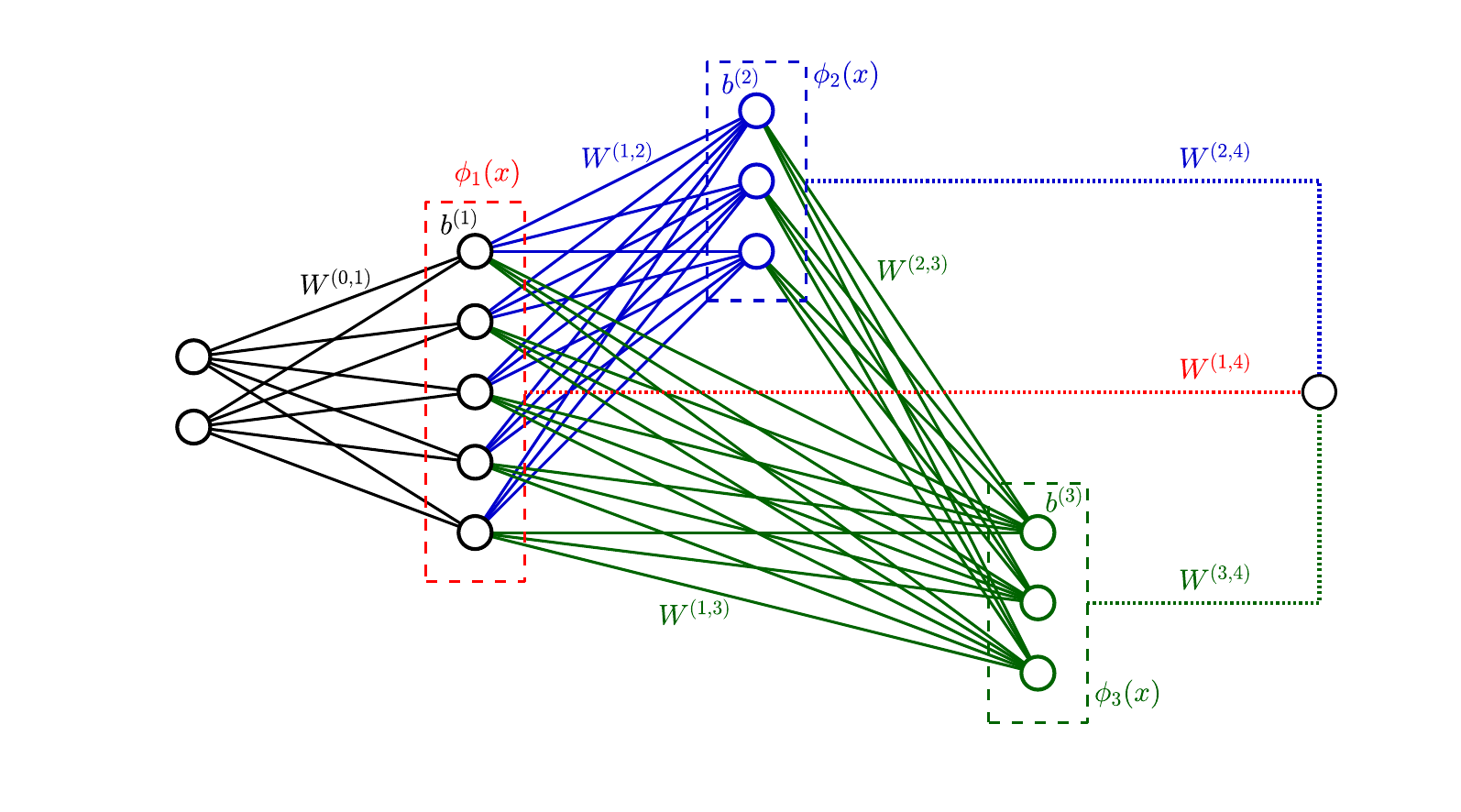}
    \vspace{-0.7cm}
    \caption{An RaNN with shortcut connections (solid lines: fixed; dotted lines: trainable). 
    Configuration: $(d,m_1,m_2,m_3,m_4)=(2,5,3,3,1)$. }
    \label{fig:RNN-shortcut}
\end{figure}

\subsection{Randomized Neural Networks}

In traditional neural networks, training involves solving a nonlinear and nonconvex problem using methods like stochastic gradient descent (SGD), which can be time-consuming and prone to large optimization errors. Randomized Neural Networks (RaNNs) address this by fixing most of the parameters and training only the final layer using least-squares methods, significantly reducing the optimization errors.

In RaNNs with shortcut connections, the trainable parameters are the weights $\set{\bW^{(l,L+1)}}_{l=1}^L$. We denote the fixed parameters as $\bcW$ and $\bcB$, respectively, where
\begin{align*}
\bcW\,=\,\set{\bW^{(0,1)},\bW^{(1,2)},\cdots,\bW^{(L-1,L)}}\quad\text{and}\quad\bcB\,=\,\set{\bb^{(1)},\bb^{(2)},\cdots,\bb^{(L)}}.
\end{align*}
Given a bounded domain $D\subset\Real^d$, we introduce basis functions $\bphi_l(\bx)=\by^{(l)}$ in the $l$-th layer, with each component representing a nonlinear function of $\bx$. The collective basis functions for an RaNN with shortcut connections are represented as $\set{\bphi_l(\bx)}_{l=1}^L$. We simplify $\bW^{(l,L+1)}$ to $\bW^{(l)}$. The function space of an RaNN is defined as:
\begin{align} \label{rnn_space}
    \cN^{\bcW,\bcB}_{\rho,m_{L+1}}(D)\,=\,\set{\bz(\bx)=\sum_{l=1}^L\,\bW^{(l)}\bphi_l(\bx)\,:\,\bx\in D\subset\Real^d, \, \forall\,\bW^{(l)}\in\Real^{m_{L+1}\times m_l}}.
\end{align}
The total number of trainable degrees of freedom (i.e., output-layer coefficients) is $M=\sum_{l=1}^Lm_l$. The network structure with shortcut connections is depicted in Figure \ref{fig:RNN-shortcut}, with solid lines representing fixed parameters and dotted lines representing trainable parameters.

In this paper, we focus on scalar-valued functions, and the RaNN function space is defined as:
\begin{align} \label{rnn_space1}
    \cN^{\bcW,\bcB}_{\rho,1}(D)\,=\,\set{z(\bx)=\sum_{l=1}^L\,\bw^{(l)}\bphi_l(\bx)\,:\,\bx\in D\subset\Real^d, \, \forall\,\bw^{(l)}\in\Real^{1\times m_l}}.
\end{align}

A fundamental challenge in the RaNN methodology lies in determining the fixed parameters, particularly in constructing appropriate basis functions to approximate the unknown target function. Igelnik and Pao (\cite{Igelnik1995Stochastic}) established the approximation property of RaNNs with a single hidden layer, proving that for any continuous function, there exist suitable activation functions and parameter distributions such that a shallow RaNN can approximate the target function with high probability. However, since the proof is nonconstructive, it does not provide a concrete algorithm for identifying an appropriate parameter distribution. Therefore, the next section introduces several constructive strategies for parameter determination and RaNN architecture design.

\section{Adaptive-Growth RaNN} \label{Sec:strategy}

\begin{figure}[!htbp]
    \centering
    \includegraphics[width=0.98\textwidth]{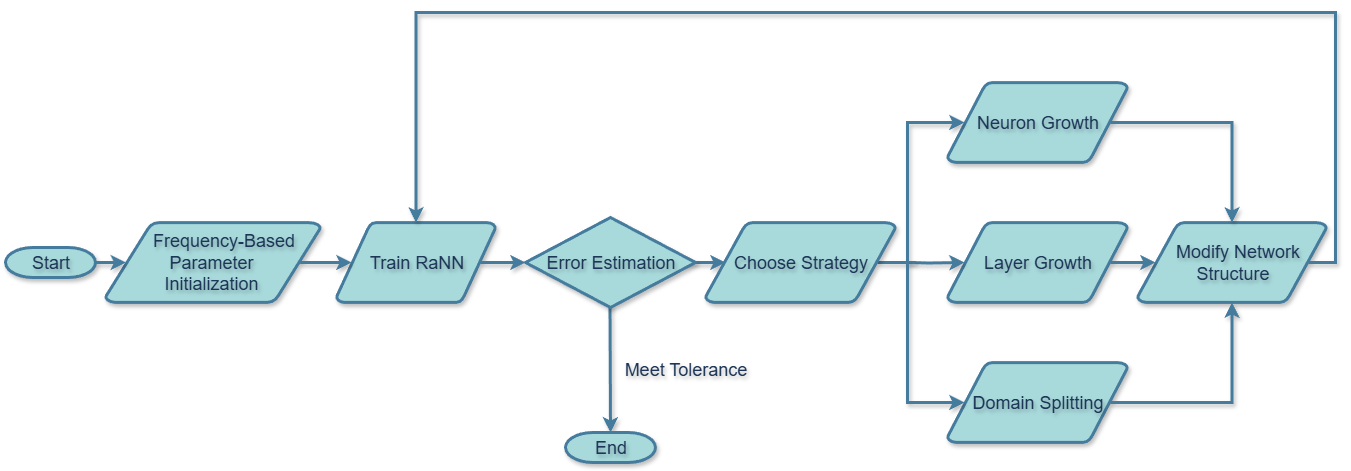}
    \caption{Flowchart of the AG-RaNN algorithm.}
    \label{fig:flow}
\end{figure}

In this section, we introduce the Adaptive-Growth RaNN (AG-RaNN) method for solving partial differential equations. The core idea behind AG-RaNN is to adapt the network structure and parameters based on the given PDE information and the previous RaNN solution, denoted by $u_\rho^0$. A flowchart illustrating the AG-RaNN algorithm is presented in Figure \ref{fig:flow}.

We introduce four strategies for constructing adaptive-growth RaNN: frequency-based parameter initialization, neuron growth, layer growth, and domain splitting. By employing these strategies, the network evolves from an initial state to a wider and deeper structure, significantly enhancing its approximation capability for the unknown solution. 

Let us consider a general partial differential equation:
\begin{align}
    \cG u&\,=\,f\quad \text{in }D, \label{PDE1}\\
    \cB u&\,=\,g\quad \text{on }\Gamma\subset\partial D. \label{PDE2}
\end{align}
Here, $\cG$ is a linear or nonlinear differential operator, and $\cB$ is the boundary operator. We approximate the unknown function $u$ by a network $u_\rho\in\cN^{\bcW,\bcB}_{\rho,1}(D)$ that minimizes the loss
\begin{align} \label{Eq:loss}
  \cL_\eta(u)\,= \,\|\cG u-f\|_{0,D}^2\,+\,\eta\|\cB u-g\|_{0,\Gamma}^2,
\end{align}
evaluated at a set of collocation points. In this paper, we use $\|\cdot\|_{s,D}$ and $\|\cdot\|_{s,\partial D}$ to denote the $H^s$-norms on the domain $D$ and its boundary $\partial D$, respectively, where $H^s$ is the standard Sobolev space (\cite{Adams1975Sobolev}). The notation $\|\cdot\|_l$ without the subscript $D$ or $\partial D$ refers to the standard vector $l$-norm. Since enforcing boundary conditions is critical for training RaNNs, we introduce an appropriate large penalty parameter $\eta$ to ensure they are satisfied.

\subsection{Frequency-Based Parameter Initialization Strategy} \label{Sec:Init}

Given a set of data points $\set{(\bx_i,y_i)}_{i=1}^N$, our goal is to find $u_\rho\in\cN^{\bcW,\bcB}_{\rho,1}(D)$ that approximates the function $u(\bx)$ such that $y_i=\cG(u_\rho)(\bx_i)$ (for instance, $y_i = f(x_i)$ in the PDE setting where $\cG u = f$). This requires constructing an appropriate set of basis functions, and we achieve this by generating initial weights using the Fourier transform.

We start with a shallow RaNN, where the RaNN space is defined as:
\begin{align} \label{rnn_space_init}
    \cN^{\bcW,\bcB}_{\rho,1}(D)\,=\,\set{z(\bx)\,=\,\bw\rho(\bW\bx+\bb)\,:\,\bx\in D\subset\Real^d, \, \forall\,\bw\in\Real^{1\times m_1}}
\end{align}
with fixed weights $\bW=(\bw_1,\bw_2,\cdots,\bw_d)\in\Real^{m_1\times d}$ and biases $\bb\in\Real^{m_1\times1}$. The weights are typically drawn from a uniform distribution $\cU(-r,r)$. Finding the optimal hyperparameter $r$, which influences approximation accuracy, can be difficult. To address this, we use a frequency-based parameter initialization technique.

We define a sufficiently large $r_\text{max}$ and divide it into intervals:
\begin{align*}
    r_k\,=\,\frac{k}{\Lambda}\,r_\text{max},\quad k=1,\cdots,\Lambda
\end{align*}
to explore different possible values. The goal is to find the optimal parameter $\br^\text{opt}=(r_1^\text{opt},\cdots,r_d^\text{opt})$ for the initial distribution. 

The candidate weight vector $\tilde{\bw}=(\tilde{w}_1,\cdots,\tilde{w}_d)$ consists of $d$ elements, where 
\begin{align} \label{Eq:candfunc}
    \tilde{w}_t = r_{k_t},\quad k_t\in\{1,\cdots,\Lambda\},\quad \forall\,1\le t\le d,
\end{align}
resulting in $\Lambda^d$ possible sets of weights. The smallest hypercube containing the domain $D$ is $H=[p_1,q_1]\times\cdots\times[p_d,q_d]$ with center $\bx_c=(\bp+\bq)/2$, where $\bp=(p_1,\cdots,p_d)^\top$ and $\bq=(q_1,\cdots,q_d)^\top$. The candidate functions are $\phi_j(\bx)=\rho(\tilde{\bw}^j(\bx-\bx_c))$ for $j=1,\cdots,\Lambda^d$. 

By performing a Fourier transform on $y$ and $\cG(\phi_j)(\bx)$ for each candidate function, we obtain $\hat{y}$ and $\widehat{\cG(\phi_j)}$ for $j=1,\cdots,\Lambda^d$. Let 
\begin{align*}
    \bxi_0\,=\,\mathop{\arg\max}\limits_{\bxi}\,\left|\hat{y}(\bxi)\right|\quad\text{and}\quad \bxi_j\,=\,\mathop{\arg\max}\limits_{\bxi}\,\left|\widehat{\cG(\phi_j)}(\bxi)\right|,
\end{align*}
from which we calculate
\begin{align*}
    j_0\,=\,\mathop{\arg\min}\limits_{j}\set{\|\bxi_j-\bxi_0\|_2\,:\,(\bxi_j)_t>(\bxi_0)_t,\,t=1,\cdots,d}.
\end{align*}
The parameter $\br^\text{opt}=\tilde{\bw}^{j_0}$ represents the near-optimal distribution parameter. Subsequently, we derive the near-optimal uniform distribution $\cU(-\br^\text{opt},\br^\text{opt})$, which is used to initialize the weights $\bW$. Specifically, $\bW\sim\cU(-\br^\text{opt},\br^\text{opt})$, meaning that each element in $\bw_i$ is drawn from a uniform distribution $\cU(-r^\text{opt}_i,r^\text{opt}_i)$.

\begin{remark}    
    The calculation of $j_0$ is explained here. Considering all candidate functions in the frequency domain via the Fourier transform, $\bxi_0$ and $\bxi_j$ denote the locations with the highest energy for different functions. To prioritize high-energy regions, $j_0$ is chosen by comparing the locations of the highest-energy frequencies of the target function and of the candidate function set.
\end{remark}

\begin{remark}
    If $\cG$ is a linear differential operator with constant coefficients, we can choose sine or cosine as the activation function, and the initialization can be achieved in a simpler way. For example, when $\cG=-\Delta$, we can obtain $-\widehat{\Delta u}$ and $\widehat{y}$:
    \begin{align*}
       - \widehat{\Delta u}(\bxi)=\|\bxi\|_2^2\,\widehat{u}(\bxi)\quad\text{and}\quad\widehat{u}(\bxi)=\frac{\widehat{f}(\bxi)}{\|\bxi\|_2^2}.
    \end{align*} 
    We can regard $\widehat{f}(\bxi)/\|\bxi\|_2^2$ as the density function to sample $W$. For simple cases, this strategy has a similar effect to the aforementioned one and is more efficient. However, due to its limitations on differential operators, we do not discuss this method in detail here.
\end{remark}

To construct biases $\bb$, we randomly sample points in the hypercube $H\supset D$. Specifically, we randomly generate $(\bb_1,\cdots,\bb_d)=\bB\in\Real^{m_1\times d}$ from the uniform distribution $\cU(\bp,\bq)$. The biases are then calculated as $\bb=-(\bW\odot\bB)\,\mathbf{1}_{d\times1}$, where $\odot$ denotes elementwise multiplication. 

\begin{remark}
When initializing an RaNN, selecting an appropriate activation function is critical. Typically, tanh is suitable for simple smooth functions, while sine or cosine can be effective for highly oscillatory functions. In this study, we also use the Gaussian function as the activation function due to its ability to approximate both low and high-frequency components of the solution by adjusting parameters. 
\end{remark}

\begin{remark} \label{re:zeroright}
When the source term vanishes, no meaningful information can be extracted for $\hat{y}(\bxi)$. In this case, we may simply set $j_0 = 1$ or choose it randomly. The subsequent approximation is then carried out using the neuron growth strategy described in the following subsection.
\end{remark}

\subsection{Neuron Growth Strategy}

After the initialization phase, a randomized neural network with a single hidden layer is created. In this subsection, we focus on growing the number of neurons in that layer. New neurons are added iteratively until the solution is accurate enough or further improvements become negligible.

Let $m_\text{add}$ denote the number of newly added neurons. The new RaNN space is:
\begin{align*}
    \cN^{\bcW,\bcB}_{\rho,1}(D)\,=\,\set{z(\bx)=\bw_1\rho(\bW\bx+\bb)+\bw_2\rho(\bW_\text{add}\bx+\bb_\text{add})\,:\,\bx\in D,\,\forall\,\bw_1\in\Real^{1\times m_1},\, \forall\,\bw_2\in\Real^{1\times m_\text{add}}},
\end{align*}
where $\bW$ and $\bb$ are determined during the initialization or the previous step. The objective is to determine the added weights $\bW_\text{add}$ and biases $\bb_\text{add}$.

By training the initial network, an approximate solution $u_\rho^0$ is obtained. While the initial network captures the high-energy frequencies, the expanded network should focus on the additional frequencies in the residual. To do this, we compute:
\begin{align*}
    y^*=\,y-\cG(u_\rho^0)(\bx),
\end{align*}
and the Fourier transform of $y^*$ can be directly computed as well. Subsequently, 
\begin{align*}
    \bxi_0^*\,=\,\mathop{\arg\max}\limits_{\bxi}\,\left|\widehat{y^*}(\bxi)\right|\quad \text{and}\quad \bxi_j\,=\,\mathop{\arg\max}\limits_{\bxi}\,\left|\widehat{\cG(\phi_j)}(\bxi)\right|
\end{align*}
can be obtained. Since $\bxi_j$ was previously calculated during initialization, only $\bxi_0^*$ needs to be calculated here, reducing computational costs. Finally, 
\begin{align*}
    j_0^*\,=\,\mathop{\arg\min}\limits_{j}\,\set{\|\bxi_j-\bxi_0^*\|_2\,:\,(\bxi_j)_t>(\bxi_0^*)_t,\,t=1,\cdots,d}
\end{align*}
is derived. Then, taking $\br^\text{opt}=\,\tilde{\bw}^{j_0^*}$, $\bW_\text{add}\sim\,\cU(-\br^\text{opt},\br^\text{opt})$, and the biases $\bb_\text{add}$ can be derived from $\cU(\bp,\bq)$ as discussed in the previous subsection. 

\begin{remark}
    Neuron growth can be regarded as a re-initialization step built upon the existing network, and it can be performed iteratively. The number of neurons added at each stage plays a crucial role: adding too few necessitates more iterations, whereas adding too many may increase the error. A practical strategy is to introduce more neurons at the beginning and gradually reduce the number in later stages.
\end{remark}

\begin{remark}
    We provide the distribution parameters $\br^\text{opt}$, while also introducing two additional hyperparameters, $r_\text{max}$ and $\Lambda$. Although these parameters have only a minor influence on the results, they can be adjusted during the neuron growth process.
\end{remark}

In addition to adding neurons, we can also prune ineffective ones after training. Neurons whose linear combination coefficients $\bw$ in \eqref{rnn_space_init} are zero are considered degenerate and can be removed. The remaining number of degrees of freedom is $m_1-m_p$, where $m_p$ is the number of degenerate neurons.

\subsection{Layer Growth Strategy} \label{Sec:layer-grow}

In this subsection, we discuss adding a new hidden layer before the output layer to improve the network's performance.

\begin{figure}[!htbp]
    \centering
    \includegraphics[width=\textwidth]{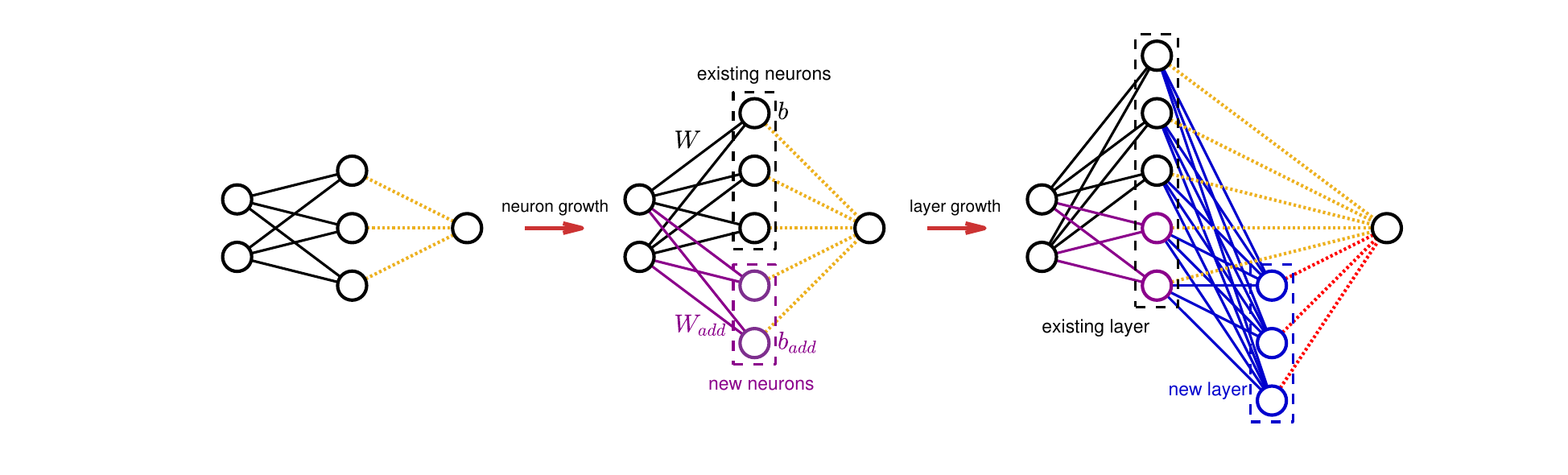}
    \caption{Initial RaNN, neuron growth, and RaNN with a shortcut connection.}
    \label{fig:RNN}
\end{figure}

Starting with an initial network that has a single hidden layer, we use shortcut connections to introduce a new hidden layer. Figure \ref{fig:RNN} illustrates the process: the left side shows the initial network, the middle shows neuron growth, and the right side shows the addition of the new hidden layer. The neurons in the new layer are depicted in blue, with solid lines representing fixed parameters and dotted lines showing the parameters to be trained.

With the aforementioned structure, several key considerations arise:
\begin{itemize}
    \item How many neurons to add to the new layer?
    \item Which activation functions to use for the new layer?
    \item How to generate the new parameters (solid blue lines) — whether to use random generation or a specific rule.
    \item Whether the existing parameters (dotted yellow lines) should be kept unchanged or trained alongside the new ones (dotted red lines).
\end{itemize}

Similar to neuron growth, additional neurons can be added to the new hidden layer in stages, balancing the number of neurons and the frequency of growth. The initial network typically captures low-frequency information, so local activation functions such as wavelets or Gaussians are useful for handling remaining errors. Training the existing parameters (dotted yellow lines) together with the new ones often yields better results, as it enriches the RaNN function space. A major focus is on determining the parameters represented by the solid blue lines.

The RaNN space after layer growth is defined as:
\begin{align} \label{rnn_space_layer}
    \cN^{\bcW,\bcB}_{\rho,1}(D)\,=\,\set{z(\bx)=\bw^{(1)}\bphi(\bx)+\bw^{(2)}\rho_2(\bW^{(1)}\bphi(\bx)+\bb^{(1)}):\bx\in D,\, \forall\,\bw^{(1)}\in\Real^{1\times m_1},\, \forall\,\bw^{(2)}\in\Real^{1\times m_2}},
\end{align}
where $\bphi(\bx)=\rho_1(\bW^{(0)}\bx+\bb^{(0)})$ represents a set of basis functions, with $\bW^{(0)}$ and $\bb^{(0)}$ determined in the previous step. The objective is to determine suitable $\bW^{(1)}$ and $\bb^{(1)}$. 

Using the RaNN, an approximate solution $u_\rho^0=\balpha\bphi$ is obtained, with $\balpha\in\Real^{1\times m_1}$ representing linear combination coefficients. The error estimator (e.g., the residual $|\cG(u_\rho^0)-y|$) is then computed at collocation points, and the top $m_2$ points $\bX_\text{err}=(\bx_1,\cdots,\bx_{m_2})\in\Real^{d\times m_2}$ with the largest errors are recorded. Each point corresponds to a neuron, and the activation function can be chosen as a local function. 
The new weights and biases are constructed as follows:
\begin{align*}
    \bW^{(1)}=\ &\bH\balpha\quad\text{and}\quad\bb^{(1)}=-\bH\odot [u_\rho^0(\bX_\text{err})]^\top,\\
    \bH=\sqrt{(\bH_0\odot\bH_0)\,\mathbf{1}_{d\times1}},\quad \bH_0&=(\bh_1^\top,\cdots,\bh_{m_2}^\top)^\top\in\,\Real^{m_2\times d},\quad \bh_i\sim\,\cU(-\br_2,\br_2),\quad\bh_i\in\,\Real^{1\times d}.
\end{align*}
Here, $\sqrt{\,\cdot\,}$ denotes the elementwise square root and $\mathbf{1}_{d\times1}$ represents the vector in $\mathbb{R}^{d\times 1}$ with all entries equal to $1$. The parameter $\br_2$ needs to be manually specified, which can significantly affect the outcome. To bypass this hyperparameter selection, we propose a criterion based on the gradient of the approximate solution $u_\rho^0$. We replace the randomly generated matrix with $\bH_0=\nabla u_\rho^0(\bX_\text{err}^\top)$ providing a more straightforward and generally effective choice.

\begin{figure}[!htbp]
    \centering
    \vspace{-0.7cm}
    \includegraphics[width=0.8\textwidth]{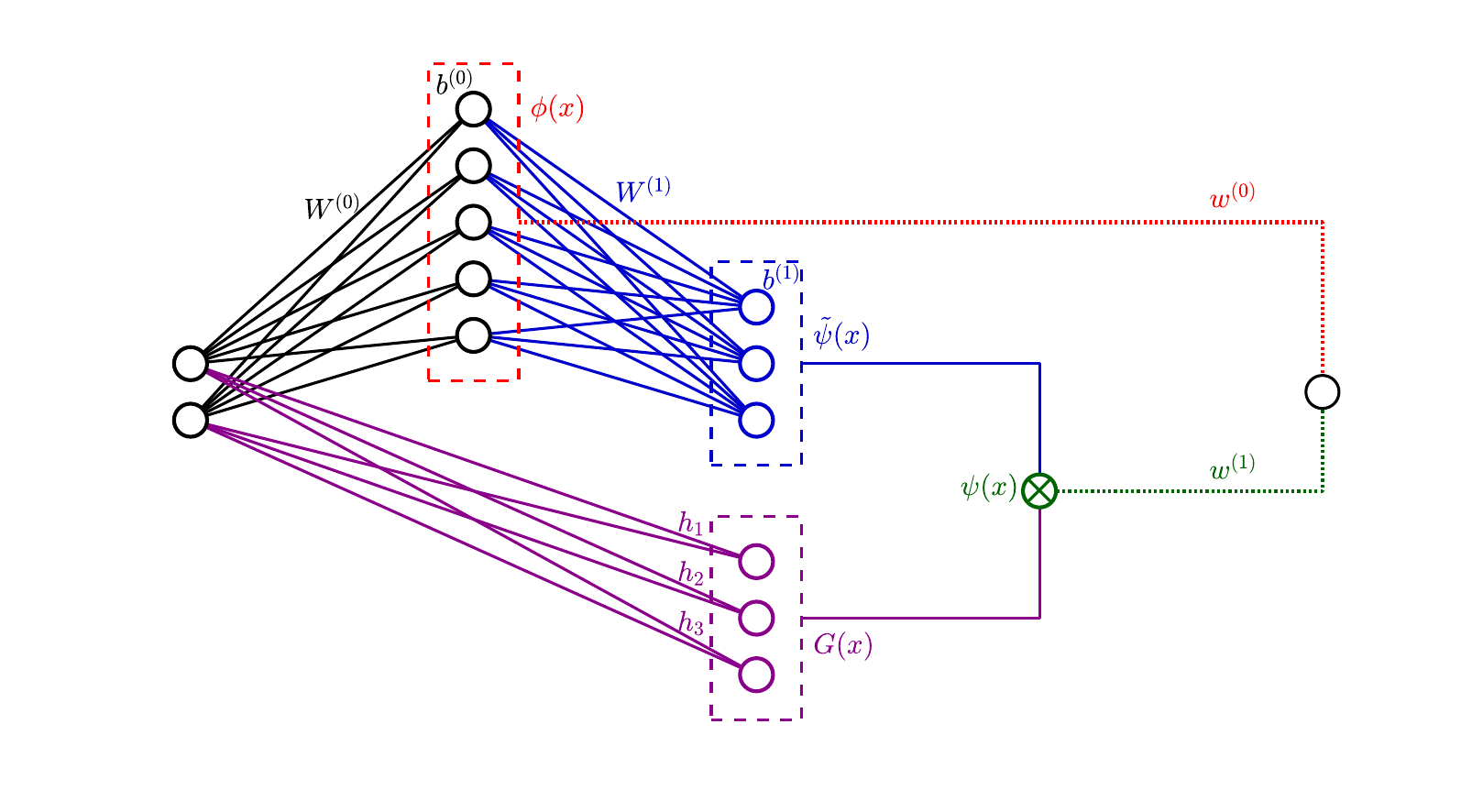}
    \vspace{-0.8cm}
    \caption{Randomized neural network with localized basis functions. }
    \label{fig:RNN-localization}
\end{figure}

Subsequently, a localization process is required. Let $\set{\tilde{\psi}_j(\bx)}_{j=1}^{m_2}$ represent the basis functions acquired through the aforementioned process. Each $\tilde{\psi}_j$ is localized by selecting a Gaussian function $G_j(\bx)$ as
\begin{align} \label{eq:localization}
    G_j(\bx)\,=\,\exp\left(-\frac{\left\|\,\bh_j \odot\,(\bx-\bx_j)^\top\right\|_2^2}{2}\right).
\end{align}
The final basis functions are $\bpsi=(\psi_1,\cdots,\psi_{m_2})^\top$, where $\psi_j(\bx)=G_j(\bx)\,\tilde{\psi}_j(\bx)$. Consequently, the modified randomized neural network space \eqref{rnn_space_layer} becomes:
\begin{align*}
    \cN^{\bcW,\bcB}_{\rho,1}(D)\,=\,\set{z(\bx)=\bw^{(1)}\bphi(\bx)+\bw^{(2)}\bpsi(\bx):\bx\in D, \,\forall\,\bw^{(1)}\in\Real^{1\times m_1},\, \forall\,\bw^{(2)}\in\Real^{1\times m_2}}.
\end{align*}

Figure \ref{fig:RNN-localization} presents a diagram of the RaNN with localized basis functions after layer growth. Compared to the network structure in the right of Figure \ref{fig:RNN}, the only difference is the purple neurons, representing the Gaussian function $G_j$. The green sign $\bigotimes$ signifies elementwise multiplication, equivalent to $\odot$. 

\subsection{Domain Splitting Strategy}

In this subsection, we develop a domain splitting strategy and utilize multiple RaNNs to approximate the unknown function.

Sometimes, the function $u(\bx)$ is too complex for the earlier growth strategies to be fully effective. For example, if $u(\bx)$ is discontinuous, a continuous NN may be insufficient. In such cases, employing multiple RaNNs becomes necessary. Each network can have its own inputs, parameters, and activation functions. Let $R$ denote the number of subnetworks. Within each subnetwork, both neuron growth and layer growth strategies can be applied to improve accuracy.

A crucial part of this strategy is partitioning the domain properly. We partition the domain by value ranges of the approximate solution $u_\rho^0(\bx)$. Since the function $u(\bx)$ is unknown, we first apply the RaNN method to get an approximate solution $u_\rho^0(\bx)$. Once the approximation becomes sufficiently informative to reveal the essential structure of the solution, we split $[v_{\min},v_{\max}]$ into $R$ segments:
\begin{align*}
    v_\text{min}=\,v_0<v_1<\cdots<v_{R-1}<v_R=\,v_\text{max}, 
\end{align*}
where $v_\text{min}=\mathop{\min}\limits_{\bx\in D}\,u_\rho^0(\bx)$ and $v_\text{max}=\mathop{\max}\limits_{\bx\in D}\,u_\rho^0(\bx)$. These segments denoted by $I_j=[v_{j-1},v_j]$ do not have to be equally spaced. Each subdomain corresponds to one of these intervals. Note that each subdomain may not be connected.

To link the RaNNs across subdomains, we sample points at the subdomain boundaries. If $u(\bx)$ is discontinuous, enforcing continuity at the interface is not required, and this step can be skipped.

However, if $u(\bx)$ is continuous, collocation points at the boundaries are needed to ensure continuity. Determining these points exactly can be challenging, so instead, we sample many collocation points $\set{\tilde{\bx}_k}_{k=1}^K$ across the domain and calculate $u_\rho^0(\tilde{\bx}_k)$. We then define the set of points at the interface between intervals $I_j$ and $I_{j+1}$ as
\begin{align*}
    \Theta_j=\set{\tilde{\bx}_k:v_j-\varepsilon_r\le u_\rho^0(\tilde{\bx}_k)\le v_j+\varepsilon_r}\quad j=1,\cdots,R-1,
\end{align*}
where $\varepsilon_r$ is a small constant. At these points, continuity can be enforced by ensuring $u_j(\tilde{\bx}_k)= u_{j+1}(\tilde{\bx}_k)$ for all $\tilde{\bx}_k\in\Theta_j$.

\section{Numerical Analysis}\label{sec:analysis}

In this section, we present an error analysis of the randomized neural network method, decomposing the total error into approximation, statistical, and optimization components. We work in an abstract Hilbert space setting and always assume that the differential operator $\cG$ is linear.

Let $V$ be a Hilbert space of candidate solutions on a bounded Lipschitz domain $D\subset\Real^d$, endowed with the norm $\|\cdot\|_V$. Let $Y$ and $Z$ be Hilbert spaces, and $\cG:V\to Y$ and $\cB:V\to Z$ be bounded linear operators representing, respectively, the PDE operator and the boundary operator. Denote by $u^*\in V$ the (unique) exact solution of
\begin{align*}
    \cG u^* = f\in Y, \qquad \cB u^* = g\in Z .
\end{align*}

We consider a single hidden layer ($L=1$) RaNN function space
\begin{align*}
  \cN_\rho(D) = \Bigl\{u(\bx)=\balpha\rho(\bW\bx+\bb):\bx\in D,\ \|\balpha\|_1\le C_N\Bigr\}\subset V,
\end{align*}
where $C_N>0$ is the constant in Theorem~\ref{Thm:H1bound-UAT}. The best approximation of $u^*$ in $\cN_\rho(D)$ with respect to the $V$-norm is
\begin{align*}
    u_a := \mathop{\arg\min}\limits_{u\in\cN_\rho(D)} \|u-u^*\|_V .
\end{align*}

\paragraph{Population and empirical losses.}
For $u\in V$ we define the (population) loss
\begin{align}\label{eq:pop-loss-general}
  \cL(u) := \|\cG u - f\|_Y^2 + \|\cB u - g\|_Z^2 .
\end{align}
In practice, the norms on $Y$ and $Z$ are evaluated (or approximated) by quadrature formulas based on interior and boundary collocation points $\{\bx^I_i\}_{i=1}^{N^I}\subset D$ and $\{\bx^B_i\}_{i=1}^{N^B}\subset\Gamma\subset\partial D$. For example, when $Y=L^2(D)$ and $Z=L^2(\Gamma)$ the corresponding empirical loss takes the familiar form
\begin{align}\label{eq:emp-loss-L2}
  \widehat{\cL}(u) := \frac{|D|}{N^I}\sum_{i=1}^{N^I}\left|\cG u(\bx^I_i)-f(\bx^I_i)\right|^2 + \frac{|\Gamma|}{N^B}\sum_{i=1}^{N^B}\left|\cB u(\bx^B_i)-g(\bx^B_i)\right|^2 ,
\end{align}
while for other choices of $Y$ and $Z$ the quadrature is understood as a numerical approximation of the corresponding Hilbert norms. The empirical minimizer
\begin{align*}
    u_\rho:= \mathop{\arg\min}\limits_{u\in\cN_\rho(D)}\,\widehat{\cL}(u)
\end{align*}
can then be obtained by least-squares.

\begin{assumption}[Graph norm equivalence]\label{ass:graph}
We assume that there exist constants $C_L,C_U>0$ such that
\begin{align}\label{eq:equiv-general}
  C_L\|v\|_V^2\;\le\;\|\cG v\|_Y^2 + \|\cB v\|_Z^2\;\le\;C_U\|v\|_V^2,\qquad \forall\,v\in V .
\end{align}
\end{assumption}

Under Assumption \ref{ass:graph}, and using that $\cL(u^*)=0$, we obtain for every $u_\rho\in\cN_\rho(D)$
\begin{align}\label{eq:equiv-urho-general}
  C_L\|u_\rho-u^*\|_V^2\;\le\;\cL(u_\rho)\;\le\;C_U\|u_\rho-u^*\|_V^2 .
\end{align}
Consequently,
\begin{align*}
  C_L\|u_\rho-u^*\|_V^2&\le \cL(u_\rho)=\cL(u_\rho)-\cL(u^*) \\
  &=[\cL(u_\rho)-\widehat{\cL}(u_\rho)]+[\widehat{\cL}(u_\rho)-\widehat{\cL}(u_a)]+[\widehat{\cL}(u_a)-\cL(u_a)]+[\cL(u_a)-\cL(u^*)] \\
  &\le[\cL(u_a)-\cL(u^*)]+2\sup_{u\in\cN_\rho(D)}\bigl|\cL(u)-\widehat{\cL}(u)\bigr|+[\widehat{\cL}(u_\rho)-\widehat{\cL}(u_a)] \\
  &\le C_U\|u_a-u^*\|_V^2+2\sup_{u\in\cN_\rho(D)}\bigl|\cL(u)-\widehat{\cL}(u)\bigr|+\underbrace{\bigl[\widehat{\cL}(u_\rho)-\widehat{\cL}(u_a)\bigr]}_{\le 0},
\end{align*}
where $\|u_a-u^*\|_V$ is the approximation error, $\sup_{u\in\cN_\rho(D)}|\cL(u)-\widehat{\cL}(u)|$ is the statistical term, and $\widehat{\cL}(u_\rho)-\widehat{\cL}(u_a)$ is the optimization term, which vanishes if $u_\rho$ is computed as the exact minimizer of $\widehat{\cL}$.

\begin{remark}[Typical choices of $(V,Y,Z)$ for PDE models]
For concrete PDEs, the spaces $V$, $Y$, and $Z$ are chosen so that the graph norm equivalence~\eqref{eq:equiv-general} follows from standard stability and regularity theory. We briefly list a few typical examples.
\begin{itemize}
    \item \textbf{Second-order elliptic problems.}For a linear elliptic operator with homogeneous Dirichlet boundary conditions, one may take
    \begin{align*}
        V = H_0^1(D),\qquad Y = H^{-1}(D),\qquad Z = \{0\},
    \end{align*}
    where $\cG:V\to Y$ is the variational operator and $\cB\equiv 0$. Then \eqref{eq:equiv-general} reduces to the usual energy estimate $\|v\|_{H^1(D)} \simeq \|\cG v\|_{H^{-1}(D)}$. For elliptic problems with $H^2$-regularity and strong form residuals, one may instead take, for example,
    \begin{align*}
        V = H^2(D),\qquad Y = L^2(D),\qquad Z = H^{3/2}(\Gamma)
    \end{align*}
    with $\cB v$ being the Dirichlet trace, in which case
    \eqref{eq:equiv-general} follows from standard $H^2$-regularity
    results.

    \item \textbf{Linear parabolic problems.}
    Consider, for instance, a heat-type equation on $Q=(0,T)\times D$ with homogeneous Dirichlet boundary conditions and prescribed initial data. A natural choice is
    \begin{align*}
        V = \left\{u:u\in L^2((0,T);H_0^1(D)),\,\partial_t u\in L^2((0,T);H^{-1}(D))\right\},
    \end{align*}
    endowed with the norm $\|u\|_V^2 = \|u\|_{L^2((0,T);H_0^1(D))}^2 +\|\partial_t u\|_{L^2((0,T);H^{-1}(D))}^2$. One may take
    \begin{align*}
        Y = L^2((0,T);H^{-1}(D)),\qquad Z = L^2(D),
    \end{align*}
    with $\cG u = \partial_t u - \cA u$ (where $\cA$ is the elliptic part) and $\cB u = u(0,\cdot)$ encoding the initial condition. Standard parabolic energy estimates then yield an inequality of the form~\eqref{eq:equiv-general} on $V$.

    \item \textbf{Linear hyperbolic problems.}
    For second-order hyperbolic equations such as the wave equation, it is convenient to rewrite the PDE as a first-order evolution system in time and to work in the corresponding energy space $V$ of trajectories $(u,\partial_t u)$. In this setting, one can choose $Y$ as an $L^2$ in time space for the residual and $Z$ for the initial data so that the classical energy estimates yield a graph norm equivalence of the form \eqref{eq:equiv-general}. We omit the precise functional setup, as our analysis only uses the abstract Assumption \ref{ass:graph}.
\end{itemize}

In all these cases, our subsequent analysis only uses the abstract form \eqref{eq:equiv-general}; the empirical loss \eqref{eq:emp-loss-L2} should be understood as a concrete realization of the corresponding $Y$- and $Z$-norms via quadrature in the particular case $Y=L^2(D)$ and $Z=L^2(\Gamma)$, or as their discrete approximations for other choices of $Y$ and $Z$.
\end{remark}

\subsection{Approximation Error}

In this part, we only consider the case where $V$ is a Sobolev space $H^s$. In order to estimate the approximation error, we first present the universal approximation theorem (UAT) for $H^{s+\varepsilon_C}$ functions with bounded weights. Consider a single hidden layer network:
\begin{align*}
  u_N(\bx)\,=\,\sum_{k=1}^N\,c_k\rho(\bw_k^\top \bx + b_k),\qquad \bw_k\in\Real^d,\ b_k,c_k\in\Real,
\end{align*}
with activation function $\rho:\Real\to\Real$.
\begin{assumption}[Activation admissibility]\label{As:act}
The activation function $\rho$ satisfies:
\begin{itemize}
    \item[(A1)] $\rho\in C^s(\Real)$, and $\|\rho^{(s_0)}\|_{L^\infty(\Real)}\,<\,\infty$ for all $s_0=0,1,\cdots,s$.
    \item[(A2)] $\widehat{\rho'}(1)\neq 0$ and $\rho^{(s_0)}\in L^1(\Real)$ for all $s_0=1,\cdots,s+1$.
\end{itemize}
\end{assumption}

\begin{remark}
(A2) is used below to justify truncation in the bias variable and to control $H^s$ errors. Condition $\widehat{\rho'}(1)\neq 0$ can be replaced by ``$\widehat{\rho'}(t_0)\neq 0$ for some $t_0>0$'' after a harmless rescaling of the $b$-variable. This class includes tanh, logistic functions and other smooth sigmoids. However, ReLU and ReLU$^k$ functions do not belong to this function class. Nevertheless, there exists an alternative proof demonstrating that ReLU-type functions also satisfy the universal approximation property (\cite{Liu2025ReLU}).
\end{remark}

\begin{lemma}[Extension, \cite{Adams1975Sobolev}]\label{lem:ext}
For a bounded Lipschitz domain $D\subset\Real^d$, there exists a bounded linear extension operator $E:H^{s+\varepsilon_C}(D)\to H^{s+\varepsilon_C}(\Real^d)$ such that $(Eu^*)|_D=u^*$ and $\|Eu^*\|_{s+\varepsilon_C,\Real^d}\,\le\,C_E(\varepsilon_C,D,d)\,\|u^*\|_{s+\varepsilon_C,D}$.
\end{lemma}

We adopt the Fourier convention
\begin{align*}
    \widehat{g}(\bxi)\,=\,\int_{\Real^d}\,g(\bx)\,e^{-i\,\bxi^\top\bx}\,\dbx, \qquad g(\bx)\,=\,(2\pi)^{-d}\int_{\Real^d}\,\widehat{g}(\bxi)\,e^{i\,\bxi^\top\bx}\,\mathrm{d}\bxi.
\end{align*}
For $s\in\Real$, the (inhomogeneous) Sobolev norm is
\begin{align*}
    \|g\|_{s,\Real^d}^2\,:=\,\int_{\Real^d}\,(1+\|\bxi\|^2_2)^s\,|\widehat{g}(\bxi)|^2\,\mathrm{d}\bxi.
\end{align*}

\begin{lemma}[Bandlimit error]\label{lem:band}
Let $v\in H^{s+\varepsilon_C}(\Real^d)$ with $\varepsilon_C>0$. For $R>0$, let $P_R$ be the Fourier projector with symbol $\mathbf{1}_{\set{\|\bxi\|_2\le R}}$ and set $v_R:=P_Rv$. Then $v_R\in H^s(\Real^d)$, $\mathrm{supp}\,\widehat{v_R}\subset\set{\|\bxi\|_2\le R}$, and
\begin{align}\label{eq:H1-band-error}
  \|v-v_R\|_{s,\Real^d}\,\le\,R^{-\varepsilon_C}\|v\|_{s+\varepsilon_C,\Real^d} .
\end{align}
\end{lemma}

\begin{proof}
By the definition of Sobolev norm,
\begin{align*}
    \|v-v_R\|_{s,\Real^d}^s\,=\,\int_{\|\bxi\|_2>R}\,(1+\|\bxi\|_2^2)^s\,|\widehat v(\bxi)|^2\,\mathrm{d}\bxi.
\end{align*}
For $\|\bxi\|_2>R$ and any $\varepsilon_C>0$,
\begin{align*}
    (1+\|\bxi\|_2^2)^s\,=\,(1+\|\bxi\|_2^2)^{s+\varepsilon_C}(1+\|\bxi\|_2^2)^{-\varepsilon_C}\,\le\,(1+\|\bxi\|_2^2)^{s+\varepsilon_C}\|\bxi\|_2^{-2\varepsilon_C}\,\le\,(1+\|\bxi\|_2^2)^{s+\varepsilon_C}R^{-2\varepsilon_C},
\end{align*}
where we used $(1+\|\bxi\|_2^2)\ge \|\bxi\|_2^2$ and $\|\bxi\|_2>R$. Therefore,
\begin{align*}
    \|v-v_R\|_{s,\Real^d}^2\,\le\,R^{-2\varepsilon_C} \int_{\Real^d}\,(1+\|\bxi\|_2^2)^{s+\varepsilon_C}\,|\widehat v(\bxi)|^2\,\mathrm{d}\bxi\,=\,R^{-2\varepsilon_C}\,\|v\|_{s+\varepsilon_C,\Real^d}^2,
\end{align*}
which implies \eqref{eq:H1-band-error} after taking square roots.
\end{proof}

We now show that a bandlimited function $v_R$ admits an integral representation with bounded ridge directions.

\begin{theorem}[Integral representation]\label{thm:ridgelet-synthesis}
Let $v_R\in H^s(\Real^d)$ satisfy $\mathrm{supp}\,\widehat{v_R}\subset\set{\bxi:\|\bxi\|_2\le R}$ with $R>0$ and $C_X:=\sup\limits_{\bx\in D}\|\bx\|_2$. Then, for any $\varepsilon_R>0$, under Assumption \ref{As:act}, there exist a bound $B_R = B_R(\varepsilon_R,s,R,d,C_X,\rho)>0$ and a finite signed Borel measure 
$\nu_R$ supported on
\begin{align*}
  \Theta_R\,:=\,\set{(\bw,b)\in\Real^d\times\Real:\|\bw\|_2\le R,\,|b|\le B_R},
\end{align*}
such that
\begin{align}\label{eq:ridge-synth}
    \left\|v_R(\bx) - \int_{\Theta_R}\!\rho(\bw^\top \bx + b)\,\mathrm{d}\nu_R(\bw,b)\right\|_{s,D}\,\le\,\varepsilon_R\quad \text{ for any fixed compact }D\subset\Real^d,
\end{align}
and
\begin{align*}
  \|\nu_R\|_{\mathrm{TV}}\,\le\,C_T\,\|v_R\|_{0,\Real^d},
\end{align*}
where $C_T>0$ depends only on $R$, $B_R$, $d$ and $\rho$. 
\end{theorem}

\begin{proof}
By inverse Fourier transform,
\begin{align}\label{eq:invFT}
    v_R(\bx)\,=\,\text{Re}\left\{\frac{1}{(2\pi)^d}\int_{\|\bw\|_2\le R}\,\widehat v_R(\bw)\,e^{i\,\bw^\top\bx}\,\mathrm{d}\bw\right\}.
\end{align}
Setting $u=\bw^\top \bx+b$ gives
\begin{align*}
  \int_{\Real}\,\rho'(\bw^\top \bx+b)e^{-i\,b}\,\mathrm{d}b\,=\,e^{i\,\bw^\top \bx}\,\int_{\Real}\,\rho'(u)\,e^{-iu}\,\mathrm{d}u\,=\,e^{i\,\bw^\top \bx}\,\widehat{\rho'}(1) ,
\end{align*}
and
\begin{align}\label{eq:1D-id}
  e^{i\,\bw^\top \bx}\,=\,\frac{1}{\widehat{\rho'}(1)}\,\int_{\Real}\,\rho'(\bw^\top \bx +b)\,e^{-i\,b}\,\mathrm{d}b
\end{align}
follows $\widehat{\rho'}(1)\neq 0$ by (A2).

By Cauchy--Schwarz and $\mathrm{supp}\,\widehat{v_R}\subset\set{\|\bw\|_2\le R}$,
\begin{align*}
    \int_{\|\bw\|_2\le R}|\widehat{v_R}(\bw)|\,\mathrm{d}\bw\,\le\,\left|\set{\|\bw\|_2\le R}\right|^{1/2}\,\|\widehat{v_R}\|_{0,\Real^d}\le C_d\,R^{d/2}\,\|v_R\|_{0,\Real^d},
\end{align*}
where $C_d$ is the volume of the $d$-dimensional unit ball and $\max\limits_d\,C_d<6$. Since $\rho'\in L^1(\Real)$ and $|e^{-ib}|=1$,
\begin{align*}
    \int_{\|\bw\|_2\le R}\int_{\Real}\,|\rho'(\bw\cdot \bx+b)|\,|\widehat{v_R}(\bw)|\,|e^{-ib}|\,\mathrm{d}b\,\mathrm{d}\bw\le \|\rho'\|_{L^1}\int_{\|\bw\|_2\le R} |\widehat{v_R}(\bw)|\,\mathrm{d}\bw\,<\,\infty.
\end{align*}
Insert \eqref{eq:1D-id} into \eqref{eq:invFT} and apply Fubini--Tonelli:
\begin{align}\label{eq:rep-infinite-b}
    v_R(\bx)\,=\,\text{Re}\left\{\frac{1}{(2\pi)^d\,\widehat{\rho'}(1)} \int_{\|\bw\|_2\le R}\widehat{v_R}(\bw)\left(\int_{\Real}\,\rho'(\bw^\top \bx + b)e^{-i\,b}\,\mathrm{d}b\right) \mathrm{d}\bw\right\}.
\end{align}

Representation \eqref{eq:rep-infinite-b} integrates $b$ over $\Real$. We now truncate $b$ to a bounded interval, with an error estimate in $H^s$ on any fixed $D\subset\Real^d$. For $(\bw,b)$ with $\|\bw\|_2\le R$ and $|b|>B$, taking $u=\bw^\top \bx+b$, 
\begin{align*}
  \int_{|b|>B} \sum_{|\bm{\gamma}|=s_0}\left|\partial_{\bgamma}\rho'(\bw^\top \bx + b)\right|\,\mathrm{d}b\,\le\,\|\bw\|_1^{s_0}\int_{|u|>B-RC_X} |\rho^{(s_0+1)}(u)|\,\mathrm{d}u\,\le\,(\sqrt{d}R)^{s_0}\int_{|u|>B-RC_X} |\rho^{(s_0+1)}(u)|\,\mathrm{d}u,
\end{align*}
where $\bgamma=(\gamma_1,\cdots,\gamma_d)$ is multiple index and $s_0=0,1,\cdots,s$. Therefore, for the $H^s(D)$ norm we obtain the uniform bound
\begin{align*}
  \sup_{\|\bw\|_2\le R}\sup_{\bx\in D}\int_{|b|>B} \sum_{s_0=0}^s\sum_{|\bm{\gamma}|=s_0}\left|\partial_{\bgamma}\rho'(\bw^\top \bx + b)\right|\mathrm{d}b\,\le\,\sum_{s_0=0}^s(\sqrt{d}R)^{s_0}\int_{|u|>B-RC_X} |\rho^{(s_0+1)}(u)|\,\mathrm{d}u.
\end{align*}
Since $\rho^{(s_0)}\in L^1(\Real)$ for $s_0=1,\cdots,s+1$, the right hand side $\to 0$ as $B\to\infty$.

For any $\varepsilon_R>0$, there is $B_R = B_R(\varepsilon_R,s,R,d,C_X,\rho)$ such that 
\begin{align*}
    v_B(\bx)&\,=\,\text{Re}\left\{\frac{1}{(2\pi)^d\,\widehat{\rho'}(1)} \int_{\|\bw\|_2\le R}\widehat{v_R}(\bw)\left(\int_{|b|\le B_R}\rho'(\bw^\top \bx + b)e^{-i\,b}\,\mathrm{d}b\right) \mathrm{d}\bw\right\}\\
    &\,=\,\text{Re}\left\{\frac{1}{(2\pi)^d\,\widehat{\rho'}(1)} \int_{\|\bw\|_2\le R}\widehat{v_R}(\bw)\set{\left[\rho(\bw^\top\bx+b)e^{-i\,b}\right]_{-B_R}^{B_R}+i\int_{|b|\le B_R}\rho(\bw^\top \bx + b)e^{-i\,b}\,\mathrm{d}b}\mathrm{d}\bw\right\}.
\end{align*}
We thus define the measure $\nu_R$ on $\Theta_R:=\set{(\bw,b): \|\bw\|_2\le R,\,|b|\le B_R}$:
\begin{align*}
    \mathrm{d}\nu_R(\bw,b):=\text{Re}\left\{\frac{1}{(2\pi)^d\,\widehat{\rho'}(1)}\,\widehat{v_R}(\bw)\left(e^{-iB_R}\,\delta_{b=B_R}-e^{iB_R}\,\delta_{b=-B_R}+ie^{-i\,b}\,\mathrm{d}b\right)\mathrm{d}\bw\right\},
\end{align*}
where $\delta$ is Dirac function. Consequently,  
\begin{align*}
    \left\|v_R - \int_{\Theta_R}\!\rho(\bw^\top \bx + b)\,\mathrm{d}\nu_R(\bw,b) \right\|_{s,D}=\,\|v_R-v_B\|_{s,D}\,\le\,\varepsilon_R.
\end{align*}

Therefore, for the total variation of $\nu_R$, we have
\begin{align*}
  \|\nu_R\|_{\mathrm{TV}}&\,\le\,\frac{1}{(2\pi)^d\,|\widehat{\rho'}(1)|}\int_{\|\bw\|_2\le R}|\widehat{v_R}(\bw)|\,\left|e^{-i\,B_R}-e^{i\,B_R}+\int_{|b|\le B_R}ie^{-i\,b}\,\mathrm{d}b\right|\mathrm{d}\bw\\
  &\,\le\,\frac{2(1+B_R)}{(2\pi)^d\,|\widehat{\rho'}(1)|}\int_{\|\bw\|_2\le R}|\widehat{v_R}(\bw)|\,\mathrm{d}\bw\,<\, \frac{2C_d(1+B_R)R^{d/2}}{(2\pi)^d\,|\widehat{\rho'}(1)|}\,\|v_R\|_{0,\Real^d},
\end{align*}
which implies $\nu_R$ is a finite signed measure.
\end{proof}

\begin{remark}
If $\widehat{\rho'}(t_0)\neq 0$ for some $t_0>0$, then the identity \eqref{eq:1D-id} reads $e^{i\,\bw^\top \bx}=\widehat{\rho'}(t_0)^{-1} \int_{\Real} \rho'(\bw^\top \bx +b)\,e^{-it_0b}\,\mathrm{d}b$. One can reduce to the $t_0=1$ case by the change of variables $b'=t_0 b$ and the reparameterization of the measure in $b$. This does not affect the bound on $\|\bw\|_2$ and only rescales the bias truncation threshold by a factor $t_0$.
\end{remark}

\begin{lemma}\label{lem:quad}
Under Assumption \ref{As:act} and Theorem \ref{thm:ridgelet-synthesis}, for any $\varepsilon_D>0$, there exist nodes $\theta_1,\cdots,\theta_N\in\Theta_R$ and weights $c_1,\cdots,c_N\in\Real$ such that
\begin{align*}
    \left\| v_R - \sum_{k=1}^N c_k\phi_{\theta_k}\right\|_{s,D} <\,\varepsilon_D,\text{ where }\phi_\theta(\bx)\,=\,\rho(\bw^\top \bx + b),\,\theta=(\bw,b)
\end{align*}
with $\max\limits_k \|\bw_k\|_2\le R$ and $\max\limits_k\,|b_k|\le B_R$. Moreover, one can ensure
\begin{align*}
    \sum_{k=1}^N\,|c_k|\,\le\,\|\nu_R\|_{\mathrm{TV}}.
\end{align*}
\end{lemma}

\begin{proof}
Since $D$ has finite measure, each $\phi_\theta$ lies in $H^s(D)$ and
\begin{align}\label{eq:H1-uniform-bound}
  \|\phi_\theta\|_{s,D}\,\le\, \sqrt{|D|}\,\sqrt{\sum_{s_0=0}^sR^{2s_0}\left\|\rho^{(s_0)}\right\|_{L^\infty,\Real}^2}\,:=\,C_D<\infty.
\end{align}

Let $\nu_R=\nu_R^+-\nu_R^-$ be the Jordan decomposition of the finite signed measure, and $|\nu_R|=\nu_R^++\nu_R^-$ its total variation measure. Let
\begin{align*}
    \sigma:=\frac{\mathrm{d}\nu_R}{\mathrm{d}|\nu_R|}
\end{align*}
be the Radon--Nikod\'ym derivative so that $|\sigma|\le 1$ $|\nu_R|$-a.e.. Define the probability measure $\mu:=|\nu_R|/\|\nu_R\|_{\mathrm{TV}}$ on $\Theta_R$. Then, in the Bochner sense in $H^s(D)$,
\begin{align}\label{eq:Bochner}
  v_B\,=\,\int_{\Theta_R} \phi_\theta\,\mathrm{d}\nu_R(\theta)\,=\,\int_{\Theta_R} \sigma(\theta)\,\phi_\theta\,\mathrm{d}|\nu_R|(\theta)\,=\,\|\nu_R\|_{\mathrm{TV}}\,\mathbb{E}_\mu\big[Z(\Theta)\big],
\end{align}
where we view the $H^s(D)$-valued random variable $Z(\Theta):=\sigma(\Theta)\,\phi_{\Theta}$ with $\Theta\sim\mu$. By \eqref{eq:H1-uniform-bound} and $|\sigma|\le 1$, $\|Z(\Theta)\|_{s,D}\le C_D$ $\mu$-a.s.

Draw independent and identically distributed samples $\theta_1,\cdots,\theta_N\sim\mu$, and define the empirical average in $H^s(D)$:
\begin{align*}
  \overline{Z}_N\,:=\,\frac{1}{N}\sum_{k=1}^N Z(\theta_k)\,=\,\frac{1}{N}\sum_{k=1}^N \sigma(\theta_k)\,\phi_{\theta_k}.
\end{align*}
Then $\mathbb{E}\,[\overline{Z}_N] = \mathbb{E}_\mu [Z(\Theta)]$ in $H^s$ and $H^s$ is a Hilbert space,
\begin{align}\label{eq:MC-variance}
  \mathbb{E}\,\left[\big\|\overline{Z}_N - \mathbb{E}_\mu\,\left[Z(\Theta)\right]\big\|_{s,D}^2\right]\,=\,\frac{1}{N}\,\mathbb{E}_\mu\,\left[\|Z(\Theta)-\mathbb{E}_\mu\,\left[Z(\Theta)\right]\|_{s,D}^2\right]\,\le\,\frac{1}{N}\,\mathbb{E}_\mu\,\left[\|Z(\Theta)\|_{s,D}^2\right]\,\le\,\frac{C_D^2}{N}.
\end{align}
Multiplying \eqref{eq:MC-variance} by $\|\nu_R\|_{\mathrm{TV}}^2$ and using \eqref{eq:Bochner},
\begin{align*}
  \mathbb{E}\left[\left\|v_B - \|\nu_R\|_{\mathrm{TV}}\overline{Z}_N \right\|_{s,D}^2\right]\,\le\,\frac{C_D^2\,\|\nu_R\|_{\mathrm{TV}}^2}{N}.
\end{align*}
Hence, by the probabilistic method, there exists a realization 
$\set{\theta_k}_{k=1}^N$ such that
\begin{align}\label{eq:exist-empirical}
  \left\| v_B - \|\nu_R\|_{\mathrm{TV}}\overline{Z}_N \right\|_{s,D}\,\le\,\frac{C_D\,\|\nu_R\|_{\mathrm{TV}}}{\sqrt{N}}.
\end{align}
Choose
\begin{align*}
  N\,\ge\,\left(\frac{2C_D\,\|\nu_R\|_{\mathrm{TV}}}{\varepsilon_D}\right)^2
\end{align*}
and take $\varepsilon_R=\varepsilon_D/2$ in Theorem \ref{thm:ridgelet-synthesis}, we obtain $\|v_R-\|\nu_R\|_{\mathrm{TV}}\overline{Z}_N \|_{s,D} \le \varepsilon_D$.

Write each sampled node as $\theta_k=(\bw_k,b_k)\in\Theta_R$ and set coefficients
\begin{align*}
  c_k:=\frac{\|\nu_R\|_{\mathrm{TV}}}{N}\,\sigma(\theta_k)\in\Real .
\end{align*}
Then
\begin{align*}
    \|\nu_R\|_{\mathrm{TV}}\,\overline{Z}_N\,=\,\sum_{k=1}^N\,c_k\phi_{\theta_k},
\end{align*}
and we have the error bound 
\begin{align*}
    \left\|v_R-\sum_{k=1}^N c_k\phi_{\theta_k}\right\|_{s,D}\,\le\,\varepsilon_D.
\end{align*}
Since $\theta_k\in\Theta_R$, we inherit $\|\bw_k\|_2\le R$ and $|b_k|\le B_R$. Moreover,
\begin{align*}
  \sum_{k=1}^N\,|c_k|\,=\,\frac{\|\nu_R\|_{\mathrm{TV}}}{N}\sum_{k=1}^N\,|\sigma(\theta_k)|\,\le\,\|\nu_R\|_{\mathrm{TV}},
\end{align*}
because $|\sigma|\le 1$. This proves the lemma.
\end{proof}

\begin{theorem}[Bounded-parameter $H^s$-UAT]\label{Thm:H1bound-UAT}
Under Assumption \ref{As:act}, for any $\varepsilon_A>0$, there exists constants $R=R(\varepsilon_A,\varepsilon_C,d,D,\|u^*\|_{s+\varepsilon_C,D})$, $B_R=B_R(\varepsilon_A,\varepsilon_C,s,d,D,\rho,\|u^*\|_{s+\varepsilon_C,D})$, $C_N=C_N(\varepsilon_A,\varepsilon_C,s,d,D,\rho,\|u^*\|_{s+\varepsilon_C,D})$ and a single hidden layer network
\begin{align*}
  u_N(\bx)\,=\,\sum_{k=1}^N\,c_k\rho(\bw_k^\top \bx + b_k),
\end{align*}
such that
\begin{align*}
  \|u^*-u_N\|_{s,D}\,<\,\varepsilon_A/2,
\end{align*}
and the parameters are bounded:
\begin{align*}
  \sum_{k=1}^N|c_k|\,\le\,C_N,\qquad \|\bw_k\|_2\,\le\,R,\qquad |b_k|\,\le\,B_R,\qquad 1\le k\le N.
\end{align*}
In particular, there is a box $\Theta=[-r,r]^{d+1}$ with $r=\max\set{R,B_R}$ that contains all inner parameters.
\end{theorem}

\begin{proof}
Fix $\varepsilon_A>0$, by Lemma \ref{lem:ext} and Lemma \ref{lem:band}, choose $R>0$ and use monotonicity of the norm to get the bounds on $H^s(D)$:
\begin{align*}
    \|u^*-v_R\|_{s,D}\,=\,\|Eu^*-v_R\|_{s,D}\,\le\,\|Eu^*-v_R\|_{s,\Real^d}\,\le\,\varepsilon_A/4.
\end{align*}
Apply Theorem \ref{thm:ridgelet-synthesis} to obtain a measure $\nu_R$ supported in $\Theta_R$ with $\|\bw\|_2\le R$ and $|b|\le B_R$. By Lemma \ref{lem:quad}, choose a finite sum (neural network) $u_N$ so that $\|v_R-u_N\|_{s,D}\le \varepsilon_A/4$ on $D$, with parameters bounded as stated. Finally,
\begin{align*}
    \|u^*-u_N\|_{s,D}\,\le\,\|u^*-v_R\|_{s,D} + \|v_R-u_N\|_{s,D}\,<\,\varepsilon_A/2,
\end{align*}
and the inner parameter bounds follow from the support $\Theta_R$.

The discretization in Lemma \ref{lem:quad} yields 
\begin{align*}
    \sum_{k=1}^N |c_k|\,\le\,\|\nu_R\|_{\rm TV}\,\le\,C_T\,\|v_R\|_{0,\Real^d}\,\le\,C_T\,\|Eu^*\|_{s+\varepsilon_C,\Real^d}\,\le\,C_T\,C_E\,\|u^*\|_{s+\varepsilon_C,D}:=C_N,
\end{align*}
which completes the proof.
\end{proof}

We rewrite $u_a$ as
\begin{align} \label{eq:rnnfunc}
    u_a(\bx)\,=\,\sum_{j=1}^M\alpha_j\rho(\tilde{\bw}_j^\top\bx+\tilde{b}_j)
\end{align}
Using the triangle inequality
\begin{align*}
    \|u^*-u_a\|_{s,D}\,\le\,\|u^*-u_N\|_{s,D}+\|u_N-u_a\|_{s,D},
\end{align*}
the subgoal for approximation error is to show
\begin{align*}
    \|u_N-u_a\|_{s,D}\,<\,\varepsilon_A/2\quad\text{with high probability.}
\end{align*}

To this end, we need the following lemma. This follows from the smoothness of $\rho$ and the continuity of the map
$(w,b)\mapsto \rho(w^\top x + b)$ as a function from the parameter box $\Theta$ into $H^s(D)$.

\begin{lemma}
    Under Assumption \ref{As:act}, for every $\theta_k=(\bw_k,b_k)$ in a bounded set $\Theta=[-r,r]^{d+1}$, there is a neighborhood $B(\theta_k,\delta_r)$ such that
    \begin{align*}
        \left\|\rho(\bw_k^\top\bx+b_k)-\rho(\tilde{\bw}_j^\top\bx+\tilde{b}_j)\right\|_{s,D}\,<\,\varepsilon_1\quad\forall\,\tilde{\theta}_j=(\tilde{\bw}_j,\tilde{b}_j)\in B(\theta_k,\delta_r).
    \end{align*}
    If $B(\theta_k,\delta_r)\cap\partial\Theta\neq\emptyset$, we can take $\Theta=[-r-\delta_r,r+\delta_r]^{d+1}$ to avoid the case, and always use $\Theta=[-r,r]^{d+1}$ for brevity.
\end{lemma}

Let $p_k$ denote the fraction of $\Theta$ occupied by the ball:
\begin{align*}
    p_k=\frac{\text{vol}(B(\theta_k,\delta_r))}{\text{vol}(\Theta)}.
\end{align*}
We have $p_1=\cdots=p_N:=p_0$. Our goal is to find some $\tilde{\theta}_j\in B(\theta_k,\delta_r)$ for all $\theta_k$ with high probability. We know that the probability of $\tilde{\theta}_j\notin B(\theta_k,\delta_r)$ for a certain $\theta_k$ is
\begin{align*}
    1-p_0.
\end{align*}
So the probability that all $M$ independent draws miss the $B(\theta_k,\delta_r)$ is
\begin{align*}
    (1-p_0)^M.
\end{align*}
For all $\theta_k$, according to Boole's inequality, the probability of the event ``at least one $B(\theta_k,\delta_r)$ does not contain some $\tilde{\theta}_j$'' is less than
\begin{align*}
    \sum_{k=1}^N\,(1-p_0)^M\,=\,N(1-p_0)^M.
\end{align*}
We require the probability of the event occurring to be sufficiently small  ($<\delta_p$), considering the opposing event, all $N$ distinct neurons $\theta_1,\cdots,\theta_N$ get matched at least one $\tilde{\theta}_j$ with probability $>1-\delta_p$. Taking
\begin{align} \label{ineq:Mcondition}
    M\,\ge\,\frac{1}{p_0}\,\ln\left(\frac{N}{\delta_p}\right)\,=\,\frac{\text{vol}(\Theta)}{\text{vol}(B(\theta_k,\delta_r))}\,\ln\left(\frac{N}{\delta_p}\right),
\end{align}
one can verify $N(1-p_0)^M<\delta_p$ with the standard inequality $(1-p)\le e^{-p}$. 

If each deterministic basis function $\rho(\bw_k^\top\bx+b_k)$ is matched by some randomized basis functions $\rho(\tilde{\bw}_j^\top\bx+\tilde{b}_j)$, we can pick $\alpha_j=c_k$ and $\alpha_j=0$ otherwise. If $\rho(\bw_s^\top\bx+b_s)$ and $\rho(\bw_t^\top\bx+b_t)$ are matched by the same $\rho(\tilde{\bw}_j^\top\bx+\tilde{b}_j)$, we can take $\alpha_j=c_s+c_t$. Summing up the basis functions, we obtain
\begin{align*}
    \|u_N-u_a\|_{s,D}\,=\,\left\|\sum_{k=1}^N\,c_k\rho(\bw_k^\top\bx+b_k)-\sum_{j=1}^M\,\alpha_j\rho(\tilde{\bw}_j^\top\bx+\tilde{b}_j)\right\|_{s,D}\,<\,C_N\,\varepsilon_1.
\end{align*}
Taking $\varepsilon_1=\varepsilon_A/2C_N$, we have
\begin{align*}
    \|u_N-u_a\|_{s,D}\,<\,\varepsilon_A/2\quad\text{with probability }1-\delta_p,
\end{align*}
where $u_a$ is \eqref{eq:rnnfunc} and $M$ satisfies \eqref{ineq:Mcondition}. 
\begin{theorem}[UAT of RaNN]\label{thm:uat_rann}
    Under Assumption \ref{As:act}, taking 
    \begin{align*}
        M\,\ge\,\frac{1}{p_0}\,\ln\left(\frac{N}{\delta_p}\right)=\frac{\text{vol}(\Theta)}{\text{vol}(B(\theta_k,\delta_r))}\,\ln\left(\frac{N}{\delta_p}\right),
    \end{align*}
    where $\Theta$, $\theta_k$ and $N$ are defined in Theorem \ref{Thm:H1bound-UAT}. With probability at least $1-\delta_p$, the randomized shallow neural network $u_a$ can approximate $u^*\in H^{s+\varepsilon_C}(D)$ to any desired accuracy $\varepsilon_A>0$:
    \begin{align*}
        \|u^*-u_a\|_{s,D}\,<\,\varepsilon_A.
    \end{align*}
\end{theorem}

\subsection{Statistical Component}
In this part, we only estimate the statistical term when $Y=L^2(D)$ and $Z=L^2(\Gamma)$, this is the most common loss for neural network methods. We begin with the necessary setup and assumptions.

\begin{lemma}[Matrix Bernstein Inequality, Theorem 1.4 in \cite{Tropp2012randmat}] \label{le:MBI}
    Let $X_1,\cdots,X_n$ be independent, random self-adjoint matrices of dimension $d$. Assume that each random matrix satisfies
    \begin{align*}
        \dE\left[X_i\right]\,=\,0\quad\text{and}\quad\|X_i\|_\text{op}\,\le\,J\quad\text{almost surely},
    \end{align*}
    where $\|\cdot\|_\text{op}$ is the operator (spectral) norm of a self-adjoint matrix. Then for all $t\ge0$,
    \begin{align*}
        \dP\left(\left\|\sum_{i=1}^n\,X_i\right\|_\text{op}\ge t\right)\,\le\,d\,\exp\left(\frac{-t^2/2}{\sigma^2+J\,t/3}\right), \quad\text{where }\sigma^2\,=\,\left\|\sum_{i=1}^n\dE\left[X_i^2\right]\right\|_\text{op}.
    \end{align*}
\end{lemma}

Inheriting the previous notation and results, each $u\in\cN_\rho(D)$ has the following form
\begin{align*}
    u(\bx)\,=\,\balpha\rho(\bW\bx+\bb)\,=\,\sum_{j=1}^M\,\alpha_j\phi_j(\bx)\quad\text{where }\phi_j(\bx)\,=\,\rho(\bw_j^\top \bx+b_j),
\end{align*}
and $\|\balpha\|_1\le C_N$, $C_N$ is defined in Theorem \ref{Thm:H1bound-UAT}. $\cN_\rho(D)$ is finite dimensional, each $\phi_j$ is bounded, then there is a constant $C_R>0$, which is independent of $M$, such that 
\begin{align*}
    |\cG u(\bx)|\,\le\,C_R,\quad|f(\bx)|\,\le\,C_R,\quad\forall\,\bx\in D,\,\forall\,u\in\cN_\rho(D),
\end{align*}
and 
\begin{align*}
    |\cB u(\bx)|\,\le\,C_R,\quad|g(\bx)|\,\le\,C_R,\quad\forall\,\bx\in\Gamma,\,\forall\,u\in\cN_\rho(D).
\end{align*}

\begin{theorem}[Statistical Component Bound] \label{Thm:statisticerr}
    Under the above conditions, with probability at least $1-\delta_s$ for the random choice of collocation points, the statistical component has bound
    \begin{align*}
        \sup_{u\in\cN_\rho(D)}\left|\cL(u)-\widehat{\cL}(u)\right|\,\lesssim\,C_N^2C_R^2M\left(\sqrt{\frac{\log(6M/\delta_s)}{N^I}}+\sqrt{\frac{\log(6M/\delta_s)}{N^B}}\right),
    \end{align*}
    where $\widehat{\cL}(u)$ is the empirical $L^2$ loss defined in \eqref{eq:emp-loss-L2}.
\end{theorem}

\begin{proof}
We begin by splitting the statistical term into domain and boundary components, 
\begin{align*}
    \cL(u)-\widehat{\cL}(u)\,=\,\left[\cL_D(u)-\widehat{\cL}_D(u)\right]+\left[\cL_\Gamma(u)-\widehat{\cL}_\Gamma(u)\right].
\end{align*}
By the triangle inequality,
\begin{align*}
    \sup_{u\in\cN_\rho(D)}\left|\cL(u)-\widehat{\cL}(u)\right|\,\le\,\sup_{u\in\cN_\rho(D)}\left|\cL_D(u)-\widehat{\cL}_D(u)\right|+\sup_{u\in\cN_\rho(D)}\left|\cL_\Gamma(u)-\widehat{\cL}_\Gamma(u)\right|.
\end{align*}
We only give the argument for the domain component, the boundary term can be done in the same way.

Next, we express the domain component in terms of the expectation and the sample mean:
\begin{align*}
    \cL_D(u)\,=\,\int_D\,\left|\cG u(\bx)-f(\bx)\right|^2\,\dbx\,=\,|D|\,\dE\left[\left|\cG u(X)-f(X)\right|^2\right],
\end{align*}
where $X\sim\cU(D)$, and 
\begin{align*}
    \widehat{\cL}_D(u)\,=\,\frac{|D|}{N^I}\,\sum_{i=1}^{N^I}\,\left|\cG u(\bx^I_i)-f(\bx^I_i)\right|^2,
\end{align*}
where $\bx_i^I$ are independent and identically distributed uniform sampling points. Hence
\begin{align} \label{eq:statistic1}
    \sup_{u\in\cN_\rho(D)}\left|\cL_D(u)-\widehat{\cL}_D(u)\right|\,=\,|D|\sup_{u\in\cN_\rho(D)}\left|\dE\left[\left|\cG u(X)-f(X)\right|^2\right]-\frac{1}{N^I}\,\sum_{i=1}^{N^I}\left|\cG u(\bx^I_i)-f(\bx^I_i)\right|^2\right|.
\end{align}

Since $\cG$ is a linear operator, each $u$ is a linear combination of a finite set of basis functions $\phi_j$, letting
\begin{align*}
    \bv(\bx)\,=\,(\cG\phi_1(\bx),\cdots,\cG\phi_M(\bx))^\top,
\end{align*}
one obtains
\begin{align*}
    \left|\cG u(\bx)-f(\bx)\right|^2\,=\,\left(\balpha\bv(\bx)-f(\bx)\right)^2=\balpha\bv(\bx)\,\bv(\bx)^\top\balpha^\top-2\balpha\bv(\bx)\,f(\bx)+\left(f(\bx)\right)^2.
\end{align*}
Define
\begin{align*}
    A\,=\,\dE\left[\bv(X)\,\bv(X)^\top\right],\quad\widehat{A}\,=\,\frac{1}{N^I}\sum_{i=1}^{N^I}\,\bv(\bx_i^I)\,\bv(\bx_i^I)^\top,
\end{align*}
and similarly
\begin{align*}
    b\,=\,\dE\left[\bv(X)f(X)\right],\quad\widehat{b}\,=\,\frac{1}{N^I}\,\sum_{i=1}^{N^I}\,\bv(\bx_i^I)\,f(\bx_i^I),
\end{align*}
then 
\begin{align} \label{eq:statistic2}
    &\sup_{u\in\cN_\rho(D)}\left|\dE\left[\left|\cG u(X)-f(X)\right|^2\right]-\frac{1}{N^I}\,\sum_{i=1}^{N^I}\,\left|\cG u(\bx^I_i)-f(\bx^I_i)\right|^2\right|\nonumber\\
    =\,&\sup_{\|\balpha\|_1\le C_N}\left|\balpha(A-\widehat{A})\balpha^\top-2\balpha(b-\widehat{b})+C_f\right|\,\le\,C_N^2\,\|A-\widehat{A}\|_\text{op}+2C_N\,\|b-\widehat{b}\|_2+C_f,
\end{align}
where 
\begin{align*}
    C_f\,=\left|\,\dE\,[f(X)^2]-\frac{1}{N^I}\,\sum_{i=1}^{N^I}\,f(\bx^I_i)^2\right|.
\end{align*}

To apply the Matrix Bernstein Inequality, we define
\begin{align*}
    Y_i\,=\,\bv(\bx_i^I)\,\bv(\bx_i^I)^\top-A,\quad\text{and} \quad \frac{1}{N^I}\,\sum_{i=1}^{N^I}\,Y_i\,=\,\widehat{A}-A.
\end{align*}
Here, each $Y_i$ is an $M\times M$ symmetric matrix and $\dE\,[Y_i]=0$. Since $|\cG u(\bx)|\le C_R$ for all $\bx\in D$ and $u\in\cN_\rho(D)$, we have
\begin{align*}
    \left\|\bv(\bx)\,\bv(\bx)^\top\right\|_\text{op}\,=\,\left\|\bv(\bx)\right\|_2^2\,\le\,C_R^2M\quad\text{and}\quad\|A\|_\text{op}\,\le\,\dE\,\left[\|\bv(X)\|_2^2\right]\,\le\,C_R^2M.
\end{align*}
Hence, 
\begin{align*}
    \|Y_i\|_\text{op}\,=\,\left\|\bv(\bx_i^I)\,\bv(\bx_i^I)^\top-A\right\|_\text{op}\,\le\,\left\|\bv(\bx_i^I)\,\bv(\bx_i^I)^\top\right\|_\text{op}+\|A\|_\text{op}\,\le\,C_R^2M+C_R^2M=2C_R^2M.
\end{align*}

Recalling Lemma \ref{le:MBI}, in our setting, $X_i=Y_i$, $n=N^I$, $d=M$, $J=2C_R^2M$ and
\begin{align*}
    \sigma^2\,=\,\left\|\sum_{i=1}^{N^I}\,\dE\left[Y_i^2\right]\right\|_\text{op}\,\le\,\sum_{i=1}^{N^I}\,\dE\left[\|Y_i^2\|_\text{op}\right]\,\le\,\sum_{i=1}^{N^I}\,\dE\left[\|Y_i\|_\text{op}^2\right]\,\le\,N^I(2C_R^2M)^2\,=\,4N^IC_R^4M^2.
\end{align*}
Apply the inequality with $\sum_{i=1}^{N^I}Y_i$, we get, for any $t>0$,
\begin{align*}
    \dP\,\left(\left\|\sum_{i=1}^{N^I}\,Y_i\right\|_\text{op}\ge t\right)\,\le\,M\,\exp\left(\frac{-t^2/2}{\sigma^2+Jt/3}\right).
\end{align*}
Since
\begin{align*}
    \left\|A-\widehat{A}\right\|_\text{op}\,=\,\frac{1}{N^I}\,\left\|\sum_{i=1}^{N^I}\,Y_i\right\|_\text{op},
\end{align*}
we have
\begin{align*}
    \dP\left(\left\|A-\widehat{A}\right\|_\text{op}\ge \frac{t}{N^I}\right)\,\le\,M\,\exp\left(\frac{-t^2/2}{\sigma^2+Jt/3}\right),\,\quad\forall\,t>0
\end{align*}
with $J=2C_R^2M$ and $\sigma^2\le4N^IC_R^4M^2$.

Take
\begin{align*}
    \frac{t}{N^I}\,>\,MC_R^2\left(\frac{2}{3}t_0+\sqrt{\frac{4}{9}t_0^2+16t_0}\right),\quad\text{where } t_0\,=\,\frac{\log(6M/\delta_s)}{N^I},
\end{align*}
we get 
\begin{align*}
    M\exp\left(\frac{-t^2/2}{\sigma^2+Jt/3}\right)\,<\,\frac{\delta_s}{6}.
\end{align*}
Therefore, we conclude that with probability $\ge1-\delta_s/6$,
\begin{align} \label{ineq:statA}
    \left\|A-\widehat{A}\right\|_\text{op}&\,<\,C_R^2M\left(\frac{2}{3}\,\frac{\log(6M/\delta_s)}{N^I}+\sqrt{\frac{4}{9}\left(\frac{\log(6M/\delta_s)}{N^I}\right)^2+16\,\frac{\log(6M/\delta_s)}{N^I}}\right)\nonumber\\
    &\lesssim C_R^2M\,\sqrt{\frac{\log(6M/\delta_s)}{N^I}},\quad\text{for large }N^I.
\end{align}
Similarly, a vector Bernstein inequality yields
\begin{align} \label{ineq:statb}
    \left\|b-\widehat{b}\right\|_2\,\lesssim\,C_R\,\sqrt{M}\,\sqrt{\frac{\log(6M/\delta_s)}{N^I}},\quad\text{for large }N^I.
\end{align}
Since
\begin{align*}
    |f(\bx)|\,\le\,C_R,\quad\forall\,\bx\in D,
\end{align*}
applying Hoeffding inequality obtains
\begin{align*}
    \dP\left(C_f\ge s\right)\,\le\,\exp\left(-\frac{2N^Is^2}{C_R^4}\right),\quad\forall\,s>0,
\end{align*}
which implies
\begin{align} \label{ineq:statf}
    C_f\,\le\,C_R^2\sqrt{\frac{\log(6/\delta_s)}{2N^I}},\quad\text{with probability }\ge1-\delta_s/6.
\end{align}
With \eqref{eq:statistic1} and \eqref{eq:statistic2}, we get with probability $\ge1-\delta_s/2$
\begin{align*}
    \sup_{u\in\cN_\rho(D)}\left|\cL_D(u)-\widehat{\cL}_D(u)\right|&\,\lesssim\,C_N^2C_R^2M\sqrt{\frac{\log(6M/\delta_s)}{N^I}}+2C_NC_R\sqrt{M}\sqrt{\frac{\log(6M/\delta_s)}{N^I}}+C_R^2\sqrt{\frac{\log(6/\delta_s)}{2N^I}}\\
    &\lesssim C_N^2C_R^2M\sqrt{\frac{\log(6M/\delta_s)}{N^I}},\quad\text{for large }N^I,
\end{align*}
where $1-\delta_s/2=1-3(1-(1-\delta_s/6))$ comes from the non-independence of the inequalities \eqref{ineq:statA}, \eqref{ineq:statb} and \eqref{ineq:statf}.

The same reasoning applies to $\cL_\Gamma(u)-\widehat{\cL}_\Gamma(u)$, the result is
\begin{align*}
    \sup_{u\in\cN_\rho(D)}\left|\cL_\Gamma(u)-\widehat{\cL}_\Gamma(u)\right|\,\lesssim\,C_N^2C_R^2M\sqrt{\frac{\log(6M/\delta_s)}{N^B}},\quad\text{with probability }\ge1-\delta_s/2\text{ for large }N^B.
\end{align*}

Finally, by combining domain and boundary components, we obtain the desired bound (with probability at least $(1-\delta_s/2)^2=1-\delta_s+\delta_s^2/4>1-\delta_s$, since the interior and boundary collocation points are independence):
\begin{align*}
    \sup_{u\in\cN_\rho(D)}\left|\cL(u)-\widehat{\cL}(u)\right|\,\lesssim\,C_N^2C_R^2M\left(\sqrt{\frac{\log(6M/\delta_s)}{N^I}}+\sqrt{\frac{\log(6M/\delta_s)}{N^B}}\right),
\end{align*}
which concludes the proof.
\end{proof}

\subsection{A convergence result}

Combining the approximation result of Theorem \ref{thm:uat_rann} with the statistical bound in Theorem \ref{Thm:statisticerr} and the graph norm equivalence Assumption \ref{ass:graph}, we obtain the following convergence theorem for RaNN-type PDE solvers.

\begin{theorem}[Convergence of RaNN-based PDE solvers]\label{thm:RaNN-convergence}
Suppose that Assumption \ref{ass:graph} holds, that $V=H^s(D)$ with $u^*\in H^{s+\varepsilon_C}(D)$ for some $\varepsilon_C>0$, and that the activation function $\rho$ satisfies Assumption \ref{As:act}. Assume further that the inner parameters $(W,b)$ of the RaNN are sampled from a fixed bounded box $\Theta\subset\mathbb R^{d+1}$, and that the output-layer coefficients are obtained as an exact minimizer of the empirical loss \eqref{eq:emp-loss-L2} based on i.i.d.\ interior and boundary collocation points.

Then, for every tolerance $\varepsilon>0$ and confidence level $\delta\in(0,1)$, there exist a network width $M$ and sample sizes $N^I,N^B\in\mathbb N$ such that the corresponding empirical minimizer $u_\rho\in\mathcal N_\rho(D)$ satisfies
\begin{align*}
    \|u_\rho-u^*\|_{s,D} \;\le\; \varepsilon
\end{align*}
with probability at least $1-\delta$ with respect to the randomness of both the collocation points and the inner parameters.
\end{theorem}

The present analysis focuses on single-hidden-layer RaNNs. The adaptive layer growth in AG-RaNN only enlarges the hypothesis space, so the approximation error bounds remain valid in the sense that deeper architectures can achieve at least the same accuracy; a rigorous analysis for fully multi-layer randomized networks is left for future work.

\section{Numerical Experiments} \label{Sec:NumerExper}

In this section, we present numerical experiments to demonstrate the performance of AG-RaNN using the different strategies introduced earlier.

All experiments are conducted in MATLAB with the following settings. To ensure reproducibility, we use the random number generator ``rng(1)''. The QR method is employed to solve the least-squares problem using ``linsolve'' with the option ``opts.RECT = true''. The discrete Fourier transform (DFT) is implemented using the ``fft'' function, with $10^4$ sampling points for 1-D problems and $100\times100$ sampling points for 2-D problems.

We adopt the hyperbolic tangent function $\text{tanh}(x)=\frac{e^x-e^{-x}}{e^x+e^{-x}}$ or the Gaussian function $G(x)=e^{-\frac{x^2}{2}}$ as the activation function $\rho(x)$. Accuracy is measured using the relative $L^2$ error, defined as $e_0(u_\rho)=\|u_\rho-u\|_{0,D}/\|u\|_{0,D}$. We approximate the $L^2$ norms by Gauss–Legendre quadrature with $200$ points for 1-D problems, $100\times100$ points for 2-D problems, and $40^3$ points for 3-D problems. 

The main parameters used in the experiments are presented below:
\begin{itemize}
\item $N^I$: number of interior points, and $N^B$: number of boundary points. 

\item $m_l$: number of neurons in the $l$-th hidden layer. $\bmm_\text{add}$ is the number of neurons added during growth, e.g., $\bmm_\text{add}=(200,100,50,50)$ means $200$ initial neurons and additional growth of $100$, $50$ and $50$ neurons over $3$ stages.

\item $r_\text{max}$ and $\Lambda$ are parameters used to generate the set of candidate functions in frequency-based parameter initialization.

\item $\br=(r_1,\cdots,r_d)$ denotes the parameters of the uniform distribution $\cU(-\br,\br)$, where $r_i$ is the parameter for the $i$-th dimension. $\br^\text{opt}$ is calculated during initialization and neuron growth. $\br_1$ and $\br_2$ are used to generate $\bW^{(0)}\sim\cU(-\br_1,\br_1)$ and $\bH_0\sim\cU(-\br_2,\br_2)$ during layer growth. 
    
\item $\eta$: penalty parameter for boundary conditions.

\item $(It_P,It_N)$: number of Picard and Newton iterations.
\end{itemize}

\subsection{Collocation Points} \label{Subsec:Point}

Next, we discuss the selection of collocation points. In this work, uniform sampling is used. For instance, in a 2-D domain $D=(p_1,q_1)\times(p_2,q_2)$, $N^I=N_1^IN_2^I$ interior points are generated as follows:
\begin{align*}
    x_{1k}&\,=\,p_1+\varepsilon_c+(k-1)\,\frac{q_1-p_1-2\varepsilon_c}{N_1^I-1}\quad \text{for }k=1,2,\cdots,N_1^I,\\
    x_{2l}&\,=\,p_2+\varepsilon_c+(l-1)\,\frac{q_2-p_2-2\varepsilon_c}{N_2^I-1}\quad \text{for }l=1,2,\cdots,N_2^I.
\end{align*}
Here, $\varepsilon_c$ is a small constant (set to $10^{-10}$ in all the examples) to avoid sampling too close to the boundary. Let $i=(k-1)\,N_2^I+l$ and $\bx_i=(x_{1k},x_{2l})$ represents the interior points for $i=1,\cdots,N^I$, and $y_i=f(\bx_i)$ from the PDE.

If $\Gamma=\partial D$, a similar sampling approach is applied. Sampling on each edge is done separately, and the total number of boundary points is $N^B=N^B_1+N^B_2+N^B_3+N^B_4$, where $N^B_i$ may vary. For example, considering $\Gamma_1=[p_1,q_1]\times\{p_2\}$:
\begin{align*}
    x_j\,=\,p_1+(j-1)\frac{q_1-p_1}{N_1^B}\quad \text{for }j=1,2,\cdots,N_1^B.
\end{align*}
The collocation points on $\Gamma_1$ are $\set{\bx_j=(x_j,p_2)}_{j=1}^{N_1^B}$, and $y_j=g(\bx_j)$. Considering all edges, the data points comprise $(\bx_j,y_j)_{j=1}^{N^B}$. 

\subsection{Numerical Examples}

\begin{example}[Poisson Equation with Oscillatory or Sharp Solution] \label{ex:Poisson}
We consider the Poisson equation:
\begin{align}
    -\Delta u &\,=\,f\quad \text{\rm in } D, \\
    u &\,=\,g\quad \text{\rm on } \partial D.
\end{align}
We investigate five cases with varying spatial dimensions $d$ and solution characteristics. 
\begin{description}
\item[Case 1:] Demonstrates the effectiveness of frequency-based parameter initialization and neuron growth for a 1-D problem.

\item[Case 2:] Examines the impact of different design choices on layer growth for a 2-D problem.

\item[Case 3:] Combines all three strategies (parameter initialization, neuron growth, and layer growth) for a 2-D problem.

\item[Case 4:] Addresses a 3-D problem using layer growth.

\item[Case 5:] Solves an equation with multiple peaks through layer growth.
\end{description}

For Case 1, the exact solution is:
\begin{align*}
    u(x)\,=\,\frac{1}{6}\,\sum_{s=1}^6\,\sin(2^s\pi x).
\end{align*}
For Cases 2, 3, and 4, the exact solution is:
\begin{align*}
    u(\bx)\,=\,4^d\left(\frac{1}{2}+\frac{1}{\pi}\arctan\left(A\left(\frac{1}{16}-\sum_{i=1}^d\left(x_i-\frac{1}{2}\right)^2\right)\right)\right)\prod_{i=1}^d\,(1-x_i)x_i,
\end{align*}
where $A = 120$. 
For Case 5, the exact solution with $P_1$ peaks is given by
\begin{align*}
    u(x_1, x_2)\,=\,\sum_{p=1}^{P_1}\,\exp\left(-1000\left(\left(x_1-\frac{\sqrt{2}}{2}\cos\left(\frac{2\pi p}{P_1}+\frac{\pi}{4}\right)\right)^2+\left(x_2-\frac{\sqrt{2}}{2}\sin\left(\frac{2\pi p}{P_1}+\frac{\pi}{4}\right)\right)^2\right)\right).
\end{align*}
And the exact solution with $P_2^2$ peaks is defined as
\begin{align*}
    u(x_1, x_2)\,=\,\sum_{p_1=1}^{P_2}\sum_{p_2=1}^{P_2}\exp\left(-1000\left(\left(x_1-\frac{2p_1-1}{2P_2}\right)^2+\left(x_2-\frac{2p_2-1}{2P_2}\right)^2\right)\right).
\end{align*}
When $P_1=1,2,4$, these cases correspond to those considered in  \cite{Zhou2023Sampling,huang2024adaptiveneuralnetworkbasis}. The domain is $D=(0,1)^d$ for the first four cases, and $D=(-1,1)^2$ for the last one. The parameters for each case are shown in Table \ref{tab:ex_Poisson_parameter}.
\end{example}

\begin{table}[!htbp]
    \centering
    \begin{tabular}{|c|c|c|c|c|c|c|c|c|c|}
    \hline
     & $d$ & $\rho_1(x)$ & $\rho_2(x)$ & $N^I$ & $N^B$ & $\eta$ \\ \hline
    \textbf{Case 1} & 1 & $G(x)$ & N/A & 2000 & 2 & 1 \\ \hline    
    \textbf{Case 2} & 2 & $\tanh(x)$ & $G(x)$ & $300\times300$ & 12000 & 3000 \\ \hline
    \textbf{Case 3} & 2 & $\tanh(x)$ & $G(x)$ & $300\times300$ & 12000 & 3000 \\ \hline
    \textbf{Case 4} & 3 & $\tanh(x)$ & $G(x)$ & $80\times80\times80$ & 9600 & 5000 \\ \hline
    \textbf{Case 5} & 2 & $\tanh(x)$ & $G(x)$ & $300\times300$ & 12000 & 1000 \\ \hline
    \end{tabular}
    \caption{Parameters used in Example \ref{ex:Poisson}. }
    \label{tab:ex_Poisson_parameter}
\end{table}

\noindent\textbf{Case 1: $d=1$, Initialization + Neuron Growth}

\begin{figure}[!htbp]
    \centering
    \includegraphics[width=\textwidth]{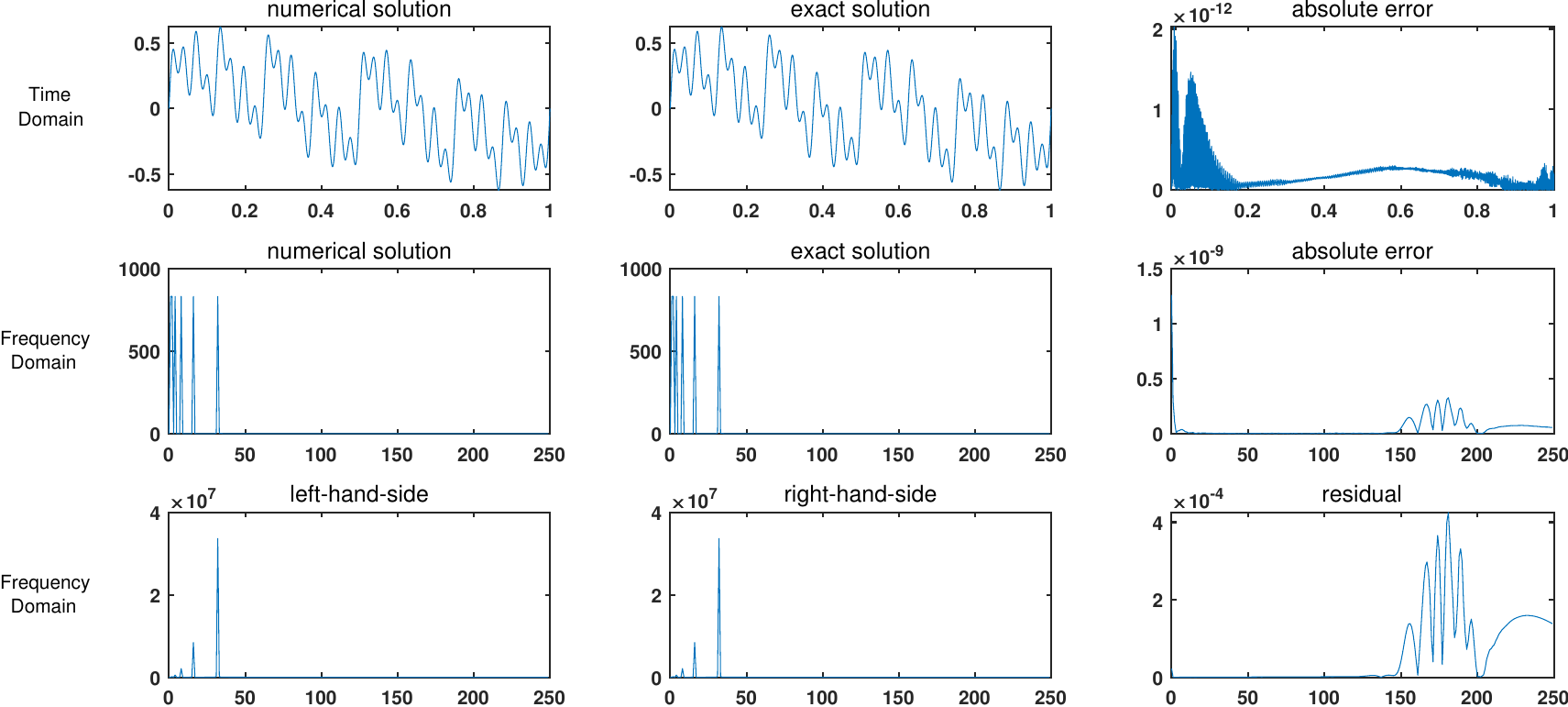}
    \caption{Best result obtained by the manually tuned RaNN in Case 1 of Example \ref{ex:Poisson}.}
    \label{fig:ex1_init-intra-1D_best}
\end{figure}

Using $r=200$ for the uniform distribution $\cU(-r,r)$ and $m_1=500$, the RaNN method achieves optimal results, as shown in Figure \ref{fig:ex1_init-intra-1D_best}, with a relative $L^2$ error of $1.07e-12$. It is noted that $r=200$ is manually selected to obtain the best performance under the uniform distribution. Parameters $N^I$, $N^B$ and $\eta$ are adjusted to optimize results for $m_1$, the number of neurons in the hidden layer of RaNN. The same parameters are used during initialization and neuron growth. During neuron growth, the sum of added neurons $\bmm_\text{add}$ matches $m_1$, ensuring a fair comparison between RaNN and AG-RaNN.

We set $r_\text{max}=800$, $\Lambda=100$ and $\bmm_\text{add}=(200,100,50,50,50,50)$, leading to $r_k=\frac{k}{\Lambda}r_\text{max}=8k$ for $k=1,\cdots,100$. Our aim is to find the near-optimal parameter $r^\text{opt}$ for the uniform distribution $\cU(-r^\text{opt},r^\text{opt})$. Set weights $\tilde{w}^j=r_j$ for $j=1,\cdots,100$, and define candidate functions as $\phi_j(x)=\rho(\tilde{w}^j(x-0.5))$, where $\rho(x)=e^{-x^2/2}$. 

\begin{figure}[!htbp]
    \begin{subfigure}{.65\textwidth}
        \centering
        \includegraphics[width=\linewidth]{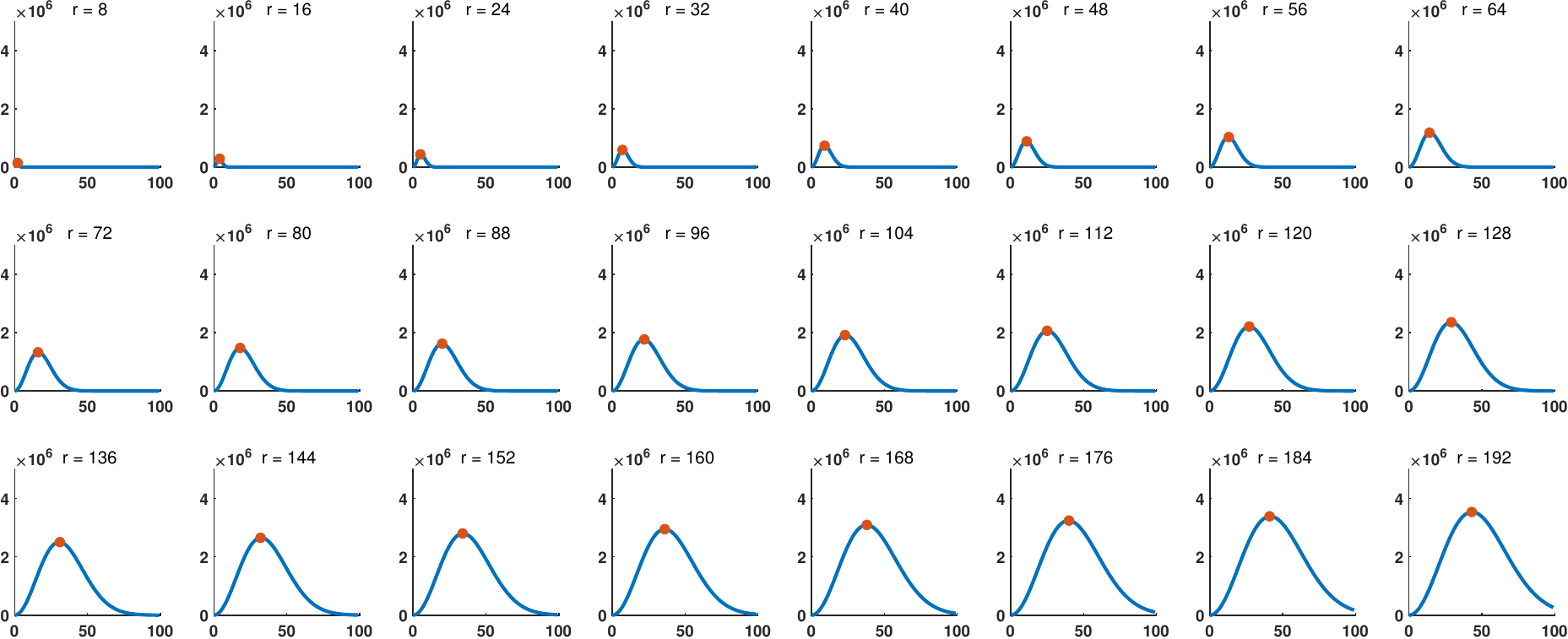}
        \caption{candidate basis functions $\left|\widehat{-\Delta \phi_j}\right|$}
        \label{fig:candidate-function}
    \end{subfigure}
    \hfil
    \begin{subfigure}{.28\textwidth}
        \centering
        \includegraphics[width=\linewidth]{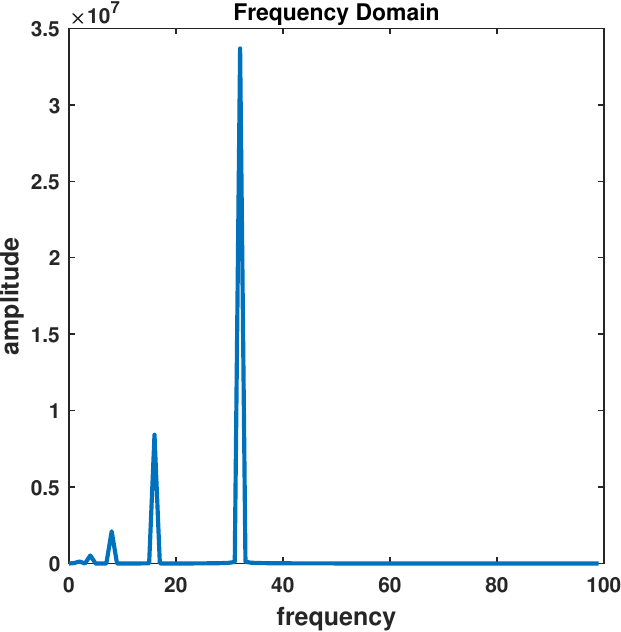}
        \caption{right-hand term $\left|\widehat{f}\right|$}
        \label{fig:Fourier-right}
    \end{subfigure}
    \caption{The discrete Fourier transform of different functions in Case 1 of Example \ref{ex:Poisson}.}
\end{figure}

Collocation points are sampled as detailed in Section \ref{Subsec:Point}, resulting in $\set{x_i}_{i=1}^{N^I}$ and $y_i=f(x_i)$. The discrete Fourier transform of $\set{f(x_i)}_{i=1}^{N^I}$ is computed. Similarly, for each $\phi_j$, we calculate the DFT of $\set{-\Delta\phi_j(x_i)}_{i=1}^{N^I}$, as illustrated in Figure \ref{fig:candidate-function}. The DFT of $f$ is shown in Figure \ref{fig:Fourier-right}.

In Figure \ref{fig:Fourier-right}, we observe $\xi_0=\mathop{\arg\max}\limits_{\xi}\left|\hat{y}(\xi)\right|=32$. From Figure \ref{fig:candidate-function}, we determine the first orange point greater than $32$, denoted as $j_0=\mathop{\arg\min}\limits_{j}\set{|\xi_j-\xi_0|:\xi_j>\xi_0}$, where $\xi_j$ represents the top orange point for $\left|-\widehat{\Delta \phi_j}\right|$. From Figure \ref{fig:candidate-function}, we find $j_0=19$ and $r^\text{opt}=\tilde{w}^{19}=152$. Therefore, the initial distribution is $\cU(-152,152)$.

The process can be explained intuitively: to obtain a good approximation, we select basis functions $\{\phi_j\}$ such that their transformed versions $\set{-\Delta \phi_j}$ closely match the highest-energy frequency of $f$. This allows us to determine the near-optimal parameter.

\begin{figure}[!htbp]
    \centering
    \includegraphics[width=\textwidth]{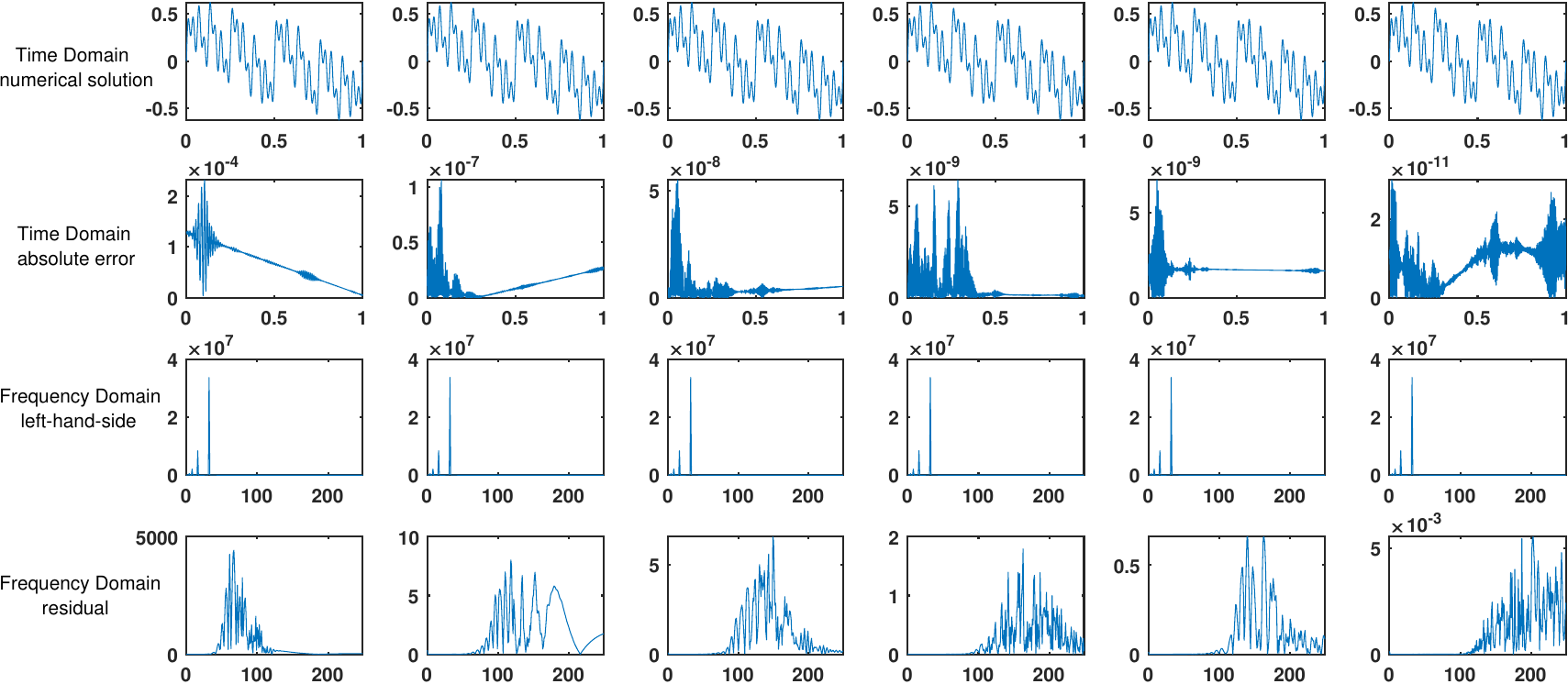}
    \caption{Results obtained by the AG-RaNN method in Case 1 of Example \ref{ex:Poisson}. The first column shows the result obtained by initialization, while the other columns show results obtained by the neuron growth strategy.}
    \label{fig:ex1_init-intra-1D}
\end{figure}

\begin{table}[!htbp]
    \centering
    \begin{tabular}{|c||c|c|c|c|c|c|c|c|c|}
    \hline
        $r^\text{opt}$ & 152 & 304 & 528 & 672 & 728 & 632 \\ \hline
        $e_0(u_\rho)$ & 2.68e-04 & 6.41e-08 & 2.85e-08 & 4.37e-09 & 5.89e-09 & 3.45e-11 \\ \hline
    \end{tabular}
    \caption{Results of the AG-RaNN with neuron growth in Case 1 of Example \ref{ex:Poisson}.}
    \label{tab:ex1_init-intra-1D}
\end{table}

For the AG-RaNN method, the numerical results are shown in Figure \ref{fig:ex1_init-intra-1D} and Table \ref{tab:ex1_init-intra-1D}. The table shows that the neuron growth strategy improves accuracy step by step, with the AG-RaNN ultimately reaching a relative $L^2$ error of 3.45e-11, nearly as accurate as the best result obtained by the manually tuned RaNN.

\noindent\textbf{Case 2: $d=2$, Layer Growth}

For AG-RaNN, we choose $\br_1=(15,15)$ for the uniform distribution $\cU(-\br_1,\br_1)$ and set the initial solution as $u_\rho^0$. To explore the influence of different factors, we obtain a result $u_\rho^4$ by using basis functions $\bpsi(\bx)$ (Section \ref{Sec:layer-grow}) in the growth layer, training both hidden layers, and generating $\bH_0$ randomly with $\br_2=(10,10)$. 

To evaluate the influence of different factors, we compare local basis functions $\bpsi(\bx)$ with nonlocal basis functions $\hat{\bpsi}(\bx)$ (without applying \eqref{eq:localization}), consider fixing the trainable parameters in the first hidden layer (represented by dotted red lines in Figure \ref{fig:RNN-localization}) versus training both layers, and examine two approaches for generating $\bH_0$: one based on the gradient $\nabla u_\rho^0$ and the other using random initialization. We obtain $u_\rho^1$ by using $\hat{\bpsi}(\bx)$, $u_\rho^2$ by training only the second hidden layer, and $u_\rho^3$ by generating $\bH_0$ using the gradient $\nabla u_\rho^0$. Since the random number generator is fixed as \texttt{rng(1)}, the value of $e_0(u^0_\rho)$ remains identical when $m_1$ is fixed and the same basis functions are used for $u_\rho^2$ and $u_\rho^4$. The corresponding results and running times are summarized in Table \ref{table:ex_Poisson_case2-rnn}, which shows that $u_\rho^4$ achieves the best overall accuracy. Choosing $\psi_i$ rather than $\tilde{\psi}_i$ as basis functions is deliberate: although $\tilde{\psi}_i$ may enhance local approximation, its global nature can introduce interference, whereas localization leads to a more stable and accurate solution. The relative $L^2$ error of $u_\rho^0$ is $5.59e-03$ as $m_1=2000$. With the same degrees of freedom, the relative $L^2$ error of $u_\rho^3$ is $9.66e-05$ for $m_1=1000$ and $m_2=1000$. These results indicate that increasing the number of neurons in the existing hidden layer can improve accuracy; however, the layer-growth strategy is more efficient and requires no additional hyperparameter tuning for computing $u_\rho^3$.

\begin{table}[!htbp]
\centering
\resizebox{\textwidth}{!}{
\begin{tabular}{|c|c|c|c|c|c|c|c|c|c|c|c|}
\hline
$m_1$ & $e_0(u_\rho^0)$ & $m_p$ & $m_2$ & $e_0(u_\rho^1)$ & time(s) & $e_0(u_\rho^2)$ & time(s) & $e_0(u_\rho^3)$ & time(s) & $e_0(u_\rho^4)$ & time(s) \\ \hline\hline
\multirow{3}{*}{1000} & \multirow{3}{*}{1.49e-02} & \multirow{3}{*}{0} & 500 & 1.05e-02 & 7.87 & 1.69e-03 & 5.00 & 3.21e-04 & 9.52 & 8.03e-05 & 11.04 \\ \cline{4-12} 
 &  &  & 1000 & 1.02e-02 & 8.43 & 1.55e-04 & 6.85 & 9.66e-05 & 12.88 & 2.00e-05 & 11.98 \\ \cline{4-12} 
 &  &  & 1500 & 1.01e-02 & 11.78 & 4.00e-05 & 14.88 & 4.09e-05 & 18.34 & 7.84e-06 & 13.00 \\ \hline\hline
\multirow{3}{*}{1500} & \multirow{3}{*}{8.39e-03} & \multirow{3}{*}{8} & 500 & 6.97e-03 & 12.07 & 1.48e-03 & 11.37 & 1.83e-04 & 14.06 & 2.86e-05 & 13.95 \\ \cline{4-12} 
 &  &  & 1000 & 6.88e-03 & 14.71 & 7.77e-05 & 13.93 & 4.76e-05 & 16.83 & 5.90e-06 & 15.86 \\ \cline{4-12} 
 &  &  & 1500 & 6.80e-03 & 22.94 & 1.73e-05 & 16.87 & 1.68e-05 & 19.94 & 1.99e-06 & 19.59 \\ \hline\hline
\multirow{3}{*}{2000} & \multirow{3}{*}{5.59e-03} & \multirow{3}{*}{47} & 500 & 5.02e-03 & 18.10 & 9.91e-04 & 14.35 & 1.25e-04 & 17.26 & 1.49e-05 & 15.27 \\ \cline{4-12} 
 &  &  & 1000 & 4.93e-03 & 20.04 & 5.28e-05 & 15.44 & 3.13e-05 & 20.78 & 2.33e-06 & 20.94 \\ \cline{4-12} 
 &  &  & 1500 & 4.92e-03 & 24.32 & 1.10e-05 & 18.48 & 1.37e-05 & 24.46 & 8.77e-07 & 22.85 \\ \hline
\end{tabular}
}
\caption{Relative $L^2$ errors and running time of AG-RaNN methods in Case 2 of Example \ref{ex:Poisson}. }
\label{table:ex_Poisson_case2-rnn}
\end{table}

\begin{figure}
    \centering
    \includegraphics[width=0.9\linewidth]{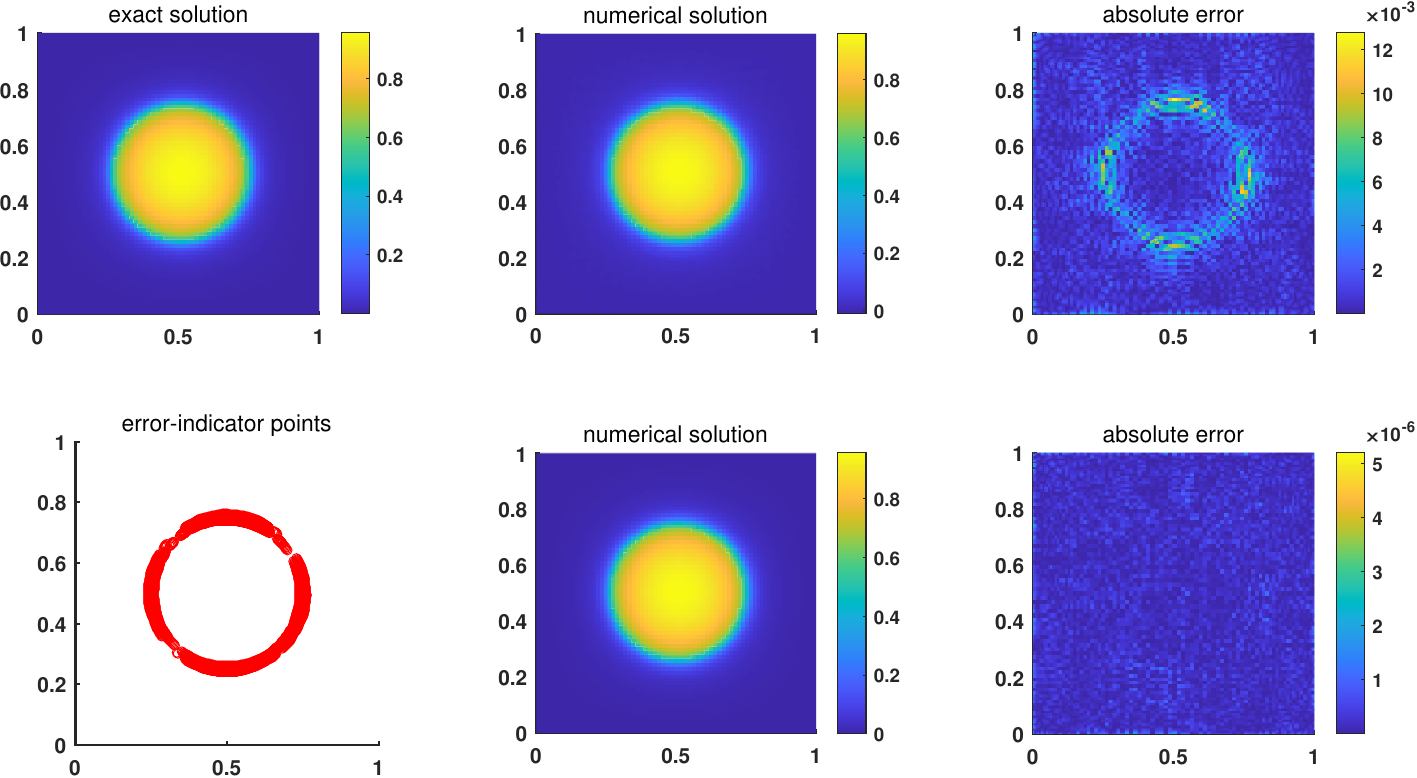}
    \caption{Exact solution, numerical solutions (middle), absolute errors (right) and error-indicator points $\bX_\text{err}$ obtained using the layer growth strategy for $m_1=2000$ and $m_2=1500$ in Case 2 of Example \ref{ex:Poisson}. }
    \label{fig:ex_Poisson_case2_layer}
\end{figure}

\begin{figure}[!htbp]
    \centering
    \includegraphics[width=0.98\textwidth]{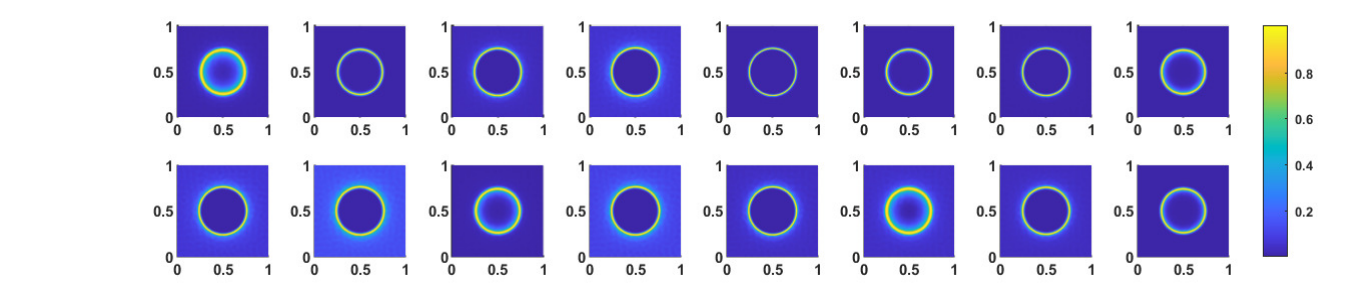}
    \includegraphics[width=0.98\textwidth]{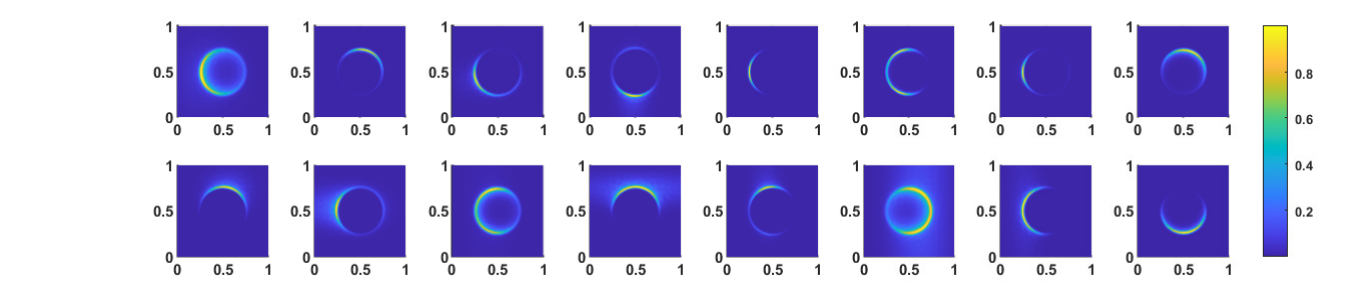}
    \caption{Examples of basis functions $\tilde{\psi}_i$ (top) and corresponding localized basis functions $\psi_i$ (bottom) in Case 2 of Example \ref{ex:Poisson}.}
    \label{fig:layergrow_basis}
\end{figure}

When $m_1=2000$ and $m_2=1500$, the figures of the exact solution $u(\bx)$, the error-indicator points $\bX_\text{err}$, as well as the numerical solutions and errors of $u_\rho^0$ and $u_\rho^4$ are presented in Figure \ref{fig:ex_Poisson_case2_layer}. Figure \ref{fig:layergrow_basis} showcases examples of $\set{\tilde{\psi}_j(\bx)}_{j=1}^{m_2}$ and $\set{\psi_j(\bx)}_{j=1}^{m_2}$. These notations are used consistently in subsequent examples.

\begin{table}[!htbp]
\centering
\begin{tabular}{|c|c|c|c|c|c||c|c|c|c|}
\hline
\multicolumn{2}{|c|}{$P_1$} & \multicolumn{2}{c|}{$P_2$} & \multicolumn{2}{c||}{$P_3$} & \multicolumn{2}{c|}{$P_1^a$} \\ 
\hline
DoF & $e_0(u_h)$ & DoF & $e_0(u_h)$ & DoF & $e_0(u_h)$ & DoF & $e_0(u_h)$ \\
\hline
1089 & 1.79e-01 & 1089 & 4.39e-02 & 625 & 7.72e-01 & 541 & 1.78e-02 \\
\hline
4225 & 1.05e-02 & 4225 & 1.09e-02 & 2401 & 1.20e-01 & 2169 & 4.86e-03 \\
\hline
16641 & 2.58e-03 & 16641 & 7.69e-04 & 9409 & 5.52e-03 & 9691 & 1.28e-03 \\
\hline
66049 & 6.67e-04 & 66049 & 9.03e-05 & 37249 & 1.66e-04 & 39641 & 3.54e-04 \\
\hline
263169 & 1.68e-04 & 263169 & 1.14e-05 & 148225 & 1.08e-05 & 155665 & 9.54e-05 \\
\hline
\end{tabular}
\caption{Relative $L^2$ errors of FEM and adaptive FEM in Case 2 of Example \ref{ex:Poisson}.}
\label{table:ex_Poisson_case2-fem}
\end{table}

To compare the performance of the AG-RaNN method with traditional numerical methods, we use the FEM and adaptive FEM to solve this example. All results are obtained using the $i$FEM package (\cite{Chen:2008ifem}). The results are presented in Table \ref{table:ex_Poisson_case2-fem}. We observe that $P_3$ FEM requires 148,225 degrees of freedom (DoF) to achieve a relative $L^2$ error of $1.08e-05$, and adaptive $P_1$ FEM requires 155,665 DoF to achieve a relative $L^2$ error of $9.54e-05$. However, using the AG-RaNN method, we achieve a relative $L^2$ error of $8.77e-07$ with only 3,453 DoF. Thus, the AG-RaNN method can obtain highly accurate numerical solutions with significantly fewer degrees of freedom.

\noindent\textbf{Case 3: $d=2$, Initialization + Neuron Growth + Layer Growth}

\begin{figure}[!htbp]
    \centering
    \includegraphics[width=0.98\textwidth]{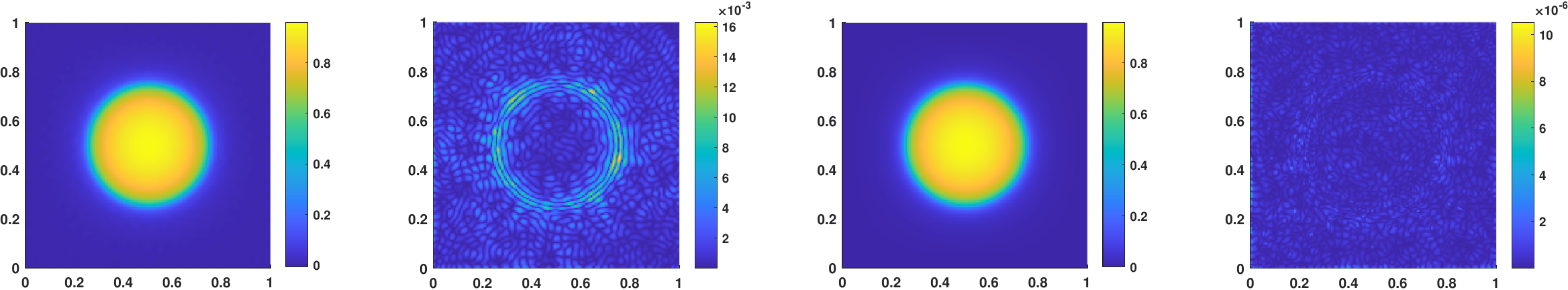}
    \caption{Numerical solutions and errors of $u_\rho^0$ (left) obtained by neuron growth and $u_\rho^4$ (right) obtained by layer growth in Case 3 of Example \ref{ex:Poisson}. }
    \label{fig:ex_Poisson_case3}
\end{figure}

This case uses frequency-based initialization and neuron growth to generate the first hidden layer, and add the second hidden layer using the layer growth strategy. We set $\bmm_\text{add}=(600,400,200,200,200,200,200)$, $m_2=1500$, $r_\text{max}=100$ and $\Lambda=100$, and generate $\bH_0$ randomly with $\br_2=(10,10)$. The relative $L^2$ errors and $\br^\text{opt}$ are presented in Table \ref{tab:ex_Poisson_case3} for neuron growth. The results obtained by layer growth are shown in Figure \ref{fig:ex_Poisson_case3} with $m_p=136$ and $e_0(u_\rho^4)=1.31e-06$. Thus, the AG-RaNN method can significantly improve accuracy using the layer growth strategy.

\begin{table}[!htbp]
    \centering
    \begin{tabular}{|c||c|c|c|c|c|c|c|c|c|c|}
    \hline
        $\br^\text{opt}$ & (11,14) & (42,27) & (19,54) & (61,32) & (24,73) & (77,37) & (43,70) \\ \hline
        $e_0(u_\rho)$ & 3.12e-02 & 1.66e-02 & 1.27e-02 & 1.06e-02 & 8.87e-03 & 7.28e-03 & 6.37e-03 \\ \hline
    \end{tabular}
    \caption{Results of the AG-RaNN with neuron growth in Case 3 of Example \ref{ex:Poisson}. }
    \label{tab:ex_Poisson_case3}
\end{table}

\noindent\textbf{Case 4: $d=3$, Layer Growth}

\begin{figure}[!htbp]
    \centering
    \includegraphics[width=\textwidth]{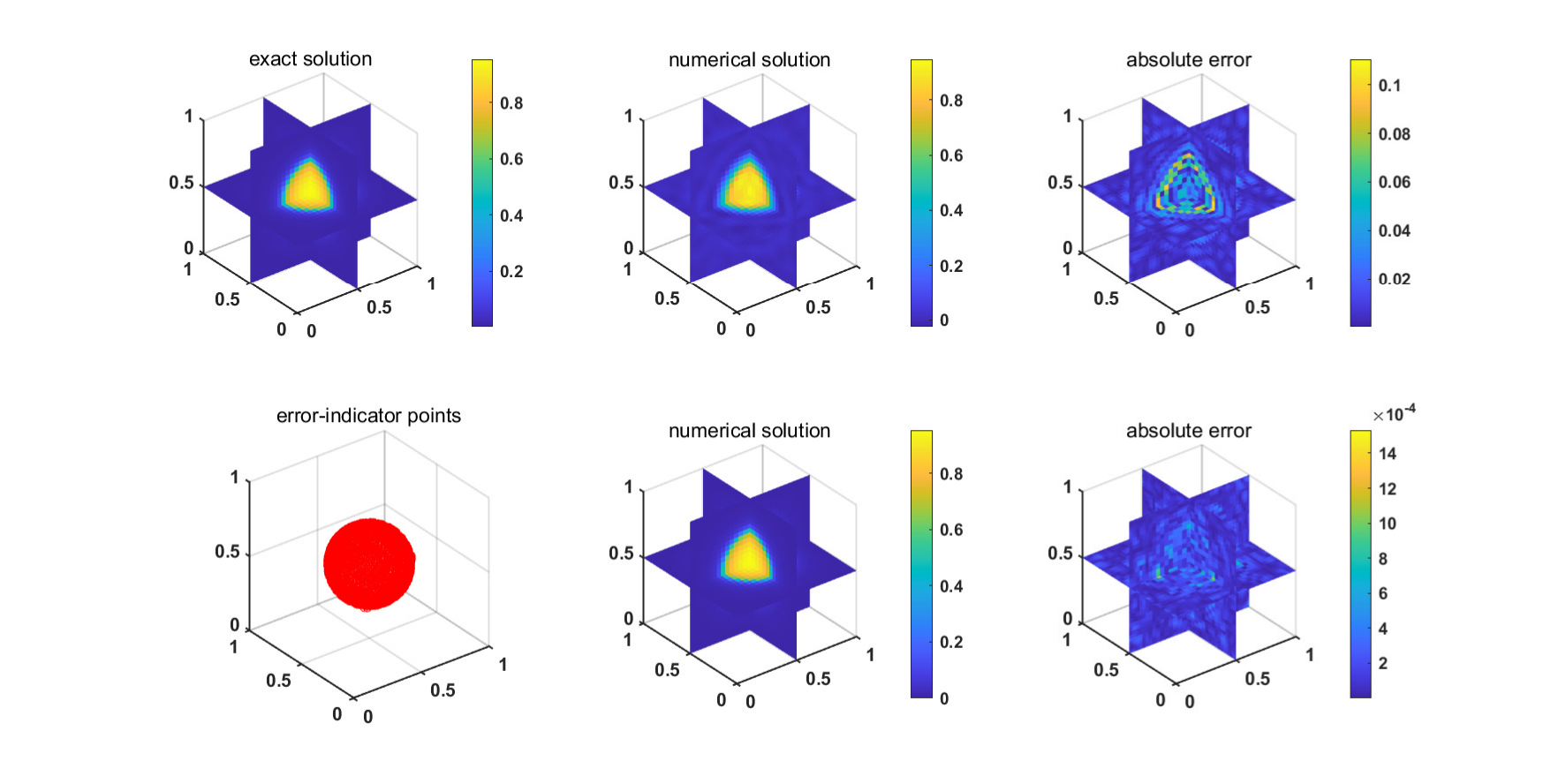}
    \caption{Exact solution, numerical solutions (middle), absolute errors (right), and error-indicator points $\bX_\text{err}$ for $m_1=2000$ and $m_2=2000$ in Case 4 of Example \ref{ex:Poisson}. }
    \label{fig:ex_Poisson_case4}
\end{figure}

We test a 3-D problem using layer growth, choosing $\br_1=(5,5,5)$ for the uniform distribution $\cU(-\br_1,\br_1)$, and generating $\bH_0$ randomly with $\br_2=(8,8,8)$. Table \ref{table:ex_Poisson_case4} shows the results for different values of $m_1$ and $m_2$. When $m_1=2000$ and $m_2=2000$, the numerical solutions and errors are shown in Figure \ref{fig:ex_Poisson_case4}. The AG-RaNN performs exceptionally well for this 3-D problem.

\begin{table}[!htbp]
\centering
\begin{tabular}{|c|c|c|c|c|c|c|c|}
\hline
 & $e_0(u_\rho^0)$ & $m_p$ & $m_2=500$ & $m_2=1000$ & $m_2=2000$ \\ \hline
$m_1=500$ & 2.07e-01 & 0 & 1.74e-02 & 1.11e-02 & 6.75e-03 \\\hline
$m_1=1000$ & 1.40e-01 & 0 & 8.07e-03 & 4.10e-03 & 2.49e-03 \\\hline
$m_1=2000$ & 1.01e-01 & 0 & 3.18e-03 & 1.67e-03 & 8.33e-04 \\\hline
\end{tabular}
\caption{Relative $L^2$ errors of AG-RaNN methods in Case 4 of Example \ref{ex:Poisson}. }
\label{table:ex_Poisson_case4}
\end{table}

\noindent\textbf{Case 5: $d=2$, Layer Growth}

\begin{figure}[!htbp]
    \centering
    \includegraphics[width=0.9\textwidth]{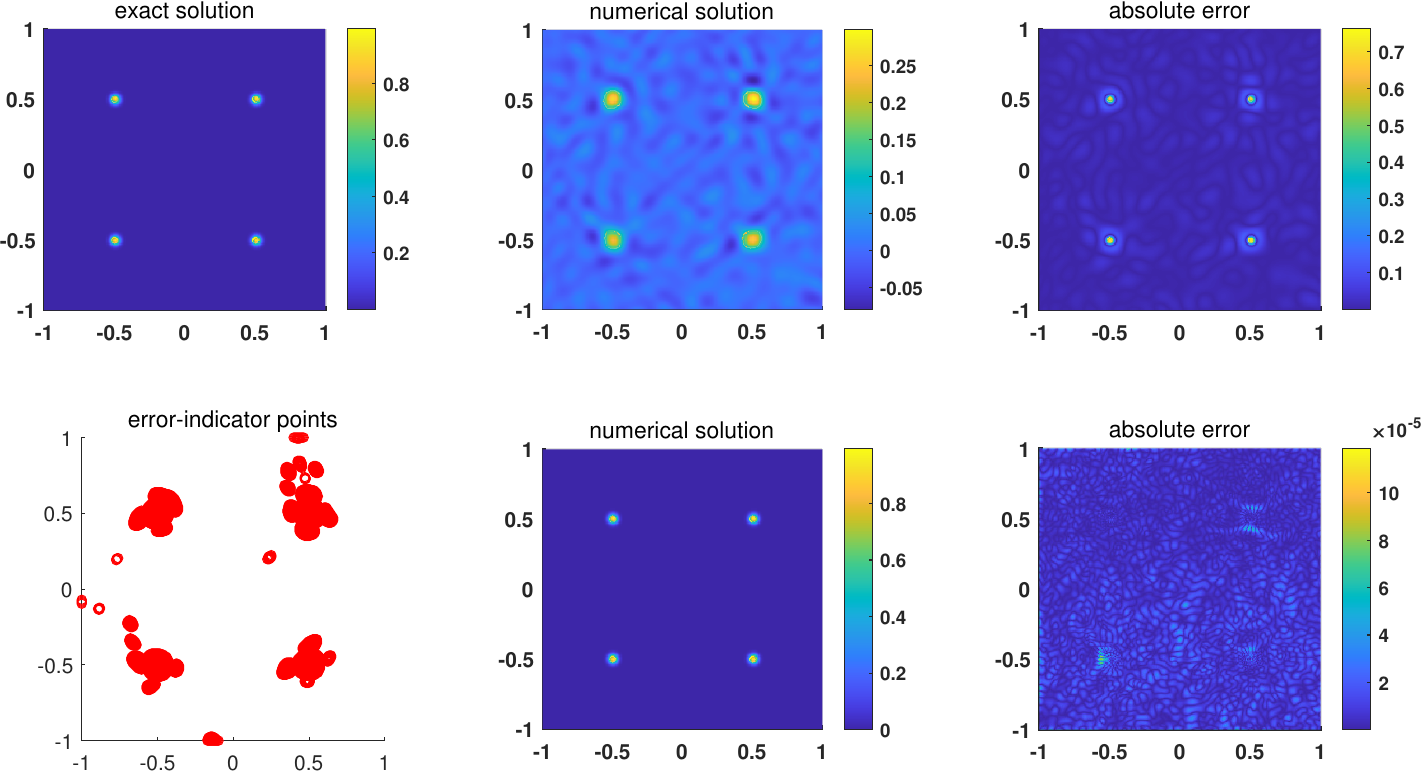}
    \caption{Exact solution, numerical solutions (middle), absolute errors (right) and error-indicator points $\bX_\text{err}$ obtained using the layer growth strategy as $P_1=4$ in Case 5 of Example \ref{ex:Poisson}. }
    \label{fig:ex_Poisson_case5_peak1}
\end{figure}

\begin{figure}[!htbp]
    \centering
    \includegraphics[width=0.9\textwidth]{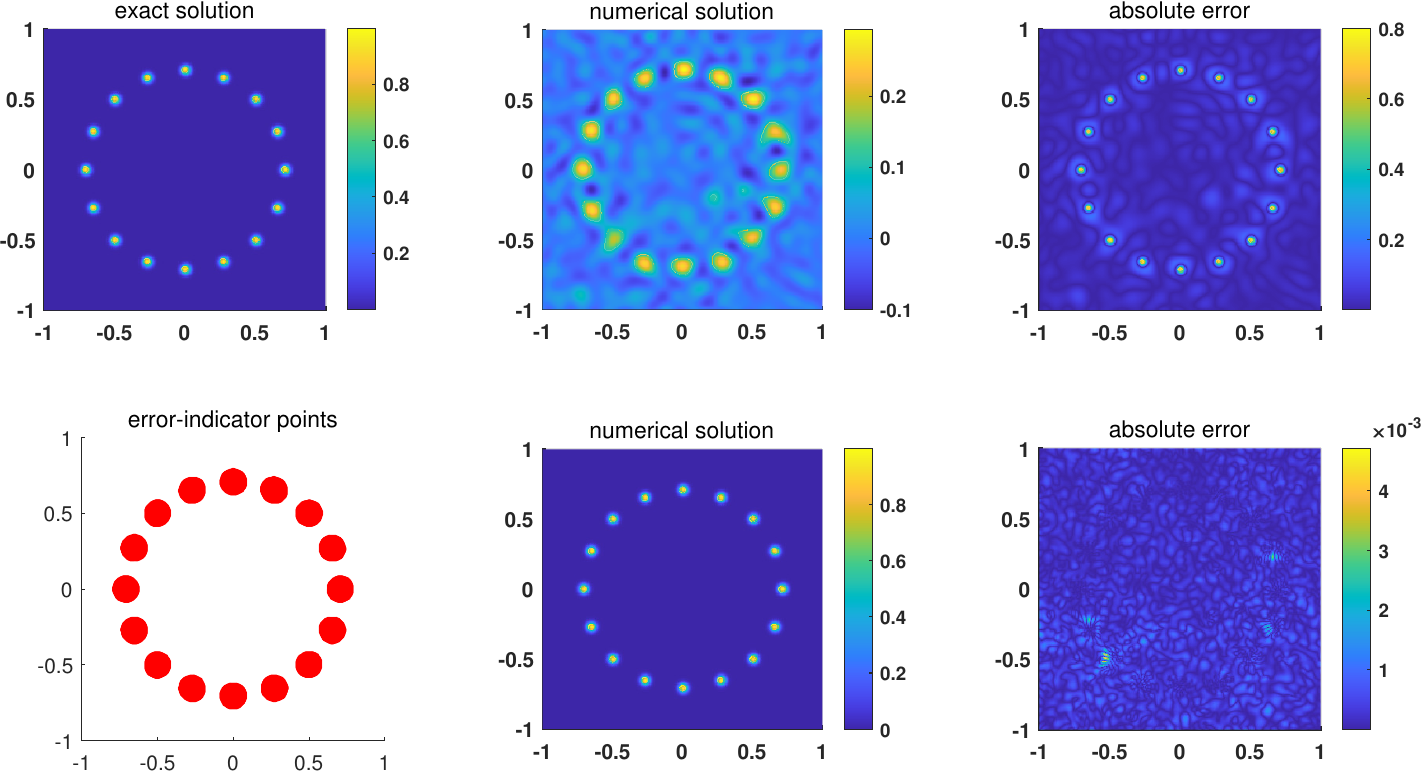}
    \caption{Exact solution, numerical solutions (middle), absolute errors (right) and error-indicator points $\bX_\text{err}$ obtained using the layer growth strategy as $P_1=16$ in Case 5 of Example \ref{ex:Poisson}. }
    \label{fig:ex_Poisson_case5_peak2}
\end{figure}

To verify the robustness of the layer growth strategy, we choose the same parameters to approximate the different exact solutions, with $m_1=500$, $m_2=3500$, $\br_1=(10,10)$ and $\br_2=(20,20)$. We consider three test cases: $P_1=4$, $P_1=16$ and $P_2=8$. The corresponding relative $L^2$ errors are shown in Table \ref{table:ex_Poisson_csae5_error}. The exact solutions, numerical solutions, absolute errors, and error-indicator points are shown in Figures \ref{fig:ex_Poisson_case5_peak1}--\ref{fig:ex_Poisson_case5_peak3}.

\begin{table}[!htbp]
\centering
\begin{tabular}{|c|c|c|c|c|c|c|c|}
\hline
 & $e_0(u_\rho^0)$ & $e_0(u_\rho^4)$ \\ \hline
$P_1=4$ & 7.6860e-01 & 2.2438e-04 \\ \hline
$P_1=16$ & 7.7185e-01 & 4.3955e-03 \\ \hline
$P_2=8$ & 7.0855e-01 & 3.2697e-02 \\ \hline
\end{tabular}
\caption{Relative $L^2$ errors of AG-RaNN methods in Case 5 of Example \ref{ex:Poisson}. }
\label{table:ex_Poisson_csae5_error}
\end{table}

\begin{figure}[!htbp]
    \centering
    \includegraphics[width=0.9\textwidth]{ex4_case3.eps}
    \caption{Exact solution, numerical solutions (middle), absolute errors (right) and error-indicator points $\bX_\text{err}$ obtained using the layer growth strategy as $P_2=8$ in Case 5 of Example \ref{ex:Poisson}. }
    \label{fig:ex_Poisson_case5_peak3}
\end{figure}

\begin{example}[Advection-Reaction Equation with Discontinuous Solution] \label{ex:AR}
We consider the following 2-D advection-reaction equation:
\begin{align}
    \nabla \cdot (\bbeta u)+u &\,=\,f\text{ in } D, \label{ex_AR_eq_pde}\\
    u &\,=\,g \text{ on } \Gamma, \label{ex_AR_eq_bdcond}
\end{align}
where $D = (0,1)^2$, $\bbeta=(0,1)^\top$ and $\Gamma = [0,1]\times \{0\}$. The exact solution is
\begin{align}
    u(x, y)\,=\,
    \begin{cases}
        1, & (x,y) \in D_1, \\
        0, & (x,y) \in D_2, 
    \end{cases}
\end{align}
where $D_1=\set{(x,y)\in D:x\in(0,0.1)\cup(0.2,0.3)\cup(0.4,0.5)\cup(0.6,0.7)\cup(0.8,0.9),y\in(0,1)}$ and $D_2=D\backslash D_1$.
\end{example}

The parameters are listed in Table \ref{tab:ex_AR_parameter}. 

\begin{table}[!htbp]
    \centering
    \begin{tabular}{|c|c|c|c|c|c|c|c|c|c|c|}
    \hline
    $\rho_1(x)$ & $N^I$ & $N^B$ & $m_1$ & $\eta$ & $\br_1$ \\ \hline
    $\tanh(x)$ & $100\times100$ & 200 & 1000 & 1 & $(20,1)$ \\ \hline
    \end{tabular}
    \caption{Parameters used in Example \ref{ex:AR}. }
    \label{tab:ex_AR_parameter}
\end{table}

\begin{figure}[!htbp]
    \centering
    \includegraphics[width=0.9\textwidth]{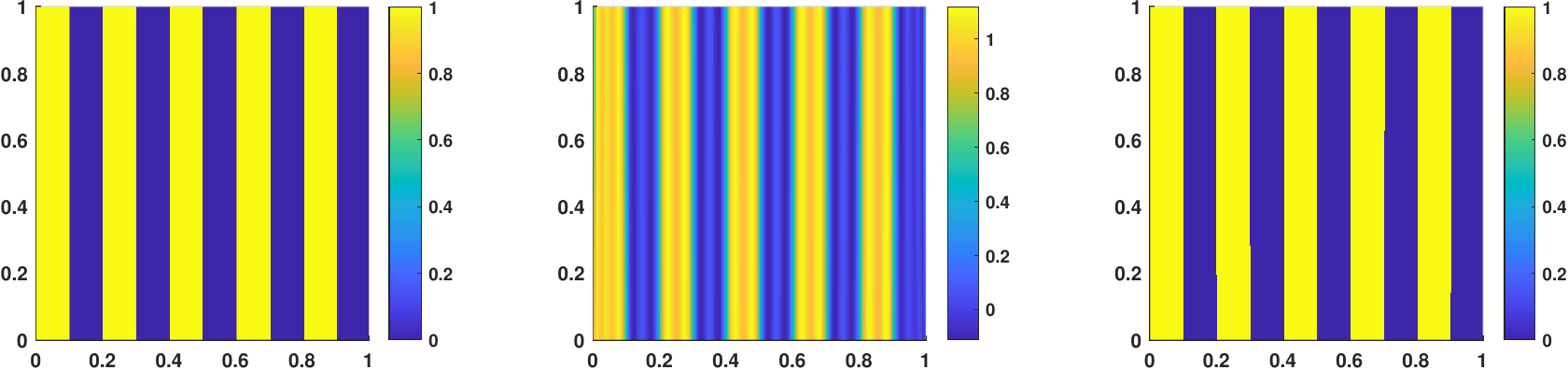}
    \caption{Exact solution (left) and numerical solutions (middle: $u_\rho^0$, right: domain splitting) in Example \ref{ex:AR}. }
    \label{fig:ex_AR_case2-domain}
\end{figure}

For more complex problems, identifying the error-indicator points $\bX_\text{err}$ for the layer growth strategy can be challenging. In such cases, the domain splitting strategy can be adopted. In this example, the solution has nine discontinuous lines (left subplot of Figure \ref{fig:ex_AR_case2-domain}), which would require at least 10 subdomains in traditional methods.

We first apply the RaNN method to obtain the approximate solution $u_\rho^0$. Based on the approximate solution $u_\rho^0$, we split the domain into two parts. For each subdomain, we use the same parameters as Table \ref{tab:ex_AR_parameter}. Using the domain splitting strategy, the numerical solutions are shown in Figure \ref{fig:ex_AR_case2-domain}. 

It is worth noting that we also attempted to use the layer growth strategy for this case. However, due to the numerous discontinuous lines, it is difficult to capture all of them through the error-indicator points $\bX_\text{err}$. 

\begin{example}[Burgers' Equation] \label{ex:Burgers}
Consider the following Burgers' equation 
\begin{align} 
    u_t+h(u)_x-\varepsilon_B u_{xx} &\,=\,f \text{ in }(0,T)\times D, \\
    u &\,=\,g \text{ on } (0,T)\times\Gamma,\\
    u(0,x) &\,=\,u_0 \text{ in } D,
\end{align}
where $h(u) = u^2/2$ and $T = 1$. We explore four cases using different strategies for each.

\noindent\textbf{Case 1: }$\varepsilon_B =0.1/\pi$, $f=0$, $g=0$, $D = (-1,1)$ and $\Gamma=\partial D$. The initial condition is
\begin{align*}
    u_0(x)\,=\,-\sin(\pi x).
\end{align*}
Using the Cole-Hopf transformation, the exact solution is:
\begin{align*}
    u(t,x)\,=\,\frac{\int_{-\infty}^{\infty}\frac{x-y}{t}\,e^{-\frac{H(x,y)}{2\varepsilon_B}}\,\dy}{\int_{-\infty}^{\infty}\,e^{-\frac{H(x,y)}{2\varepsilon_B}}\,\dy}\quad \text{and}\quad H(x,y)\,=\,\int_0^y\,u_0(s)\,\ds+\frac{(x-y)^2}{2t}.
\end{align*}

\noindent\textbf{Case 2: }$\varepsilon_B = 0.01$, $D = (0,1)$ and $\Gamma=\partial D$. The exact solution is
\begin{align}
    u(t,x)\,=\,\frac{1}{1+e^{\frac{x-t}{2\varepsilon_B}}}.
\end{align}

\noindent\textbf{Case 3: }$\varepsilon_B = 0$, $f=0$, $D = (-1,2)$ and $\Gamma=\Gamma^L\cup\Gamma^R=\set{(t,-1):t\in[0,T]}\cup\set{(t,2):t\in[0,T]}$. It deals with a rarefaction wave with the boundary and initial conditions
\begin{align}
    g(t,x)\,=\,
    \begin{cases}
        0 & \text{ on }\Gamma^L,\\
        1 & \text{ on }\Gamma^R,
    \end{cases} \quad \text{and}\quad    
    u_0(x)\,=\,
    \begin{cases}
        0 & x<0, \\
        1 & x>0.
    \end{cases}
\end{align}

\noindent\textbf{Case 4: }$\varepsilon_B = 0$, $f=0$, $D = (-1,2)$ and $\Gamma=\Gamma^L\cup\Gamma^R=\set{(t,-1):t\in[0,T]}\cup\set{(t,2):t\in[0,T]}$. It deals with a shock wave with the boundary and initial conditions
\begin{align}
    g(t,x) = 
    \begin{cases}
        1 & \text{ on }\Gamma^L,\\
        0 & \text{ on }\Gamma^R,
    \end{cases} \quad \text{and}\quad    
    u_0(x)=
    \begin{cases}
        1 & x<0, \\
        0 & x>0.
    \end{cases}
\end{align}

\end{example}

Parameters for all cases are presented in Table \ref{tab:ex_Burgers_parameter}.
\begin{table}[!htbp]
    \centering
    \begin{tabular}{|c|c|c|c|c|c|c|c|c|c|}
    \hline
     & $\rho_1(x)$ & $\rho_2(x)$ & $m_1$ & $N^I$ & $N^B$ & $\eta$ & $\br_1$ \\ \hline
    \textbf{Case 1} & $G(x)$ & N/A & 1000 & $100\times200$ & 600 & 200 & N/A \\ \hline
    \textbf{Case 2} & $\tanh(x)$ & $G(x)$ & N/A & $200\times200$ & 600 & 200 & (20,20) \\ \hline
    \textbf{Case 3} & $\tanh(x)$ & $G(x)$ & 1000 & $66\times200$ & 332 & 10 & (3,15) \\ \hline
    \textbf{Case 4} & $\tanh(x)$ & N/A & 1000 & $66\times200$ & 332 & 20 & (15,15) \\ \hline
    \end{tabular}
    \caption{Parameters used in Example \ref{ex:Burgers}. }
    \label{tab:ex_Burgers_parameter}
\end{table}

\noindent\textbf{Case 1: Initialization + Neuron Growth}

\begin{figure}[!htbp]
    \centering
    \includegraphics[width=0.8\textwidth]{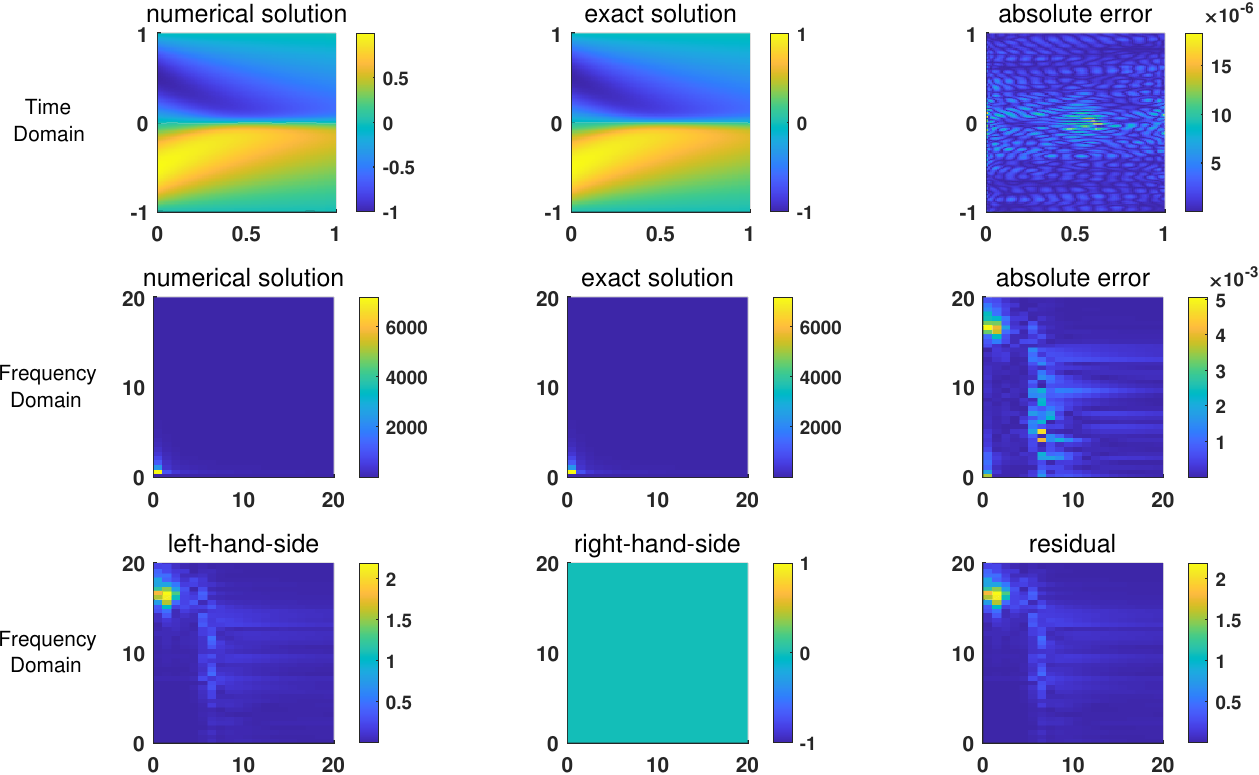}
    \caption{Best result obtained by the manually tuned RaNN in Case 1 of Example \ref{ex:Burgers}.}
    \label{fig:ex_Burgers_case1_init-intra-Burgers_best}
\end{figure}

For this time-dependent problem, we employ space-time RaNNs, treating the time variable as equivalent to a spatial variable. To verify that the frequency-based parameter initialization and neuron growth strategies can achieve results comparable to the best manually tuned RaNN with $m_1 = 1000$, we identify the near-optimal uniform distribution parameter is $\br=(8,25)$ and iteration counts $(It_P,It_N)=(20,0)$. The relative $L^2$ error is $4.84e-06$, as shown in Figure \ref{fig:ex_Burgers_case1_init-intra-Burgers_best}. 

\begin{figure}[!htbp]
    \centering
    \includegraphics[width=0.98\textwidth]{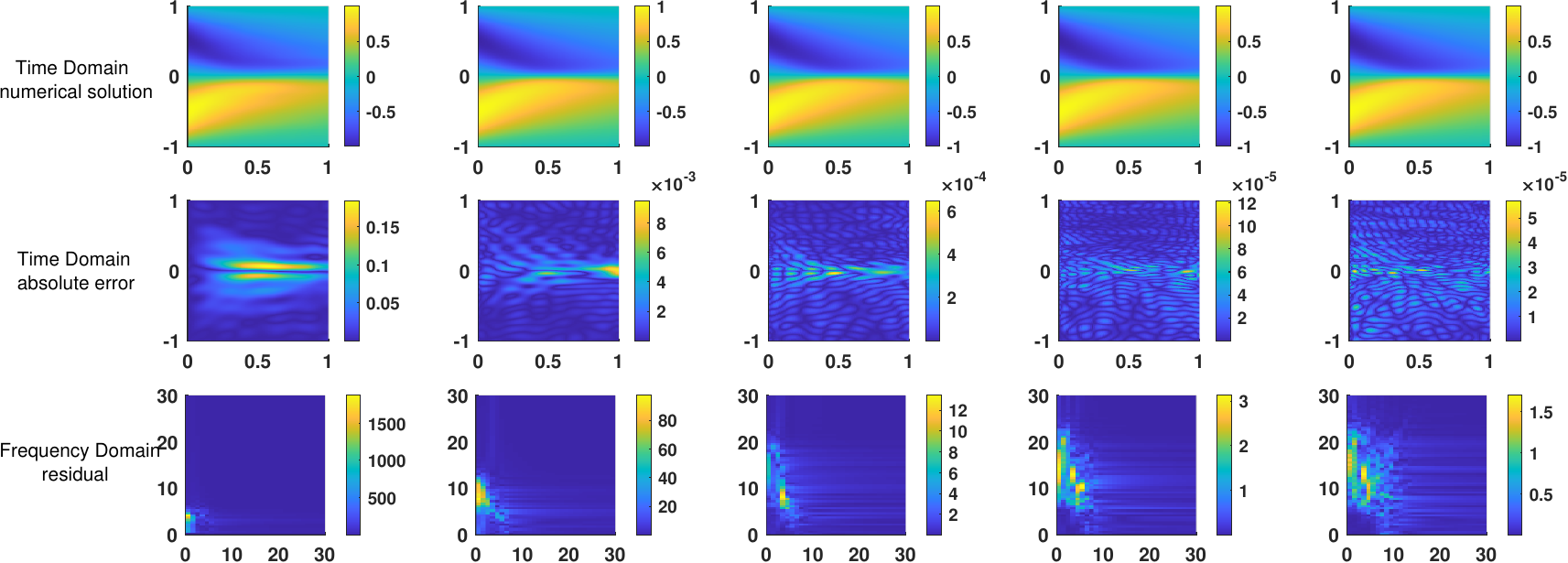}
    \caption{Results obtained by the AG-RaNN method in Case 1 of Example \ref{ex:Burgers}. The first column shows the result obtained by initialization, while the other columns show results obtained by the neuron growth strategy.}
    \label{fig:ex_Burgers_case1_init-intra-Burgers}
\end{figure}

For the adaptive-growth RaNN method, we set $r_\text{max}=100$, $\Lambda=100$ and $\bmm_\text{add}=(200,200,200,200,200)$. Since the right-hand side is zero, providing no useful information, as stated in Remark \ref{re:zeroright}, we choose $\br^\text{opt}=(1,1)$ as the initial distribution before neuron growth, with iteration counts $(It_P,It_N)=(10,0)$. For each growth stage, we use the final approximate solution as the initial solution for the next growth stage, with iteration counts $(It_P,It_N)=(5,0)$. The relative $L^2$ errors and $\br^\text{opt}$ are presented in Table \ref{tab:ex_Burgers_case1_init-intra-Burgers}, with the results shown in Figure \ref{fig:ex_Burgers_case1_init-intra-Burgers}. The table shows that the AG-RaNN with initialization and neuron growth achieves a relative $L^2$ error of $2.83e-05$.

\begin{table}[!htbp]
    \centering
    \begin{tabular}{|c||c|c|c|c|c|c|c|c|c|c|}
    \hline
        $\br^\text{opt}$ & (1,1) & (3,11) & (3,30) & (12,40) & (16,33) \\ \hline
        $e_0(u_\rho)$ & 1.64e-01 & 1.03e-02 & 5.47e-04 & 1.52e-04 & 2.83e-05 \\ \hline
    \end{tabular}
    \caption{Results of the AG-RaNN with neuron growth in Case1 of Example \ref{ex:Burgers}.}
    \label{tab:ex_Burgers_case1_init-intra-Burgers}
\end{table}

\noindent\textbf{Case 2: Layer Growth}

\begin{figure}
    \centering
    \includegraphics[width=0.9\linewidth]{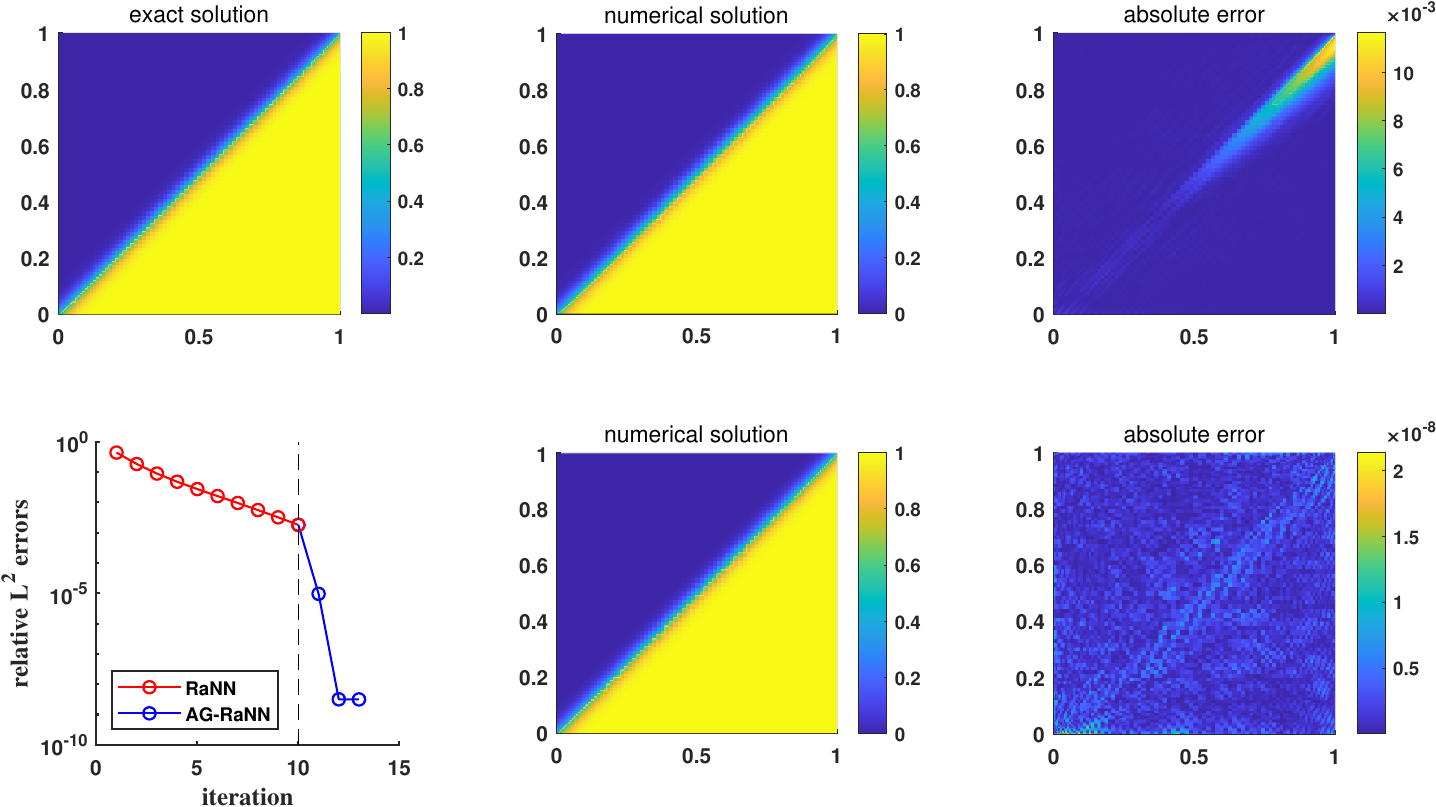}
    \caption{Exact solution, numerical solutions (middle) and absolute errors (right), and relative $L^2$ errors for $u_\rho^0$ and $u_\rho^4$ in Case 2 of Example \ref{ex:Burgers}.}
    \label{fig:ex_Burgers_case2_layer}
\end{figure}

We use $(It_P,It_N)=(10,0)$ for the first training and $(It_P,It_N)=(0,3)$ during layer growth. Comparing $u_\rho^3$ and $u_\rho^4$ (as defined in Case 2 of Example \ref{ex:Poisson}), $\bH_0$ is generated randomly with parameters $\br_2=(7,7)$. The results are presented in Table \ref{table:ex_Burgers_case2-rnn}. The numerical solutions and errors for $u_\rho^0$ and $u_\rho^4$ with $m_1=2000$, $m_2=2000$ are shown in Figure \ref{fig:ex_Burgers_case2_layer}. We observe that the AG-RaNN method can achieve high precision for this problem, which has a sharp solution.

\begin{table}[!htbp]
\centering
\begin{tabular}{|c|c|c|c|c|c|c|c|}
\hline
$m_1$ & $e_0(u_\rho^0)$ & $m_p$ & $m_2$ & $e_0(u_\rho^3)$ & $e_0(u_\rho^4)$ \\ \hline\hline
\multirow{4}{*}{500} & \multirow{4}{*}{2.9285e-02} & \multirow{4}{*}{0} & 500 & 1.7374e-03 & 3.5720e-05  \\ \cline{4-6} 
    &  &  & 1000 & 1.5419e-03 & 3.5981e-06 \\ \cline{4-6} 
    &  &  & 1500 & 1.1241e-03 & 1.9748e-06 \\ \cline{4-6} 
    &  &  & 2000 & 9.2282e-04 & 1.1140e-06 \\ \hline\hline
\multirow{4}{*}{1000} & \multirow{4}{*}{7.5540e-03} & \multirow{4}{*}{0} & 500 & 5.3504e-04 & 9.8944e-06  \\ \cline{4-6} 
    &  &  & 1000 & 4.7899e-04 & 1.3535e-06 \\ \cline{4-6} 
    &  &  & 1500 & 4.2683e-04 & 7.9742e-07 \\ \cline{4-6} 
    &  &  & 2000 & 3.6866e-04 & 5.6555e-07 \\ \hline\hline
\multirow{4}{*}{1500} & \multirow{4}{*}{1.8014e-03} & \multirow{4}{*}{0} & 500 & 1.1080e-05 & 1.4255e-07 \\ \cline{4-6} 
    &  &  & 1000 & 8.8926e-06 & 3.6183e-08 \\ \cline{4-6} 
    &  &  & 1500 & 8.7247e-06 & 2.3615e-08 \\ \cline{4-6} 
    &  &  & 2000 & 8.9395e-06 & 1.7544e-08 \\ \hline\hline
\multirow{4}{*}{2000} & \multirow{4}{*}{1.8058e-03} & \multirow{4}{*}{1} & 500 & 7.1958e-06 & 4.7314e-08 \\ \cline{4-6} 
    &  &  & 1000 & 4.2425e-06 & 7.7380e-09 \\ \cline{4-6} 
    &  &  & 1500 & 4.9113e-06 & 4.6066e-09 \\ \cline{4-6} 
    &  &  & 2000 & 5.3040e-06 & 3.1493e-09 \\ \hline
\end{tabular}
\caption{Relative $L^2$ errors of AG-RaNN methods in Case 2 of Example \ref{ex:Burgers}. }
\label{table:ex_Burgers_case2-rnn}
\end{table}

\noindent\textbf{Case 3: Layer Growth}

\begin{figure}[!htbp]
    \centering
    \includegraphics[width=0.9\textwidth]{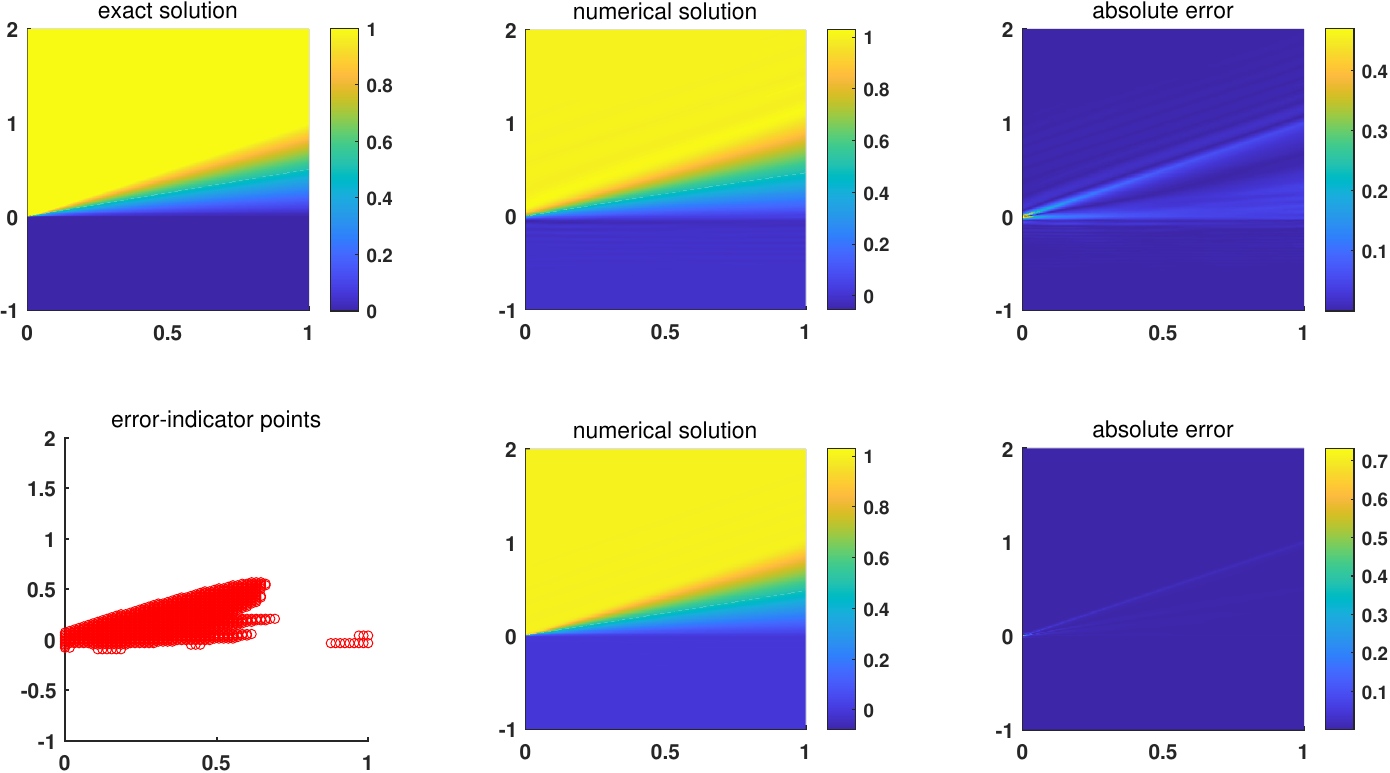}
    \caption{Exact solution, numerical solutions (middle), absolute errors (right) and error-indicator points $\bX_\text{err}$ obtained using the layer growth strategy in Case 3 of Example \ref{ex:Burgers}. }
    \label{fig:ex_Burgers_case3_rarefaction-layer}
\end{figure}

We use $\|\nabla u_\rho^0\|$ as the error indicator, which is a typical gradient-based shock detector. Other detectors, such as entropy-based shock detectors, can also be considered. For the initial solution $ u_\rho^0 $, we set $ (It_P, It_N) = (10, 0) $, and for the final approximation, $ (It_P, It_N) = (0, 5) $. We take $m_2=1000$ and $\br_2=(20,20)$ for layer growth. The numerical solutions and errors are shown in Figure \ref{fig:ex_Burgers_case3_rarefaction-layer}, with relative $L^2$ errors of $2.83e-02$ for $u_\rho^0$ and $4.31e-03$ for $u_\rho^4$.

\noindent\textbf{Case 4: All Strategies}

\begin{figure}[!htbp]
    \centering
    \includegraphics[width=0.9\textwidth]{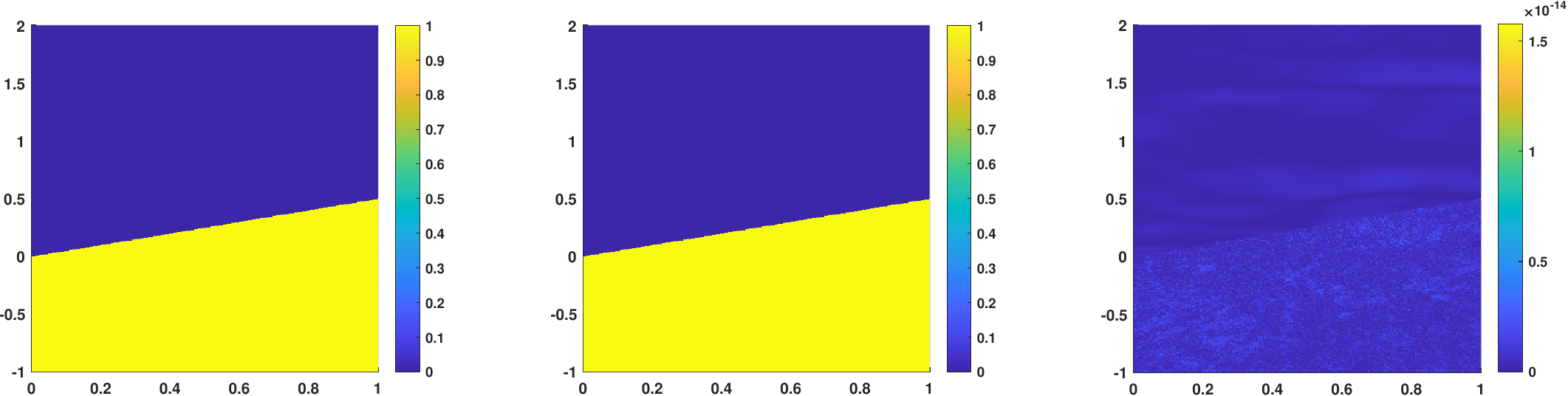}
    \caption{Numerical solution and absolute error using characteristic line domain splitting in Case 4 of Example \ref{ex:Burgers}. }
    \label{fig:ex_Burgers_case4-known}
\end{figure}

In this case, the domain is decomposed based on characteristic lines and boundary conditions. The parameters are shown in Table \ref{tab:ex_Burgers_parameter}. The numerical solution and absolute error are presented in Figure \ref{fig:ex_Burgers_case4-known}, with a relative $L^2$ error of $9.78e-16$.

Finally, we combine all proposed strategies to solve this example. To make neuron growth feasible, we use the $L^1$ norm as a stopping criterion:
\begin{align} \label{eq:l1}
    \left\|u_\rho^\text{new}-u_\rho^\text{old}\right\|_{L^1((0,T)\times D)}<10^{-4},
\end{align}
where $u_\rho^\text{new}$ and $u_\rho^\text{old}$ are the solutions from consecutive growth stages.

\begin{figure}[!htbp]
    \centering
    \hspace{-1.0cm}
    \includegraphics[width=1.05\textwidth]{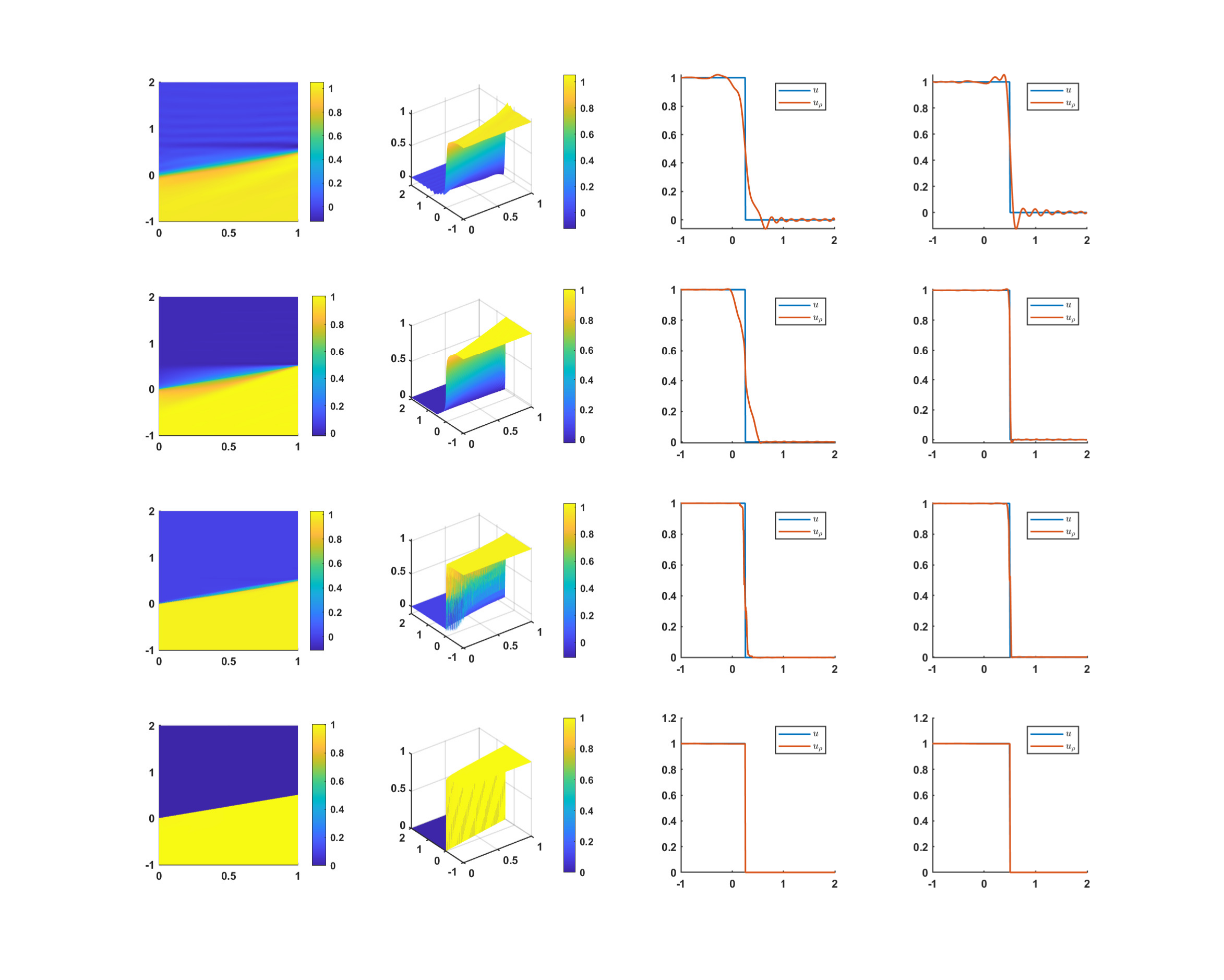}
    \caption{Results obtained by AG-RaNN method in Case 4 of Example \ref{ex:Burgers}. The first two columns show the numerical solutions obtained by different strategies, while the subsequent columns show the results at $t=0.5$ (third column) and $t=1$ (last column). }
    \label{fig:ex_Burgers_case4-all}
\end{figure}

Initially, we employ the frequency-based parameter initialization strategy to generate a network with 100 neurons. The parameters are set as $ N^I = 300 \times 100 $, $ N^B = 666 $, $ \eta = 20 $, $ \br_1 = (1, 1) $ and $ (It_P, It_N) = (10, 0) $. The activation function is $ G(x) $. Using the neuron growth strategy with $ \Lambda = 100 $, $ r_{\text{max}} = 100 $ and $ (It_P, It_N) = (5, 0) $, the network grows twice, adding 100 neurons each time based on the stopping criterion \eqref{eq:l1}. The relative $ L^2 $ errors and optimal hyperparameters $ \br^{\text{opt}} $ are listed in Table \ref{tab:ex_Burgers_case4_init-intra}, with the results displayed in the first row of Figure \ref{fig:ex_Burgers_case4-all}. We observe that the Gibbs phenomenon arises near discontinuities, making it difficult for the single layer RaNN to accurately capture the shock wave.

\begin{table}[!htbp]
    \centering
    \begin{tabular}{|c||c|c|c|c|c|c|c|c|c|c|}
    \hline
        $\br^\text{opt}$ & (1,1) & (3,8) & (3,8) \\ \hline
        $e_0(u_\rho)$ & 1.8389e-01 & 1.3625e-01 & 1.2972e-01 \\ \hline
    \end{tabular}
    \caption{Results of the AG-RaNN with neuron growth in Case 4 of Example \ref{ex:Burgers}.}
    \label{tab:ex_Burgers_case4_init-intra}
\end{table}

Next, we apply the layer growth strategy to add a new hidden layer to the previously grown network. For this, we select $ m_2 = 700 $, $ \br_2 = (10, 10) $, $ (It_P, It_N) = (0, 5) $, and use $ \rho_2(x) = G(x) $ for the layer growth. The error indicator remains $ \|\nabla u_\rho^0\| $. The numerical solution is shown in the second row of Figure \ref{fig:ex_Burgers_case4-all}, with a relative $ L^2 $ error of  $1.24e-01$. Due to the presence of a shock wave, the solution is discontinuous along a specific characteristic line, which limits the effectiveness of the layer growth strategy. Let us denote this numerical solution as $ \tilde{u}_\rho $. 

Therefore, we use a discontinuous activation function for layer growth, defined as:
\begin{align*}
    \rho_2(x)\,=\,\begin{cases}
        0, &x\le0,\\x+1, &x>0.
    \end{cases}    
\end{align*}
Using the same parameters as before, the second hidden layer with the above discontinuous activation function is added following neuron growth. The numerical results are displayed in the third row of Figure \ref{fig:ex_Burgers_case4-all}, with a relative $ L^2 $ error of $8.36e-02$, showing a significant improvement in performance.

When all hidden layers use continuous activation functions, approximating a discontinuous solution becomes challenging. In such cases, domain splitting proves more effective. We set $ R = 2 $, $ (It_P, It_N) = (2, 0) $, $ \br = (2, 2) $, $ m_1 = 500 $, and $ \rho(x) = G(x) $ for each subdomain. To aid domain splitting, we use the numerical solution $ \tilde{u}_\rho $ to split the domain. The results are shown in the last row of Figure \ref{fig:ex_Burgers_case4-all}, achieving a relative $ L^2 $ error of $ 3.91 \times 10^{-2} $. This approach yields a noticeable improvement in handling discontinuities, for example, when $ t = 0.5 $ and $ t = 1 $, the method captures the shock wave with high precision.

\begin{example}[Allen-Cahn Equation] \label{ex:allencahn}
We consider the following Allen-Cahn equation
\begin{align} 
    u_t-\varepsilon_a u_{xx}+\kappa(u^3-u)&=0\text{ in }(0,T)\times (-1,1),\\
    u(t,-1)&=u(t,1) \text{ on }(0,T),\\
    u_x(t,-1)&=u_x(t,1) \text{ on }(0,T),\\
    u(0,x)&=u_0(x) \text{ in }(-1,1),
\end{align}
where $\varepsilon_a=10^{-4}$, $\kappa=5$ and $T=1$. The initial condition is given by
\begin{align*}
    u_0(x)=x^2\cos(\pi x).
\end{align*}
\end{example}

In this example, we investigate the influence of random seeds on the AG-RaNN method. For different random
seeds, we use the same set of hyperparameters, summarized in Table~\ref{tab:ex_AllenCahn_parameter}.
\begin{table}[!htbp]
    \centering
    \begin{tabular}{|c|c|c|c|c|c|c|c|c|c|}
    \hline
    $\rho_1(x)$ & $\rho_2(x)$ & $m_1$ & $m_2$ & $N^I$ & $N^B$ & $\eta$ & $\br_1$ & $\br_2$ \\ \hline
    $G(x)$ & $G(x)$ & 3000 & 1500 & $150\times300$ & 2000 & 10 & (5,70) & (2,2)\\ \hline
    \end{tabular}
    \caption{Parameters used in Example \ref{ex:allencahn}. }
    \label{tab:ex_AllenCahn_parameter}
\end{table}

Since an exact solution is not available for this problem, we adopt as reference the high-fidelity solution provided in~\cite{Urban2025PINN}. We first apply the RaNN method with $(It_P,It_N) = (0,15)$ to obtain an initial approximation $u_\rho^0$. Then we perform layer growth with $(It_P,It_N) = (0,5)$ to obtain the final approximation $u_\rho^4$.

\begin{figure}[!htbp]
    \centering
    \includegraphics[width=0.8\textwidth]{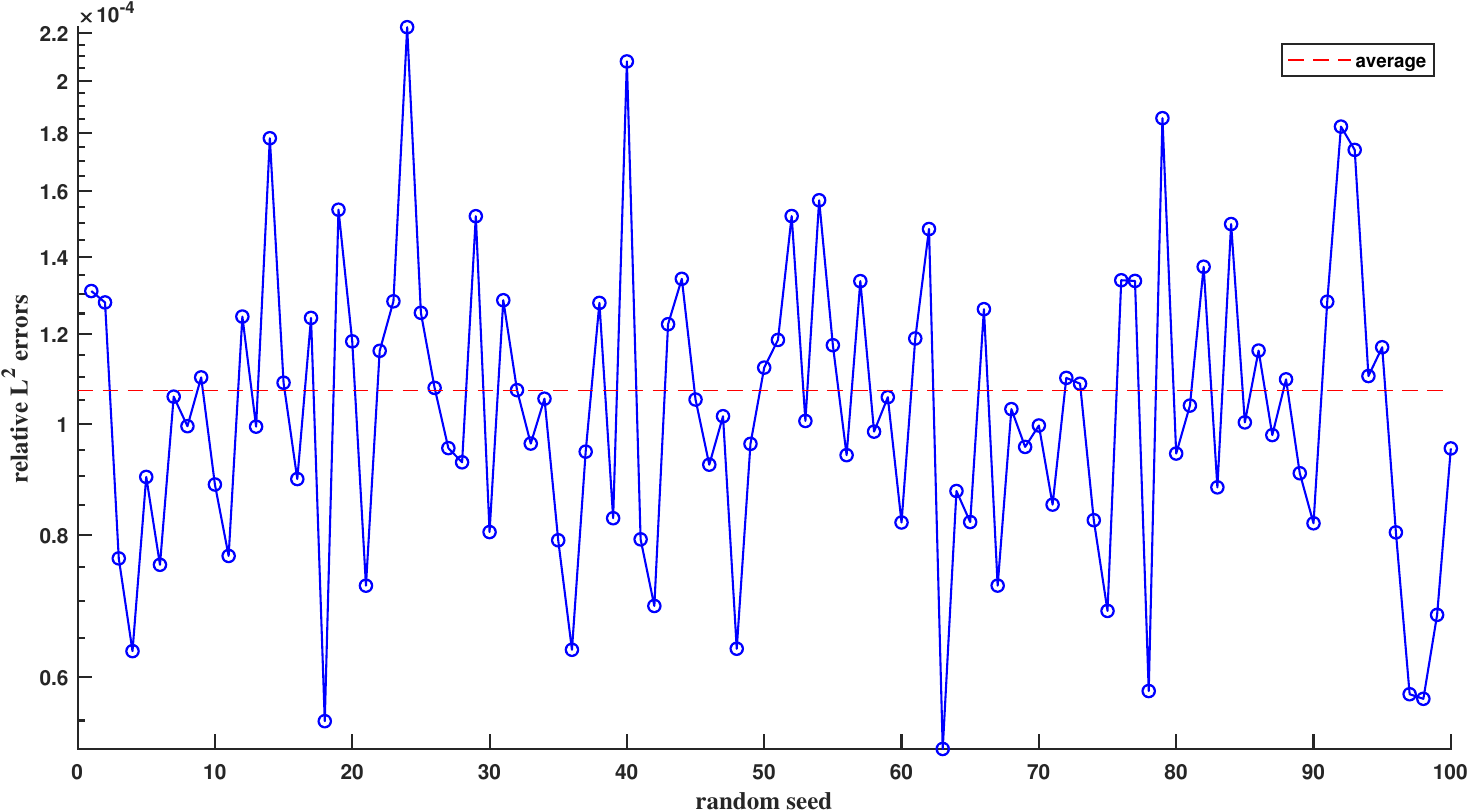}
    \caption{Relative $L^2$ errors $e_0(u_\rho^4)$ obtained by the AG-RaNN method with 100 different random seeds in Example~\ref{ex:allencahn}. The horizontal line indicates the error of the ensemble-averaged solution.}
    \label{fig:ex_allencahn_average}
\end{figure}

\begin{figure}[!htbp]
    \centering
    \includegraphics[width=1\textwidth]{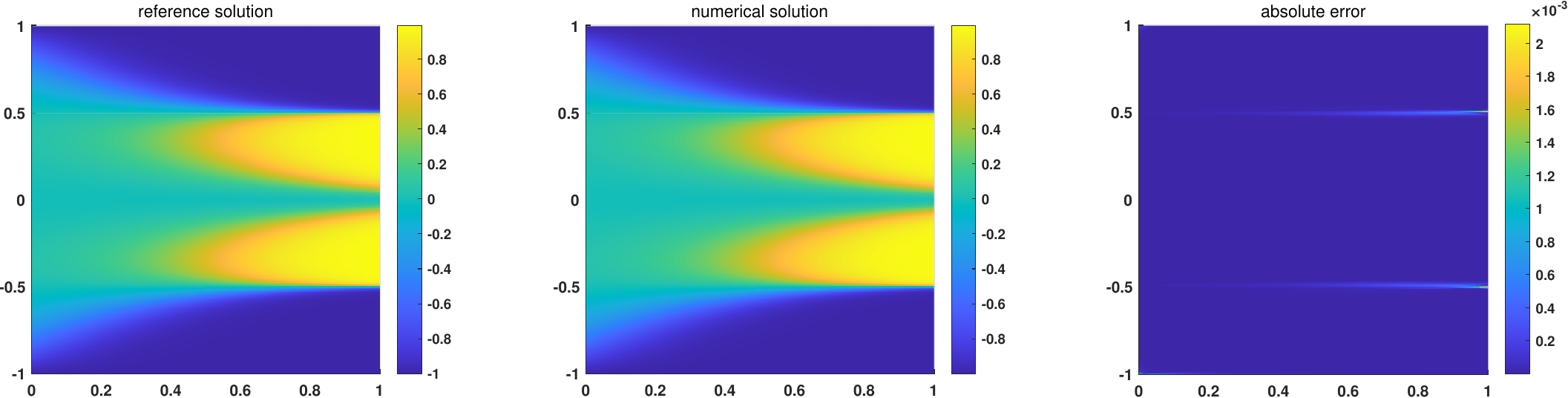}
    \caption{Reference solution (left), ensemble-averaged AG-RaNN solution (middle), and pointwise absolute error (right) in Example~\ref{ex:allencahn}.}
    \label{fig:ex_allencahn_solution}
\end{figure}

To assess the sensitivity with respect to random initialization, we set \texttt{rng($i$)} for $i = 1,\dots,100$ and compute the corresponding relative $L^2$ errors $e_0(u_\rho^4)$, which are reported in Figure~\ref{fig:ex_allencahn_average}. Following the ensemble-averaging strategy used in \cite{Breiman1996}, we further average all numerical solutions to obtain a mean solution. This ensemble-averaged AG-RaNN solution attains a relative $L^2$ error of $6.14\times 10^{-5}$ and is shown in Figure~\ref{fig:ex_allencahn_solution}, together with the reference solution and the pointwise absolute error.

\section{Summary}
\label{sec:sum}

The Adaptive-Growth Randomized Neural Network framework offers a novel, genuinely constructive approach to solving partial differential equations by dynamically expanding the network's capacity through neuron and layer growth strategies. This method is particularly effective for handling PDEs with sharp or discontinuous solutions, where traditional neural networks or numerical methods may struggle.

AG-RaNN employs four key strategies: frequency-based parameter initialization, neuron growth, layer growth, and domain splitting. The frequency-based parameter initialization optimizes initial weight distribution, ensuring the network is well-prepared to capture the solution's significant frequency components. Neuron growth progressively adds neurons to existing hidden layers, refining the solution by addressing residual errors. When a single hidden layer is insufficient, layer growth introduces additional layers to capture more complex patterns, often using local basis functions. For discontinuous solutions, domain splitting strategy divides the domain into subdomains, each handled by separate networks, effectively addressing discontinuities.

Within an abstract Hilbert-space setting, we derive an error decomposition that separates approximation, statistical, and optimization errors, and we prove convergence of the RaNN method under mild assumptions. Through various numerical examples, the AG-RaNN has demonstrated superior accuracy compared to traditional methods like finite element method, especially for complex problems with sharp gradients or discontinuities. By dynamically adapting its structure, AG-RaNN achieves highly accurate solutions with fewer DoF.

Future research will proceed along several directions. First, the present analysis focuses on single-hidden-layer RaNNs. The adaptive layer growth in AG-RaNN only enlarges the hypothesis space, so our approximation error bounds remain valid in the sense that deeper architectures can achieve at least the same accuracy as the shallow model.
However, the current results establish only convergence in a qualitative sense; a rigorous study of the convergence order and computational complexity of the fully adaptive procedure, in particular for genuinely multi-layer randomized networks, is left for future work. Beyond these theoretical questions, it is important to further enhance AG-RaNN in practice by extending its application to higher-dimensional and multi-physics problems, developing adaptive hyperparameter selection and optimization strategies, and refining the domain-splitting techniques for a more robust treatment of discontinuities and moving interfaces. In addition, extending the present ideas to other scientific computing tasks---such as inverse problems, operator learning, and data-driven model reduction---would be a valuable direction for future research.

\end{document}